\documentclass[11pt]{article}
\usepackage[margin=1in]{geometry}
\usepackage{amsmath,amssymb,amsthm}
\usepackage{booktabs,array}
\usepackage[numbers]{natbib}
\usepackage[T1]{fontenc}
\usepackage{lmodern,microtype}
\usepackage[hidelinks]{hyperref}
\hypersetup{
	colorlinks=false
}
\usepackage[nameinlink,noabbrev]{cleveref}

\newcommand{\E}{\mathbb E}
\newcommand{\Pp}{\mathbb P}
\newcommand{\R}{\mathbb R}
\newcommand{\N}{\mathbb N}
\newcommand{\ind}{\mathbf 1}
\newcommand{\cL}{\mathcal L}
\newcommand{\norm}[1]{\lVert #1\rVert}
\newcommand{\Law}{\operatorname{Law}}
\newcommand{\dd}{\mathrm{d}}
\newtheorem{theorem}{Theorem}
\newtheorem{lemma}{Lemma}
\newtheorem{proposition}{Proposition}
\newtheorem{corollary}{Corollary}
\theoremstyle{definition}
\newtheorem{assumption}{Assumption}
\newtheorem{example}{Example}
\newtheorem{remark}{Remark}
\crefname{theorem}{Theorem}{Theorems}
\crefname{assumption}{Assumption}{Assumptions}

\title{Routing Information and Stationary Limits\\in Many-Server Load Balancing}
\author{Wei You\thanks{Department of Industrial Engineering and Decision Analytics, The Hong Kong University of Science and Technology. Email: \href{mailto:weiyou@ust.hk}{\texttt{weiyou@ust.hk}}.}}
\date{\today}

\begin{document}
\maketitle
\begin{abstract}
We study how routing information determines stationary congestion in many-server systems with Poisson arrivals and homogeneous exponential servers.
We consider policies that prioritize idle servers and select a queue stochastically no longer than a uniformly sampled queue.
In nondegenerate slowdown, where spare capacity is fixed, we identify population laws from the fraction of queues with waiting jobs at each population level.
We establish these laws for one-job priority, join-idle-queue, and idle-first sampling.
One-job priority attains the join-the-shortest-queue population law and mean delays; further queue information can improve limiting delays in convex order without changing their means.
In the super-Halfin--Whitt and Halfin--Whitt regimes, every such policy has the same stationary limits.
Our estimates bound idleness-weighted occupancy tails in terms of ordinary tails, uniformly over queue-length thresholds and policies.
They also yield process convergence and limit interchange for initial laws uniformly dominated by stationarity, with convergent scaled initial states.
\end{abstract}

{\small
\par
}

\section{Introduction}
\label{sec:introduction}

Routing arriving jobs among parallel queues requires deciding which queue-length information to use.
A dispatcher may locate idle servers, additionally identify queues with one job, or compare longer queues.
For homogeneous exponential servers, the optimality properties of join-the-shortest-queue (JSQ) make it a benchmark for dispatching policies~\cite{Winston1977}.
The search for policies that retain its performance with less information has led to randomized sampling~\cite{Mitzenmacher1999} and join-idle-queue (JIQ)~\cite{LuEtAl2011}, among many other schemes.
We ask which of these distinctions determine stationary congestion as the number of servers grows, and when policies using different information have the same limit.

When aggregate spare capacity grows on the square-root scale of the number of servers, the Halfin--Whitt diffusion limits of JIQ and JSQ agree on finite time intervals~\cite{MukherjeeBorstVanLeeuwaardenWhiting2016}.
A different picture emerges in nondegenerate slowdown (NDS), where spare capacity remains fixed as the number of servers grows.
Gupta and Walton~\cite{GuptaWalton2019} proved a JSQ diffusion limit and proposed approximations for JIQ and idle-one-first (I1F).
JIQ uses an idle server whenever possible and otherwise routes uniformly; I1F additionally gives priority to queues containing one job.
These approximations predict greater congestion under JIQ, while I1F has the same population law and mean delay as JSQ.
Their concluding discussion called for rigorous derivations of these proposed limits, steady-state limit interchange, and a systematic treatment of occupancy dynamics.

We study homogeneous parallel queues with Poisson arrivals, exponential service, and infinite buffers.
Our basic routing conditions give priority to idle servers and require the chosen queue to be stochastically no longer than a uniformly selected queue, conditional on the queue state (Assumption~\ref{ass:base}).
Which information matters depends on spare capacity.
In the super-Halfin--Whitt and Halfin--Whitt regimes, where spare capacity grows below or on the square-root scale, idle-server information suffices: JIQ and every other policy in this class have the same scaled stationary limits as JSQ (Theorems~\ref{thm:superhw} and~\ref{hw:stationary}).
In the NDS regime, JIQ has a more congested limiting population law than JSQ, while one-job priority suffices to recover the JSQ population law and limiting mean waiting and sojourn times.
We also identify the limiting stationary laws for idle-first sampling with a fixed number of samples.
Once one-job priority is imposed, further information about longer queues can improve the delay distribution without changing these aggregate limits.

To compare delay distributions, we consider a policy that chooses a shortest queue whenever the minimum queue length is at most a specified threshold and uses a fixed routing rule otherwise.
We call this remaining routing choice the \emph{fallback rule}; uniform routing and choosing a shortest queue among a fixed number of independent uniform samples are two examples.
For every fixed fallback rule in the family studied here, increasing the priority threshold beyond one job decreases limiting waiting and sojourn times in convex order while preserving their means.
Thus, for every target, it reduces the expected amount by which delay exceeds that target.
Every diverging priority threshold recovers the full JSQ stationary profile and delay laws when the remaining routing satisfies Assumption~\ref{ass:base}.

The mechanism behind the NDS results is the relation between queue occupancy and unused service capacity.
Flow balance fixes the same mean number of idle servers for every stable policy at a given arrival rate.
With service rates normalized to one, the population drift equals the number of idle servers minus spare capacity, so its average at a given population depends on the expected number of idle servers conditional on that population.
In NDS, the number of idle servers remains on the scale of the fixed spare capacity even though the idle fraction vanishes.
Service completions at one-job queues create idle servers, while arrivals remove idle servers whenever one is available.
A larger limiting fraction of one-job queues at a given population therefore yields greater conditional mean idleness.

By filling one-job queues before extending longer queues, one-job priority produces the same limiting fraction of one-job queues, and hence the same conditional idleness, as JSQ.
JIQ can instead send an arrival to a longer queue while a one-job queue remains, producing more conditional idleness at the same population level.
With growing spare capacity, jobs beyond the second position become negligible on the relevant congestion scale.
Routing among longer queues then disappears from the limiting population and idleness dynamics, yielding the common stationary limits.

The general NDS result formalizes this mechanism through an approximation linking the fraction of queues with waiting jobs to the total population.
The approximation must remain accurate after its error is multiplied by the number of idle servers (Assumption~\ref{ass:weighted-occupancy}).
Together with the basic routing conditions, it determines the stationary population law and limiting mean delays (Theorem~\ref{thm:regular-nds}).
Policies with the same limiting occupancy relation therefore have the same population law and mean delays.
Under these two conditions, the limiting population is stochastically bounded below by the JSQ limit and above by the JIQ limit.

Theorem~\ref{lem:tail-transfer} controls idleness-weighted occupancy tails through ordinary tails, uniformly over queue-length thresholds and eligible policies.
Proposition~\ref{prop:weighted-occupancy-criterion} combines these bounds with stationary idleness balance to derive the weighted approximation from ordinary occupancy convergence and ordinary queue-tail control.
We verify the condition directly under one-job priority, uniformly over the eligible routing choices among longer queues.
For the fixed fallback family in Example~\ref{ex:fallback-nds}, which includes JIQ and idle-first sampling, we establish ordinary profile convergence and tail control and apply this criterion.

Theorem~\ref{thm:common-interchange} gives process convergence to the corresponding diffusion in each regime on every finite time interval.
The conclusion holds in stationarity and from initial laws on full queue configurations bounded by a common constant multiple of the corresponding stationary laws, provided the scaled initial observations converge in distribution.
For these initial laws, taking the many-server and long-time limits in either order gives the same invariant law, including convergence of scaled population means.

\paragraph{Relation to earlier work.}
Table~\ref{tab:regime-breadth} compares the available process and stationary limits.
For this comparison, let $N$ be the number of rate-one servers, let $\Delta_N$ denote aggregate spare capacity, and fix $\beta>0$.

\begin{table}[htbp]
\centering
\small
\renewcommand{\arraystretch}{1.2}
\newcommand{\regime}[2]{\multicolumn{4}{@{}l@{}}{\textbf{#1}\quad$#2$}}
\begin{tabular*}{\textwidth}{@{\extracolsep{\fill}}>{\raggedright\arraybackslash\leftskip=1em}p{0.42\textwidth}ccc@{}}
\toprule
\multicolumn{1}{@{}l}{Routing class} & Process limit & Stationary convergence & Limit interchange\\
\midrule
\regime{Nondegenerate slowdown}{\Delta_N=\beta}\\
JSQ & \cite{GuptaWalton2019} & Theorem~\ref{thm:regular-nds} & Theorem~\ref{thm:common-interchange}\\
Assumptions~\ref{ass:base} and~\ref{ass:weighted-occupancy}, including I1F, JIQ, and idle-first sampling & Theorem~\ref{thm:common-interchange} & Theorem~\ref{thm:regular-nds} & Theorem~\ref{thm:common-interchange}\\
\addlinespace
\regime{Super-Halfin--Whitt}{1\ll\Delta_N\ll\sqrt N}\\
JSQ with power-law $\Delta_N$ & \cite{ZhaoBanerjeeMukherjee2026} & \cite{ZhaoBanerjeeMukherjee2026} & \cite{ZhaoBanerjeeMukherjee2026}\\
All policies under Assumption~\ref{ass:base} & Theorem~\ref{thm:common-interchange} & Theorem~\ref{thm:superhw} & Theorem~\ref{thm:common-interchange}\\
\addlinespace
\regime{Halfin--Whitt}{\Delta_N=\beta\sqrt N}\\
JSQ & \cite{EschenfeldtGamarnik2018} & \cite{Braverman2020} & \cite{Braverman2020}\\
Idle-first sampling, including JIQ & \cite{MukherjeeBorstVanLeeuwaardenWhiting2016} & Theorem~\ref{hw:stationary} & Theorem~\ref{thm:common-interchange}\\
All policies under Assumption~\ref{ass:base} & Theorem~\ref{thm:common-interchange} & Theorem~\ref{hw:stationary} & Theorem~\ref{thm:common-interchange}\\
\bottomrule
\end{tabular*}
\caption{Process limits, convergence of finite-system stationary distributions, and limit interchange, by regime and routing class.}
\label{tab:regime-breadth}
\end{table}

Following Atar's NDS framework for pooled systems~\cite{Atar2012}, Gupta and Walton~\cite[Section~5]{GuptaWalton2019} identified the role of one-job queues in capacity loss and explained how information about longer queues can improve higher delay moments while preserving the limiting mean.
Our NDS results establish stationary convergence to their JSQ invariant law and to the laws underlying their JIQ and I1F approximations, together with a general occupancy criterion for population limits.

In the super-Halfin--Whitt window, Liu and Ying~\cite{LiuYing2022} obtained stationary scaling bounds for a class of policies with finite buffers.
Zhao, Banerjee, and Mukherjee~\cite{ZhaoBanerjeeMukherjee2026} established the JSQ diffusion and its stationary Gamma limit under power-law spare-capacity scaling.
Our common estimates extend this stationary limit to every eligible policy sequence and to arbitrary spare-capacity sequences within the window.

In the Halfin--Whitt regime~\cite{HalfinWhitt1981}, Eschenfeldt and Gamarnik~\cite{EschenfeldtGamarnik2018} obtained the two-dimensional JSQ diffusion, and Mukherjee, Borst, van Leeuwaarden, and Whiting~\cite{MukherjeeBorstVanLeeuwaardenWhiting2016} established the same transient limit for idle-first sampling, including JIQ.
Braverman~\cite{Braverman2020} proved the JSQ stationary approximation and diffusion ergodicity, and later obtained convergence rates for finite buffers~\cite{Braverman2023}.
Our Halfin--Whitt result establishes stationary convergence and limit interchange throughout the class in Assumption~\ref{ass:base}.
For ordinary shortest-of-sample routing, the coupling analysis of Mukherjee et al.~\cite{MukherjeeBorstVanLeeuwaardenWhiting2018} also gives the JSQ diffusion when the sample size grows sufficiently quickly.
This sampling rule can miss idle servers outside its sample and need not satisfy our exact idle-priority condition.

At fixed subcritical load, Stolyar~\cite{Stolyar2015} proved that the stationary probability of waiting vanishes under pull-based routing, while Foss and Stolyar~\cite{FossStolyar2017} obtained the same conclusion for JIQ with general service times when the load is below one half.
With finite buffers and sub-Halfin--Whitt power-law scaling, Liu and Ying~\cite{LiuYing2020} proved vanishing stationary waiting probabilities for a policy class including JSQ and JIQ.
Liu, Gong, and Ying~\cite{LiuGongYing2022} obtained analogous guarantees for Coxian-2 service times; both analyses use Stein's method and state-space collapse.
The present paper studies the more congested NDS, super-Halfin--Whitt, and Halfin--Whitt regimes, with infinite buffers and nontrivial scaled congestion limits.

Beyond the regimes studied here, Jhunjhunwala, Hurtado-Lange, and Maguluri~\cite{JhunjhunwalaHurtadoMaguluri2024} obtain a stationary exponential limit for JSQ when spare capacity is $N^{1-\alpha}$ with $\alpha>1$.
Banerjee, Budhiraja, and Estevez~\cite{BanerjeeBudhirajaEstevez2026} obtain a rank-gap diffusion and stationary interchange for a fixed number of queues.
For growing systems under extreme congestion, Banerjee, Budhiraja, and Loeser~\cite{BanerjeeBudhirajaLoeser2026} obtain an infinite-dimensional rank-gap limit.
These rank-gap models allow routing without exact idle priority.

The uniform estimates build on mixed workload--unused-service bounds in steady-state heavy traffic~\cite{EryilmazSrikant2012,ZhouShroff2020,HurtadoLangeMaguluri2022} and the idleness--congestion estimates of Liu and Ying~\cite{LiuYing2022}.
Related moment estimates support stationary approximations for multiclass networks: Ye and Yao~\cite{YeYao2018} justify interchange under moment conditions on interarrival and service times, while Braverman, Dai, and Miyazawa~\cite{BravermanDaiMiyazawa2025} use moment state-space collapse to control boundary-conditioned transforms under static buffer priority.
The fixed-fallback profile analysis connects with Mitzenmacher's closed-supermarket equations~\cite{Mitzenmacher1999} and uses invariant-measure convergence in the spirit of Bena\"im and Le Boudec~\cite{BenaimLeBoudec2011}.
Our cubic potential and product identity give tail bounds uniform over thresholds and policies, preserving the profile mean and its contribution to the population drift.

\paragraph{Organization.}
Section~\ref{sec:model} introduces the model and presents the main results on stationary congestion, limit interchange, and customer delay.
Section~\ref{sec:common} develops the uniform stationary estimates and the methods for identifying limiting stationary laws and averaging idleness in the population dynamics.
Section~\ref{sec:main-proofs} applies these tools to prove the main results, with supporting technical proofs collected in the appendices.

\section{Routing information and stationary performance}
\label{sec:model}

There are $N$ homogeneous single-server queues with infinite buffers.
Arrivals form a Poisson process of rate $\lambda_N=N-\Delta_N$, where $0<\Delta_N<N$.
Service times are independent rate-one exponentials, independent of the arrival process.
An arrival is routed immediately, remains at its chosen server, and receives first-come, first-served (FCFS) service.
Routing probabilities depend only on the current queue-length vector and $N$ and are time-homogeneous.
For $i=1,\ldots,N$, let $Q_i^N(t)$ be the number of jobs at server $i$ at time $t$, including any job in service, and write $Q^N(t)=(Q_i^N(t))_{i=1}^N$.
We write $Q^N=(Q_i^N)_{i=1}^N$ for a \textit{stationary} queue-length vector and define
\[
 S_N=\sum_{i=1}^N Q_i^N,\qquad I_N=\#\{i:Q_i^N=0\},\qquad
 Z_{m,N}=\#\{i:Q_i^N\ge m\},\qquad X_N=S_N/N.
\]
Here $S_N$ is the total population, $I_N$ is the number of idle servers, $Z_{m,N}$ counts queues with at least $m$ jobs, and $X_N$ is the population per server.
The index $m$ ranges over the nonnegative integers $\N_0=\{0,1,\ldots\}$, so $Z_{0,N}=N$ and $Z_{1,N}=N-I_N$.
The total number of waiting jobs is $W_N=S_N-N+I_N$.
We write $S_N(t)$ and $I_N(t)$ for the total population and number of idle servers at time $t$, with analogous notation for the other quantities.
We use $\ell^1$ for the space of absolutely summable sequences, with norm $\|u\|_1=\sum_{k\ge1}|u_k|$.
The empirical tail profile of waiting jobs per queue is $u^N=(u_k^N)_{k\ge1}$, where $u_k^N=Z_{k+1,N}/N$ and $\|u^N\|_1=W_N/N$; set $u_0^N=1$.

For a real number $x$, write $x^+=\max\{x,0\}$, and for an event $B$, let $\ind_B$ denote its indicator.
We write $\Rightarrow$ for weak convergence and $W_p$ for Wasserstein distance of order $p\ge1$, using the Euclidean metric unless another metric is stated.
For a deterministic queue state $q=(q_i)_{i=1}^N\in\N_0^N$, notation such as $I_N(q)$ and $Z_{m,N}(q)$ denotes the same statistics evaluated at $q$.
Let $\iota$ denote the server selected by an arrival in that state.
All routing rules satisfy the following condition.

\begin{assumption}[Idle priority and comparison with uniform routing]
\label[assumption]{ass:base}
In every queue state $q$, an arrival joins an idle server whenever $I_N(q)>0$.
For every $q\in\N_0^N$ and every integer $m\ge0$,
\begin{equation}
 \Pp(q_\iota\ge m\mid q)\le Z_{m,N}(q)/N.
 \label{eq:domination}
\end{equation}
\end{assumption}

The inequality~\eqref{eq:domination} means that, conditional on the queue state, the chosen queue is stochastically no longer than a uniformly selected queue.
JSQ, I1F, JIQ, and idle-first shortest-of-$d$ sampling with any integer $d\ge1$ satisfy Assumption~\ref{ass:base}.
When all servers are busy, uniform routing satisfies~\eqref{eq:domination} with equality.
Choosing a shortest queue, or a one-job queue when one is available, selects a queue of minimum length and therefore satisfies the inequality.
Choosing a shortest queue among $d\ge1$ independent uniform samples gives conditional probability $(Z_{m,N}(q)/N)^d\le Z_{m,N}(q)/N$ of selecting a queue with at least $m$ jobs, for every $m\ge0$.
All these policies require locating an idle server whenever one exists.
The communication needed to implement this requirement depends on how the dispatcher obtains and updates queue information.

Under Assumption~\ref{ass:base}, the queue-length process has a unique stationary law $\Pi_N$ with finite population moments of every order for each fixed $N$ (Lemma~\ref{lem:foundations}).

\subsection{Stationary limits across loading regimes}
\label{sec:performance}

We compare stationary congestion as spare capacity increases from the NDS scale to the Halfin--Whitt scale.

\paragraph{The NDS regime.}
The spare capacity is fixed: $\Delta_N=\beta>0$.
Flow balance gives $\E I_N=\beta$ for every policy considered here.
This common average utilization leaves room for different congestion distributions because unused capacity can occur at different population levels.
Since the total-population drift is $I_N-\beta$, the relevant quantity is $\E[I_N\mid X_N]$, the mean number of idle servers conditional on the population per server.

The fraction of queues containing at least two jobs provides the routing information needed to identify the corresponding conditional mean in the many-server limit.
Let $A:(1,\infty)\to(0,1]$ be a candidate approximation, fixed as $N$ varies.
For each compact set $K\subset(1,\infty)$, the integrand below is defined to be zero when $X_N\notin K$.

\begin{assumption}[Weighted occupancy approximation]
\label[assumption]{ass:weighted-occupancy}
The function $A$ is locally Lipschitz and, for every compact $K\subset(1,\infty)$,
\begin{equation}
 \E\left[(I_N+1)\left|\frac{Z_{2,N}}N-A(X_N)\right|\ind_{\{X_N\in K\}}\right]\longrightarrow0.
 \label{path:weighted-coefficient}
\end{equation}
\end{assumption}

Thus $A(x)$ approximates the fraction of queues with waiting jobs when the population per server is $x$.
The factor $I_N+1$ requires ordinary and idleness-weighted mean errors to vanish on each compact population set, controlling the approximation in the population drift.
Proposition~\ref{prop:weighted-occupancy-criterion} derives this weighted approximation from ordinary occupancy convergence and ordinary queue-tail control.
For example, JSQ and I1F satisfy this condition with $A(x)=\min\{x-1,1\}$, while JIQ satisfies it with $A(x)=(x-1)/x$.
Examples~\ref{ex:priority-nds}--\ref{ex:sampling-nds} give the resulting laws and the extension to fixed fallback rules.

Since $I_N/N\to0$, the corresponding fraction of one-job queues is $1-A(x)$.
A completion at one of these queues creates an idle server, while an arrival removes an idle server whenever one exists.
The normalized birth and death rates of idleness are therefore $1-A(x)$ and one, giving a geometric equilibrium with conditional mean
\[
 g(x)=\frac{1-A(x)}{A(x)},\qquad x>1.
\]
Although the idle fraction vanishes, the number of idle servers remains on the same scale as the spare capacity $\beta$.
Its conditional mean contributes $g(x)-\beta$ to the population drift.
On the time scale $Nt$, the normalized arrival and service fluctuations have total variance rate two.
These drift and variance coefficients lead to the candidate stationary population density and probability law
\begin{equation}
 p_{g,\beta}(x)=\frac{1}{C_{g,\beta}}
 \exp\left\{-\beta x+\int_2^xg(u)\,\dd u\right\},\qquad x>1,
 \qquad \pi_{g,\beta}(\dd x)=p_{g,\beta}(x)\,\dd x,
 \label{eq:generic-density}
\end{equation}
where $C_{g,\beta}$ is the normalizing constant.
The density $p_{g,\beta}$ solves the stationary balance equation $p_{g,\beta}'(x)=\{g(x)-\beta\}p_{g,\beta}(x)$.
Let $T_N^{\mathrm w}$ and $T_N^{\mathrm s}$ denote the waiting and sojourn times of a stationary arrival.
For a random population $X>1$, define the idleness variable $I$ conditionally by
\[
 \Pp(I=i\mid X=x)=A(x)\{1-A(x)\}^i,\qquad i\in\N_0,
\]
with $I=0$ when $A(X)=1$.

The next theorem uses the weighted occupancy approximation in Assumption~\ref{ass:weighted-occupancy} to identify the joint stationary limit of population and idleness.
It also supplies the moment convergence needed to recover limiting mean waiting and sojourn times from the population law.

\begin{theorem}[NDS stationary congestion and conditional capacity loss]
\label[theorem]{thm:regular-nds}
Under Assumptions~\ref{ass:base} and~\ref{ass:weighted-occupancy},
\begin{equation}
 \frac{x-1}{x}\le A(x)\le\min\{x-1,1\},\qquad x>1.
 \label{nds:occupancy-bounds}
\end{equation}
The density \eqref{eq:generic-density} is well defined, and $X_N$ converges to $X\sim\pi_{g,\beta}$ in Wasserstein distance of every order $1\le p<2$.
Moreover,
\[
 (X_N,I_N)\Rightarrow(X,I),
\]
and $\E[I\mid X=x]=g(x)$, $\E I=\beta$, and $(I_N)$ is uniformly integrable.
The stationary customer means and immediate-service probability satisfy
\[
 \E T_N^{\mathrm w}\longrightarrow\E(X-1),\qquad
 \E T_N^{\mathrm s}\longrightarrow\E X,\qquad
 \Pp(T_N^{\mathrm w}=0)\longrightarrow\E[1-A(X)].
\]
\end{theorem}

The corresponding conditional mean idleness satisfies
\begin{equation}
 \frac{(2-x)^+}{x-1}\le g(x)\le\frac1{x-1},\qquad x>1,
 \label{nds:idleness-bounds}
\end{equation}
and $\E g(X)=\beta$.
The estimate \eqref{path:weighted-coefficient} also justifies the drift averaging used to prove process convergence and interchange in Theorem~\ref{thm:common-interchange}.
Consequently, policies with the same $A$ have the same limiting population law, population diffusion, and mean delays.

\begin{example}[One-job priority: JSQ and I1F]
\label[example]{ex:priority-nds}
A one-job priority policy satisfies Assumption~\ref{ass:base} and routes to a one-job queue whenever no server is idle and a one-job queue exists.
This class includes JSQ and I1F, which uses uniform routing when every queue has at least two jobs.
More generally, the routing used when every queue has at least two jobs may depend on the full state and on $N$, subject to \eqref{eq:domination}.
A state-dependent choice of rank $R$ among queues ordered from shortest to longest satisfies~\eqref{eq:domination} whenever $\Pp(R\le r\mid q)\ge r/N$ for every rank $r$.
Such a rule belongs to this class when it also gives priority to idle and one-job queues.
Every sequence of such policies verifies Assumption~\ref{ass:weighted-occupancy} with
\[
 A(x)=\min\{x-1,1\},\qquad
 g(x)=g_*(x):=\frac{(2-x)^+}{x-1}.
\]
Proposition~\ref{pri:weighted-occupancy} proves the required approximation uniformly over the eligible fallback rules.
Thus JSQ, I1F, and state- or system-size-dependent one-job priority policies share the population density
\begin{equation}
 p_{g_*,\beta}(x)=\frac{\beta(1+\beta)^2}{1+\beta e^{1+\beta}}
 \begin{cases}
 (x-1)e^{-(\beta+1)(x-2)},&1<x\le2,\\
 e^{-\beta(x-2)},&x\ge2.
 \end{cases}
 \label{eq:priority-density}
\end{equation}
\end{example}

\begin{example}[Idle priority with a fixed fallback function]
\label[example]{ex:fallback-nds}
This family uses an idle server whenever possible and a fixed rule for choosing among busy servers otherwise.
We describe the latter rule through the tail probabilities of the selected queue length.
Let $F:[0,1]\to[0,1]$ be a nondecreasing Lipschitz function satisfying
\[
 F(0)=0,\qquad F(1)=1,\qquad F(z)\le z\quad(0\le z\le1).
\]
A policy has fallback rule $F$ if it gives exact priority to idle servers and, in every state $q$ with $I_N(q)=0$, satisfies
\begin{equation}
 \Pp(q_\iota\ge m\mid q)=F(Z_{m,N}(q)/N),\qquad m\ge1.
 \label{eq:F-routing}
\end{equation}
The function $F$ is fixed as $N$ varies.

At a fixed population level, let $a$ be the first coordinate of a candidate equilibrium profile.
The geometric idleness balance above gives limiting probability $a$ that all servers are busy.
Fallback is therefore active for a fraction $a$ of arrivals, and balance across queue-length thresholds yields the recurrence below.
For $y=x-1>0$ waiting jobs per server, define $q^F(y)=(q_k^F(y))_{k\ge1}$ using the unique parameter $a_F(y)\in(0,1)$ satisfying
\[
 q_0^F(y)=1,\qquad q_1^F(y)=a_F(y),\qquad
 q_{k+1}^F(y)=a_F(y)F(q_k^F(y))\quad(k\ge1),\qquad
 \sum_{k\ge1}q_k^F(y)=y.
\]
At $y=0$, set $a_F(0)=0$, $q_0^F(0)=1$, and $q_k^F(0)=0$ for $k\ge1$.
Theorem~\ref{thm:profile} establishes existence, uniqueness, and attraction of this equilibrium profile.
Proposition~\ref{reg:weighted-law} verifies Assumption~\ref{ass:weighted-occupancy} with
\[
 A(x)=a_F(x-1),\qquad
 g(x)=g_F(x):=\frac{1-a_F(x-1)}{a_F(x-1)},
\]
and also gives the full stationary profile limit
\[
 (X_N,u^N)\Rightarrow(X,q^F(X-1))
 \quad\text{in }\R\times\ell^1,
 \qquad X\sim\pi_{g_F,\beta}.
\]
The profile coordinate $q_1^F(x-1)=a_F(x-1)$ determines population congestion through Theorem~\ref{thm:regular-nds}; the full profile identifies customer delay.
\end{example}

\begin{example}[Join-idle-queue]
\label[example]{ex:jiq-nds}
JIQ, also called idle-queue-first (IQF), gives priority to idle servers and uses uniform routing when every server is busy, so $F(z)=z$.
It uses the idle-server information underlying the low-communication policy of Lu et al.~\cite{LuEtAl2011}.
Its equilibrium profile is geometric:
\[
 a_F(y)=\frac{y}{1+y},\qquad q_k^F(y)=\left(\frac{y}{1+y}\right)^k.
\]
Hence $A(x)=(x-1)/x$, $g(x)=1/(x-1)$, and
\begin{equation}
 X-1\sim\operatorname{Gamma}(2,\beta).
 \label{eq:jiq-law}
\end{equation}
Here and below, the Gamma parameters specify shape and rate, in that order.
The JIQ law \eqref{eq:jiq-law} and the one-job priority density \eqref{eq:priority-density} validate the stationary population laws underlying the I1F and IQF approximations of Gupta and Walton~\cite{GuptaWalton2019}, alongside their rigorous JSQ diffusion limit.
For example, at $\beta=1$, the limiting mean waiting time is approximately $1.2384$ under one-job priority and $2$ under JIQ, an increase of about $61.5\%$.
\end{example}

\begin{example}[Idle-first sampling]
\label[example]{ex:sampling-nds}
Choosing an idle server whenever possible and otherwise the shortest of $d$ independent uniform samples with replacement gives $F(z)=z^d$ for a fixed integer $d\ge1$.
More generally, each arrival that finds all servers busy may draw an integer-valued sample size $D\ge1$ from a fixed distribution independent of the queue state and then choose the shortest of $D$ independent uniform samples.
If $\E D<\infty$, this gives $F(z)=\E z^D$.
The bound $|z^d-w^d|\le d|z-w|$ for $z,w\in[0,1]$ makes $F$ Lipschitz with constant at most $\E D$, and $D\ge1$ gives $F(z)\le z$.
These functions, including fixed convex mixtures, satisfy the conditions on $F$ in Example~\ref{ex:fallback-nds}, so its conclusions apply.
The use of global idle-server information is the same as in the JIQ($d$) architecture of Mukherjee et al.~\cite{MukherjeeBorstVanLeeuwaardenWhiting2016}.
The usual JSQ($d$) rule can miss an idle server outside its sample and need not satisfy Assumption~\ref{ass:base}.
\end{example}

The bounds on conditional mean idleness also yield a comparison of limiting population laws for policies satisfying Assumptions~\ref{ass:base} and~\ref{ass:weighted-occupancy}.
At a fixed $\beta$, let $X_{\mathrm{JSQ}}$, $X_g$, and $X_{\mathrm{JIQ}}$ have the corresponding limiting population laws.
Equations~\eqref{eq:generic-density} and~\eqref{nds:idleness-bounds} give
\[
 X_{\mathrm{JSQ}}\le_{\mathrm{lr}}X_g\le_{\mathrm{lr}}X_{\mathrm{JIQ}}.
\]
Here $\le_{\mathrm{lr}}$ denotes \emph{likelihood ratio order}: the ratio of the larger law's density to the smaller law's density is nondecreasing.
This order implies the usual stochastic order, denoted by $\le_{\mathrm{st}}$.
Indeed, the derivative of $\log(p_{g,\beta}/p_{g_*,\beta})$ is $g-g_*\ge0$, and the upper comparison follows in the same way.
Thus JSQ and JIQ bound the limiting population and mean delay over this class, and every one-job priority policy attains the JSQ population law.

Within the fallback family, $F\le G$ gives $g_F\le g_G$ and the same likelihood ratio ordering of population laws; Proposition~\ref{prof:ordering} proves this comparison.
In particular, increasing the number $d$ of fallback samples decreases stationary congestion.

For growing spare capacity, Assumption~\ref{ass:base} yields a common stationary limit in each of the two regimes below.
We call a stationary limit \emph{universal} over a routing class when every eligible sequence of policies has the same limiting law for the stated scaled observable.

\paragraph{The super-Halfin--Whitt regime.}
Assume
\begin{equation}
 \Delta_N\longrightarrow\infty,\qquad \Delta_N/\sqrt N\longrightarrow0.
 \label{shw:scaling}
\end{equation}
The waiting population is then of order $N/\Delta_N=o(N)$.
When all servers are busy, the bound~\eqref{eq:domination} limits the probability of joining a queue that already contains waiting jobs by $Z_{2,N}/N$.
Lemma~\ref{shw:collapse} shows that jobs beyond the second position are negligible on the waiting-population scale and under the idleness weight entering the population equation.
The averaged capacity balance can therefore be written in terms of the waiting population alone, yielding a common population limit for all policies satisfying Assumption~\ref{ass:base}.

Set $Y_N=\Delta_N(S_N-N)/N$.
Let $G$ have the $\operatorname{Gamma}(2,1)$ law, denoted by $\pi_{\mathrm{shw}}$.

\begin{theorem}[Stationary universality in the super-Halfin--Whitt window]
\label[theorem]{thm:superhw}
Under \eqref{shw:scaling}, every sequence of policies satisfying Assumption~\ref{ass:base} obeys
\begin{equation}
 Y_N\Rightarrow G,\qquad \E Y_N\longrightarrow2,
 \label{eq:superhw-main}
\end{equation}
with convergence in Wasserstein distance of every order $1\le p<2$.
The same conclusions hold for $\Delta_NW_N/N$ and $\Delta_NZ_{2,N}/N$.
\end{theorem}

Theorem~\ref{thm:superhw} extends the JSQ Gamma limit of Zhao, Banerjee, and Mukherjee~\cite[Theorem~2.7]{ZhaoBanerjeeMukherjee2026} to this routing and scaling scope.

The NDS laws recover this Gamma law as their spare-capacity parameter grows.
For fixed $g$ satisfying~\eqref{nds:idleness-bounds}, let $X^{(\beta)}\sim\pi_{g,\beta}$ and put $r(y)=\int_2^{1+y}g(u)\,\dd u-\log y$.
The bounds give $r(y)\le1$ for $y>0$ and $-1\le r'(y)\le0$ for $0<y\le1$, so $r(0+)$ is finite.
The density of $\beta(X^{(\beta)}-1)$ is proportional to $ze^{-z}\exp\{r(z/\beta)\}$, and dominated convergence yields $\beta(X^{(\beta)}-1)\Rightarrow G$ as $\beta\to\infty$.

\paragraph{The Halfin--Whitt regime.}
Set $\Delta_N=\beta\sqrt N$ for fixed $\beta>0$.
Lemma~\ref{hw:collapse} shows that jobs beyond the second position are negligible on the $\sqrt N$ scale, while idleness and the number of queues with waiting jobs both remain visible on that scale.
Their reflected equations depend only on the two conditions in Assumption~\ref{ass:base}.
Let $U_N=-I_N/\sqrt N$ and $V_N=Z_{2,N}/\sqrt N$.
The limiting reflected diffusion on $(-\infty,0]\times[0,\infty)$ is driven by a standard Brownian motion $B$ and a nondecreasing regulator $R$ with $R(0)=0$:
\begin{align}
 U(t)&=U(0)+\sqrt2 B(t)-\beta t+
       \int_0^t\{-U(s)+V(s)\}\,\dd s-R(t),\label{hw:limit-u}\\
 V(t)&=V(0)-\int_0^t V(s)\,\dd s+R(t).\label{hw:limit-v}
\end{align}
The regulator increases only when $U(t)=0$, as expressed by
\[
 \int_0^\infty\ind_{\{U(t)<0\}}\,\dd R(t)=0.
\]
The diffusion has a unique invariant law, denoted by $\Pi_\beta$.
The solution map is well posed by Eschenfeldt and Gamarnik~\cite[Lemma~1]{EschenfeldtGamarnik2018}, and existence, uniqueness, and exponential ergodicity of its invariant law follow from Braverman~\cite[Theorem~3 and Corollary~1]{Braverman2020}.
Let $(U,V)$ have law $\Pi_\beta$.

\begin{theorem}[Halfin--Whitt stationary universality]
\label[theorem]{hw:stationary}
For every sequence of policies satisfying Assumption~\ref{ass:base},
\[
 (U_N,V_N)\Rightarrow(U,V),
\]
in stationarity, with convergence in Wasserstein distance of every order $1\le p<2$.
Consequently,
\[
 \frac{S_N-N}{\sqrt N}\Rightarrow U+V,\qquad
 \frac{\E(S_N-N)}{\sqrt N}\longrightarrow\E_{\Pi_\beta}(U+V).
\]
\end{theorem}

The stationary theorem covers every policy satisfying Assumption~\ref{ass:base}.
It complements the transient universality for idle-first sampling established by Mukherjee et al.~\cite{MukherjeeBorstVanLeeuwaardenWhiting2016} and the JSQ stationary approximation of Braverman~\cite{Braverman2020}.

\subsection{Process limits and interchange of limits}
\label{sec:process-limits}

The invariant laws of the diffusion limits below coincide with the corresponding stationary limits in Section~\ref{sec:performance}.
Under the initialization conditions specified below, taking the many-server limit before the long-time limit yields the same invariant law as taking them in the reverse order.
For NDS and super-Halfin--Whitt, respectively, define
\[
 \widehat X_N(t)=S_N(Nt)/N,\qquad
 \widehat Y_N(t)=\frac{\Delta_N}{N}
 \left\{S_N\left(\frac{N}{\Delta_N^2}t\right)-N\right\}.
\]
Their limiting diffusions are
\begin{align}
 \dd X(t)&=\sqrt2\,\dd B(t)+\{g(X(t))-\beta\}\,\dd t,\qquad X(t)>1,
 \label{path:sde}\\
 \dd Y(t)&=\sqrt2\,\dd B(t)+\left\{\frac1{Y(t)}-1\right\}\,\dd t,\qquad Y(t)>0,
 \label{shw:limit-sde}
\end{align}
with invariant laws $\pi_{g,\beta}$ and $\pi_{\mathrm{shw}}$, respectively.
In the Halfin--Whitt regime, define $(U_N(t),V_N(t))=(-I_N(t),Z_{2,N}(t))/\sqrt N$ on the physical time scale; its limit is \eqref{hw:limit-u}--\eqref{hw:limit-v}.
Write $\mathcal X_N$ for the appropriate scaled process, $\mathcal X$ for its limiting diffusion, and $\pi$ for the corresponding invariant law.
Process convergence is in $D([0,T],\R)$ for the one-dimensional processes and $D([0,T],\R^2)$ for the pair, with the Skorokhod $J_1$ topology.

Recall that $\Pi_N$ is the stationary law of the full queue vector.
For nonstationary initial laws $\zeta_N$ on full queue vectors, assume independence from future arrival, service, and routing randomness and
\begin{equation}
 \zeta_N(\mathcal A)\le C\Pi_N(\mathcal A)\quad\text{for every set of queue vectors }\mathcal A,
 \qquad \Law_{\zeta_N}(\mathcal X_N(0))\Rightarrow\nu,
 \label{eq:common-dominated-initialization}
\end{equation}
where $C<\infty$ is independent of $N$.
This condition bounds the initial probability of every full-configuration event by a constant times its stationary probability, including events that specify how jobs are distributed among servers.
For example, in NDS, conditioning $\Pi_N$ on $X_N$ lying in a fixed interval with positive limiting probability satisfies this condition.

\begin{theorem}[Stationary processes and interchange of limits]
\label[theorem]{thm:common-interchange}
Suppose Assumption~\ref{ass:base} holds and, in NDS, Assumption~\ref{ass:weighted-occupancy} also holds.
In each of the three regimes, the stationary processes $\mathcal X_N$ converge on every finite time interval to the stationary diffusion $\mathcal X$.
Under \eqref{eq:common-dominated-initialization}, they converge instead to $\mathcal X$ with initial law $\nu$ and Brownian motion independent of its initial state, and
\begin{equation}
 \lim_{N\to\infty}\lim_{t\to\infty}\Law_{\zeta_N}(\mathcal X_N(t))
 =\pi
 =\lim_{t\to\infty}\lim_{N\to\infty}\Law_{\zeta_N}(\mathcal X_N(t)).
 \label{eq:common-two-orders}
\end{equation}
For every $1\le p<2$, the stationary marginal limits, the fixed-time limits just described, and the inner and iterated limits in \eqref{eq:common-two-orders} hold in $W_p$.
These convergences preserve expectations of continuous functions bounded in absolute value by $C'(1+\|x\|^a)$ for some $C'<\infty$ and $0<a<2$.
In particular, both orders of limits of the scaled coordinate means equal their means under $\pi$.
\end{theorem}

\subsection{Queue-length information and delay distributions in NDS}
\label{sec:delay-distribution}

Fix $\Delta_N=\beta>0$.
Under one-job priority, Example~\ref{ex:priority-nds} and Theorem~\ref{thm:regular-nds} give the same limiting population law and mean delays as JSQ.
We compare JSQ with the fixed fallback family of Example~\ref{ex:fallback-nds} (including JIQ and idle-first sampling) and its threshold-priority extensions.
For a fixed fallback, increasing priority beyond one job reduces expected delay in excess of every target while preserving the mean.

For a fixed integer $r\ge1$ and a fixed $F$ as in Example~\ref{ex:fallback-nds}, a threshold-priority rule chooses a shortest queue whenever the minimum queue length is at most $r$ and otherwise uses the fallback tail probabilities $F(Z_{m,N}(q)/N)$ in state $q$.
Uniform fallback gives the idle-$r$-first rule $\mathrm{I}_r\mathrm{F}$ discussed by Gupta and Walton~\cite[Section~5.2]{GuptaWalton2019}.

\paragraph{From occupancy to delay.}
For each policy, let $X$ have its limiting population law and write $q(x)=(q_m(x))_{m\ge1}$ for its limiting occupancy profile at population level $x>1$.
The coordinate $q_m(x)$ is the fraction of queues with at least $m$ waiting jobs, corresponding to $u_m^N=Z_{m+1,N}/N$, and $\sum_{m\ge1}q_m(x)=x-1$.
Let $D_N$ be the number of jobs already at the server selected by a stationary arrival, including any job in service.
Stationary level balance matches the rate of arrivals joining queues with at least $m$ jobs to the rate of service completions from queues with at least $m+1$ jobs.
For these policies, this balance holds conditional on the limiting population and identifies $q_m(x)$ as the limiting probability of at least $m$ jobs ahead.
Let $E_n$ be the sum of the first $n$ variables in a sequence of independent rate-one exponential random variables, with $E_0=0$.
By FCFS and memorylessness, the limiting number of jobs ahead $D$ and the waiting and sojourn times satisfy
\[
 \begin{gathered}
 \Pp(D\ge m\mid X=x)=q_m(x),\qquad m\ge1,\\
 \Law(T^{\mathrm w}\mid X,D=d)=\Law(E_d),\qquad
 \Law(T^{\mathrm s}\mid X,D=d)=\Law(E_{d+1}),\qquad d\ge0.
 \end{gathered}
\]
Lemma~\ref{ext:arrival-balance} proves convergence of the laws of $D_N$ and the two delays to these laws in Wasserstein distance of order one.

The policy-specific profiles make these delay laws explicit.
Idle priority with fallback $F$ has $q_m(x)=q_m^F(x-1)$.
For $x>1$, the JSQ profile is
\begin{equation}
 q_m^{\mathrm{JSQ}}(x)=\min\{1,(x-m)^+\},\qquad m\ge1,
 \label{ext:jsq-profile}
\end{equation}
and we extend it by zero for $x\le1$.
For priority threshold $r$ and fallback $F$, the profile is
\begin{align}
 q_m^{r,F}(x)&=\min\{1,(x-m)^+\},&&1\le m\le r,
 \notag\\
 q_{r+k}^{r,F}(x)&=q_k^F((x-r-1)^+),&&k\ge1.
 \label{ext:priority-profile}
\end{align}
The first $r$ coordinates coincide with JSQ; the remaining coordinates form the fallback profile for the excess population above level $r+1$.
Proposition~\ref{ext:profiles} proves the profile limits for JSQ and these threshold policies.

\paragraph{Comparing priority levels.}
At population level $x=3$, the limiting JSQ waiting time is the sum of two independent rate-one exponentials, whereas the I1F waiting time is exponential with mean two.
Both conditional means equal two.
Write $T_P^{\mathrm w}$ and $T_P^{\mathrm s}$ for the limiting waiting and sojourn times under policy $P$.
Use $T_{r,F}^{\mathrm w},T_{r,F}^{\mathrm s}$ for threshold $r$ with fallback $F$, and $T_F^{\mathrm w},T_F^{\mathrm s}$ for idle priority with fallback $F$.
For integrable variables, \emph{increasing convex order}, written $U\le_{\mathrm{icx}}V$, means $\E(U-t)^+\le\E(V-t)^+$ for every real $t$; \emph{convex order}, written $U\le_{\mathrm{cx}}V$, additionally requires equal means~\cite[Sections~3.A and~4.A]{ShakedShanthikumar2007}.
These orders compare expected delay in excess of every target.

\begin{proposition}[Delay comparison and information about longer queues]
\label[proposition]{ext:delay-order}
For every $F$ as in Example~\ref{ex:fallback-nds},
\[
 T^{\mathrm w}_{1,F}\le_{\mathrm{icx}}T^{\mathrm w}_{F},\qquad
 T^{\mathrm s}_{1,F}\le_{\mathrm{icx}}T^{\mathrm s}_{F}.
\]
For the same $F$ and every integer $r\ge1$,
\[
 T^{\mathrm w}_{\mathrm{JSQ}}
 \le_{\mathrm{cx}}T^{\mathrm w}_{r+1,F}
 \le_{\mathrm{cx}}T^{\mathrm w}_{r,F},\qquad
 T^{\mathrm s}_{\mathrm{JSQ}}
 \le_{\mathrm{cx}}T^{\mathrm s}_{r+1,F}
 \le_{\mathrm{cx}}T^{\mathrm s}_{r,F}.
\]
As $r\to\infty$, both threshold-priority delay laws converge in Wasserstein distance of order one to their JSQ counterparts.
\end{proposition}

The mean reaches its JSQ value at one-job priority; further priority levels reduce expected excess delay while preserving that mean.

\paragraph{Growing priority thresholds.}
The preceding proposition first takes the many-server limit at fixed $r$ and then lets $r\to\infty$.
For a threshold growing with the system, let $r_N\ge1$ be integers with $r_N\to\infty$, and let $X_*$ have density \eqref{eq:priority-density}.
On $\R\times\ell^1$, use the metric $|x-x'|+\|u-u'\|_1$.

\begin{corollary}[Growing shortest-queue priority thresholds]
\label[corollary]{ext:growing-thresholds}
Suppose a policy chooses a shortest queue whenever the minimum queue length is at most $r_N$ and otherwise uses any state- and $N$-dependent fallback satisfying~\eqref{eq:domination}.
Then
\[
 (X_N,u^N)\Rightarrow(X_*,q^{\mathrm{JSQ}}(X_*))
 \quad\text{in }\R\times\ell^1,
\]
with convergence also in Wasserstein distance of every order $1\le p<2$.
The stationary laws of the number of jobs ahead, waiting time, and sojourn time converge in Wasserstein distance of order one to their JSQ limiting laws.
\end{corollary}

\section{Uniform estimates for stationary limits and interchange}
\label{sec:common}

The total-population drift $I_N-\Delta_N$ makes it necessary to control how unused service capacity varies with congestion.
The balance equations relating idleness to occupancy contain products such as $I_NZ_{2,N}/N$.
Replacing $Z_{2,N}/N$ by a population-dependent approximation therefore produces an error multiplied by $I_N$.
The estimates below control queue tails and population moments in ordinary expectation and after multiplication by idleness.
In NDS, the tail estimates combine with stationary idleness balance to turn ordinary occupancy convergence and ordinary queue-tail control into the weighted approximation in Assumption~\ref{ass:weighted-occupancy}, as shown in Proposition~\ref{prop:weighted-occupancy-criterion}.

\subsection{Uniform stationary estimates}
\label{sec:mixed-potential}

In fixed-system calculations throughout the proofs, we suppress the index $N$, writing $Q=(Q_i)_{i=1}^N$ for $Q^N$, $\lambda=\lambda_N$, $\Delta=\Delta_N$, and similarly $S,I,Z_m,X,W$ for the corresponding stationary quantities.
For a queue vector $q$ and $j\ge0$, let $M_j=\#\{i:q_i=j\}$, so $M_0=I$ and $M_1=N-I-Z_2$.
For integers $k\ge1$ and $r\ge1$, write
\[
 H_k=\sum_{i=1}^N(q_i-k)^+,
 \qquad P_{r,k}=\sum_{i=1}^N\bigl((q_i-k)^+\bigr)^r,
 \qquad P_r=P_{r,1}.
\]
The statistic $H_k$, written $H_{k,N}$ when the system size varies, counts jobs above queue-length level $k$; in particular, $H_1=W=S-N+I$.
The potential $P_{2,k}$ is the sum of the squared excess queue lengths above level $k$.
Write $\rho=\lambda/N$.
The generator of the full queue-length chain is denoted by $\cL$.

\begin{lemma}[Fixed-system foundations]
\label[lemma]{lem:foundations}
Under \cref{ass:base}, the queue-length chain is irreducible and positive recurrent for every $0<\lambda<N$.
It has a unique stationary distribution with all polynomial moments of $S$ finite.
All polynomial generator identities used below are valid in stationarity.
\end{lemma}

The Foster--Lyapunov argument and the justification of the stationary generator identities are given in Appendix~\ref{app:common-technical}.

Applying the generator to population and idleness gives
\[
 \cL S=I-\Delta,
 \qquad
 \cL I=M_1-\lambda\ind_{\{I>0\}}
          =\Delta-I-Z_2+\lambda\ind_{\{I=0\}}.
\]
In stationarity, the first identity gives $\E I=\Delta$.
For the squared population, $\cL S^2=2S(I-\Delta)+\lambda+N-I$.
Taking expectations yields
\begin{equation}
 \Delta\E S=N-\Delta+\E[SI].
 \label{eq:population-square}
\end{equation}
This balance also underlies the many-server analysis of Hurtado-Lange and Maguluri~\cite{HurtadoLangeMaguluri2022}.
Since $S=N-I+W$, it can be written in terms of the waiting population as
\begin{equation}
 \Delta\E W=\E[IW]+N-\Delta+\Delta^2-\E I^2
 \le\E[IW]+N+\Delta^2.
 \label{eq:waiting-balance}
\end{equation}
Thus a bound on mean waiting population requires control of $\E[IW]$, which measures the coexistence of waiting jobs and idle servers.

For a general threshold $k$, the corresponding quantity is $\E[IH_k]$: it weights jobs above level $k$ by the number of idle servers in the same configuration.
We also control $\E[IH_k^2]$, the second moment of this excess population under the idleness weight.
Together with $\E I=\Delta$, this stronger bound gives the first-moment bound by Cauchy--Schwarz; at $k=1$, it also supplies the mixed moment needed to bound the population second moment.
The estimates keep their dependence on $N$ and $\Delta$ explicit and apply in NDS, super-Halfin--Whitt, and Halfin--Whitt.

\begin{theorem}[Uniform stationary tail estimates]
\label[theorem]{lem:tail-transfer}
Under \cref{ass:base}, for every $0<\Delta<N$ and $k\ge1$,
\begin{align}
 (\E[IH_k])^2&\le\Delta(2N+\Delta^2)\E H_k,
       \label{eq:main-weighted-estimate}\\
 \E[IH_k^2]&\le(2N+\Delta^2)\E H_k.
       \label{eq:tail-square}
\end{align}
Moreover,
\begin{align}
 \frac\Delta N\E P_{2,k}+\rho\E[P_{2,k}\ind_{\{I>0\}}]&\le2\E H_k, \label{eq:tail-cubic}\\
 \E[IZ_2P_{2,k}]&\le(2N+\Delta^2)\E H_k.
       \notag
\end{align}
The constants and coefficients are independent of the threshold $k$ and the routing policy.
\end{theorem}

To prove Theorem~\ref{lem:tail-transfer}, we use the following generator identity.
Let $\Phi$ be a nonnegative function of the queue state such that every transition changing $I$ leaves $\Phi$ unchanged.
Assume that $\cL\Phi=-C_\Phi$ on $\{I>0\}$, where $C_\Phi\ge0$, and that the stationary generator identity for $I(I-1)\Phi$ is justified.
These conditions hold for every nonnegative polynomial potential that depends only on the waiting-job counts $((Q_i-1)^+)_{i=1}^N$ and is nondecreasing in each coordinate.
Write $B_\Phi=\E[\Phi\ind_{\{I>0\}}]$.

\begin{lemma}[Generator identity for products]
\label[lemma]{mix:identity}
In stationarity,
\begin{equation}
 \E[IZ_2\Phi]
 =\lambda B_\Phi+\E[(\Delta I-I^2)\Phi]
       -\frac12\E[I(I-1)C_\Phi]
 \le\left(\lambda+\frac{\Delta^2}{4}\right)B_\Phi.
 \label{mix:identity-equation}
\end{equation}
\end{lemma}

\begin{proof}
No transition changes both $I$ and $\Phi$, so the cross term in the generator product rule vanishes.
For $G(i)=i(i-1)/2$, the upward increment is $i$ and the downward increment is $1-i$, giving
\[
 \cL G(I)=\Delta I-I^2-IZ_2+\lambda\ind_{\{I>0\}}.
\]
Stationarity of $G(I)\Phi$ proves the identity in \eqref{mix:identity-equation}.
Since $I$ is integer-valued, $I(I-1)C_\Phi\ge0$, and
$(\Delta I-I^2)\Phi\le(\Delta^2/4)\Phi\ind_{\{I>0\}}$.
These inequalities prove the bound.
\end{proof}

\begin{proof}[Proof of Theorem~\ref{lem:tail-transfer}]
Fix $k\ge1$ and write $m_k=\E H_k$ and $B_k=\E[P_{2,k}\ind_{\{I>0\}}]$.
Let $b_k(q)=\Pp(q_\iota\ge k\mid q)$ and $Z=Z_{k+1}$.
The exact linear tail generator is $\cL H_k=\lambda b_k-Z$.
For $u=(q_i-k)^+$, the arrival increment of $u^3-u$ is $3u^2+3u$ and the service increment is $-3u^2+3u$.
The subtraction of $u$ makes the increment from $u=0$ to $u=1$ vanish, so creating the first excess job adds no constant term.
This keeps the bound proportional to the expected number of jobs above the threshold, with coefficients independent of that threshold.
Applying~\eqref{eq:domination} after the cancellation gives
\[
 \cL(P_{3,k}-H_k)
 \le-3(1-\rho\ind_{\{I=0\}})P_{2,k}
        +3(1+\rho\ind_{\{I=0\}})H_k.
\]
Taking stationary expectations gives
\[
 B_k+\frac\Delta N\E[P_{2,k}\ind_{\{I=0\}}]\le2m_k,
\]
which is \eqref{eq:tail-cubic} and implies $B_k\le2m_k$.

Put $C_k=2H_k-Z\ge0$.
The coordinates $I$ and $P_{2,k}$ never change at the same jump, and $\cL P_{2,k}=-C_k$ whenever $I>0$.
Lemma~\ref{mix:identity} therefore gives
\[
 \E[IZ_2P_{2,k}]
 \le(\lambda+\Delta^2/4)B_k
 \le(2N+\Delta^2)m_k.
\]
Pointwise Cauchy--Schwarz gives $H_k^2\le Z_{k+1}P_{2,k}\le Z_2P_{2,k}$, proving \eqref{eq:tail-square}.
Weighted Cauchy--Schwarz and $\E I=\Delta$ then give
\[
 (\E[IH_k])^2\le\E I\,\E[IH_k^2]
       \le\Delta(2N+\Delta^2)m_k.
\]
This proves \eqref{eq:main-weighted-estimate}.
\end{proof}

For fixed $\Delta=\beta$, \eqref{eq:main-weighted-estimate} gives
\begin{equation}
 \frac{\E[IH_k]}N
 \le\left\{\beta\left(2+\frac{\beta^2}{N}\right)
                  \frac{\E H_k}{N}\right\}^{1/2}.
 \label{eq:tail-transfer}
\end{equation}
Uniformity over $k$ therefore yields the implication
\[
 \lim_{k\to\infty}\limsup_{N\to\infty}\E[H_k/N]=0
 \quad\Longrightarrow\quad
 \lim_{k\to\infty}\limsup_{N\to\infty}\E[IH_k/N]=0.
\]
The following bounds give tightness of the scaled stationary populations in all three regimes and the uniform integrability needed to obtain limiting mean congestion and delay.
The mixed moments provide the corresponding control for the idleness-weighted population measures used in NDS and the super-Halfin--Whitt regime.
The idleness moments also bound the error terms in the stationary balance equations and process approximations.
The proof is given in Appendix~\ref{app:population-moments}.

\begin{proposition}[Population and idleness moment bounds]
\label[proposition]{prop:common-moments}
For every $0<B<\infty$, there is a constant $C_B$ such that, under \cref{ass:base}, the following bounds hold for all $N\ge1$ and $0<\Delta<N$ with $\Delta\le B\sqrt N$:
\begin{align}
 \E I&=\Delta,
 &\E I^2&\le\Delta^2+\sqrt{2N}\,\Delta,
 &\E I^3&\le C_B N\Delta, \notag\\
 \E W&\le C_B N/\Delta,
 &\E[IW]&\le C_B N,
 &\E[IW^2]&\le C_B N^2/\Delta, \label{eq:waiting-moments}
\end{align}
and
\begin{equation}
 \E W^2+\E(S-N)^2\le C_B N^2/\Delta^2.
 \label{eq:waiting-second}
\end{equation}
For fixed $\Delta=\beta>0$, one also has
\begin{equation}
 \sup_N\E X_N^2<\infty,
 \qquad \sup_N\E[I_NX_N^2]<\infty,
 \qquad \E I_N^2=O(\sqrt N),
 \label{eq:NDS-moments}
\end{equation}
and
\begin{equation}
 \E[(N-Z_2)P_2]\le C_\beta N^2.
 \label{eq:NDS-profile-moment}
\end{equation}
\end{proposition}

In NDS, these estimates and stationary idleness balance allow us to verify the weighted occupancy approximation from ordinary occupancy convergence and queue-tail control.
We fix $\Delta_N=\beta>0$, take a locally Lipschitz function $A:(1,\infty)\to(0,1]$, and write $a_N=Z_{2,N}/N$.

\begin{proposition}[Criterion for the weighted occupancy approximation]
\label[proposition]{prop:weighted-occupancy-criterion}
Under \cref{ass:base}, Assumption~\ref{ass:weighted-occupancy} holds if
\[
 \begin{aligned}
 &\E\bigl[|a_N-A(X_N)|\ind_{\{X_N\in K\}}\bigr]\longrightarrow0
 \quad\text{for every compact }K\subset(1,\infty),\\
 &\lim_{k\to\infty}\limsup_{N\to\infty}\E[H_{k,N}/N]=0.
 \end{aligned}
\]
\end{proposition}

\begin{proof}
For general spare capacity $\Delta_N$, the idleness birth and death rates are $M_1=N-I-Z_2$ and $\lambda_N\ind_{\{I>0\}}$, so $G_N(I)=I(I+1)/(2N)$ satisfies
\begin{equation}
 \cL G_N
 =1-a_N-Ia_N+\frac{(\Delta_N-1)I-I^2}{N}.
 \label{quad:generator}
\end{equation}
For $\Delta_N=\beta$ and a bounded Lipschitz function $h$ on $\R\times[0,1]$, stationarity and the generator product rule give
\begin{equation}
 \left|\E[h(X_N,a_N)(1-a_N-I_Na_N)]\right|
 \le\frac{C_h}{N}\E(I_N+1)^2=O_h(N^{-1/2}).
 \label{quad:stationary-coefficient}
\end{equation}
Each jump changes the arguments of $h$ by $O(N^{-1})$, so $G_N\cL h$, the generator cross term, and the last term of~\eqref{quad:generator} multiplied by $h$ are bounded by $C_h(I_N+1)^2/N$.
The idleness moment bound in~\eqref{eq:NDS-moments} gives the final order.

Let $\widetilde\mu_N=\Law(X_N,a_N)$ and $\widetilde\eta_N(B)=\E[I_N\ind_{\{(X_N,a_N)\in B\}}]$.
The moment bounds~\eqref{eq:NDS-moments} make these families tight on $\R\times[0,1]$.
For any joint subsequential limit, \eqref{quad:stationary-coefficient} gives $a\widetilde\eta=(1-a)\widetilde\mu$.
Ordinary convergence places $\widetilde\mu$ on $a=A(x)$ for $x>1$, and this identity places $\widetilde\eta$ on the same graph wherever $a>0$.
Fix $x_0>1$, $k\ge2$, and $0<\varepsilon\le(x_0-1)/(2(k-1))$, and let $G_\varepsilon=\{(x,a):x>x_0,\ a<\varepsilon\}$.
On $\{(X_N,a_N)\in G_\varepsilon\}$, one has $H_{k,N}\ge W_N-(k-1)Z_{2,N}\ge(x_0-1)N/2$.
The bound~\eqref{eq:tail-square} therefore gives
\[
 \widetilde\eta_N(G_\varepsilon)
 \le\frac{4\E[I_NH_{k,N}^2]}{(x_0-1)^2N^2}
 \le\frac{4(2+\beta^2/N)}{(x_0-1)^2}\E[H_{k,N}/N].
\]
Portmanteau for the open set $G_\varepsilon$ gives $\widetilde\eta\{a=0,\ x>x_0\}\le\liminf_N\widetilde\eta_N(G_\varepsilon)$ along this subsequence, and letting $k\to\infty$ excludes this remaining weighted mass.
Testing both limiting measures against $\chi(x)|a-A(x)|$, extended by zero outside $(1,\infty)$, with $\chi\in C_c((1,\infty))$ nonnegative and equal to one on $K$, proves~\eqref{path:weighted-coefficient} along every convergent subsequence and hence for the full sequence.
\end{proof}

\subsection{Identifying stationary population laws}
\label{sec:scalar-stationary}

We identify the stationary population laws in NDS and the super-Halfin--Whitt regime from two balance equations.
Population balance relates the ordinary and idleness-weighted measures, while an occupancy relation from idleness balance expresses one in terms of the other.
The Halfin--Whitt limit retains idleness and the number of queues with waiting jobs as separate coordinates.
Its stationary law is the unique invariant law of the reflected diffusion~\eqref{hw:limit-u}--\eqref{hw:limit-v}; the proof of Theorem~\ref{hw:stationary} identifies every stationary path limit with this diffusion.

\paragraph{Identification from balance equations.}
The ordinary measure describes the distribution of population, while its idleness-weighted counterpart records where unused service capacity occurs.
The limiting measures are supported on $[\ell,\infty)$, with $\ell=1$ for the NDS population per server and $\ell=0$ for the scaled excess population in the super-Halfin--Whitt regime.
For $x>\ell$, combining the occupancy relation with population balance determines the population density.
The population balance holds against smooth, compactly supported test functions on $\R$, including those that are nonzero at $x=\ell$, and therefore also constrains any mass assigned to that point by the idleness-weighted measure.
The theorem uses the population density's behavior as $x\downarrow\ell$ to show that this mass is zero.
For a measure $\nu$, write $\nu(h)=\int h\,\dd\nu$; its distributional derivative $D\nu$ is defined by $(D\nu)(h)=-\nu(h')$ for smooth compactly supported test functions $h$.

Fix $\ell\in\R$, a continuous function $b$ on $[\ell,\infty)$, and continuous functions $\phi>0$ and $\psi\ge0$ on $(\ell,\infty)$.
Put $g=\psi/\phi$ and, for a fixed $x_0>\ell$, define
\[
 w(x)=\exp\left\{\int_{x_0}^x\{g(v)-b(v)\}\,\dd v\right\},
 \qquad x>\ell,
 \qquad Z=\int_\ell^\infty w(x)\,\dd x.
\]
When $Z<\infty$, write $p=Z^{-1}w$ on $(\ell,\infty)$ and extend $p$ and $gp$ by zero to $(-\infty,\ell]$.
Let $\mu$ be a probability measure and $\eta$ a finite nonnegative measure on $\R$, both supported on $[\ell,\infty)$, with $\int |b|\,\dd\mu<\infty$.

\begin{theorem}[Identification from stationary balance equations]
\label[theorem]{thm:scalar-identification}
Suppose that
\begin{equation}
 D\mu=\eta-b\mu\quad\text{on }\R,
 \qquad
 \phi\eta=\psi\mu\quad\text{on }(\ell,\infty).
 \label{eq:scalar-two-identities}
\end{equation}
The first identity is distributional, and the second is an equality of locally finite measures.
Then $Z<\infty$, and $\mu$ has density $p$, strictly positive and locally absolutely continuous on $(\ell,\infty)$, with no atom at $\ell$.
If, for some $c>0$,
\begin{equation}
 g(\ell+y)=\frac c y+O(1)
 \qquad\text{as }y\downarrow0,
 \label{eq:scalar-boundary-pole}
\end{equation}
then
\begin{equation}
 \mu(\dd x)=p(x)\ind_{\{x>\ell\}}\,\dd x,
 \qquad
 \eta(\dd x)=g(x)p(x)\ind_{\{x>\ell\}}\,\dd x.
 \label{eq:scalar-identified-measures}
\end{equation}
In this case, neither measure has an atom at $\ell$, $p(\ell+y)\sim C y^c$ for some $C\in(0,\infty)$, and
\begin{equation}
 \eta(\R)=\int_\ell^\infty g(x)p(x)\,\dd x
          =\int_\ell^\infty b(x)p(x)\,\dd x.
 \label{eq:scalar-flow-balance}
\end{equation}
\end{theorem}

\begin{proof}
The signed measure $\sigma=\eta-b\mu$ is finite and is supported on $[\ell,\infty)$.
Let $P(x)=\sigma((\! -\infty,x])$.
This function has bounded variation and distributional derivative $\sigma$.
Hence $D(\mu-P(x)\,\dd x)=0$ on $\R$, so $\mu-P(x)\,\dd x$ is a constant multiple of Lebesgue measure.
Both terms vanish on $(-\infty,\ell)$, forcing that constant to be zero.
Thus $\mu$ has a density $p$ of bounded variation and has no atoms.

On every compact subinterval of $(\ell,\infty)$, the second identity in \eqref{eq:scalar-two-identities} gives $\eta=g\mu$ because $\phi$ is bounded away from zero there.
The first identity therefore implies
\[
 p'=(g-b)p\qquad\text{on }(\ell,\infty)
\]
in the distributional sense.
Its right-hand side is locally integrable, so $p$ has a locally absolutely continuous representative.
An integrating factor then gives $p=C_0w$ throughout this interval.
Since $\mu$ is a probability measure with no atom at $\ell$, $C_0>0$ and $Z=C_0^{-1}<\infty$.

Continuity of $b$ at $\ell$ and \eqref{eq:scalar-boundary-pole} give
\[
 g(\ell+y)-b(\ell+y)=\frac c y+r(y),
 \qquad |r(y)|\le C_1
\]
for sufficiently small $y>0$.
Integrating this identity shows that $w(\ell+y)/y^c$ converges to a finite strictly positive constant.
Consequently $p(\ell+)=0$.

The interior measure identity and the support assumptions leave only the possibility
\[
 \eta=g(x)p(x)\ind_{\{x>\ell\}}\,\dd x+\kappa\delta_\ell,
 \qquad \kappa\ge0.
\]
The extension of $p$ by zero to $(-\infty,\ell)$ has a distributional derivative whose atom at $\ell$ is its jump $p(\ell+)-p(\ell-)=0$.
Since $b\mu$ also has no atom at $\ell$, the identity $D\mu=\eta-b\mu$ on $\R$ forces $\kappa=0$.
This proves \eqref{eq:scalar-identified-measures}.

Finally, a function of bounded variation on $\R$ that is integrable has limit zero at both infinities.
The total mass of its distributional derivative is therefore zero.
Taking total masses in $D\mu=\eta-b\mu$ proves \eqref{eq:scalar-flow-balance}.
\end{proof}

Sections~\ref{sec:nds-proof} and~\ref{sec:shw-proof} derive the NDS and super-Halfin--Whitt balance equations, respectively.

\paragraph{Application to NDS.}
For fixed $\Delta_N=\beta>0$, let $\mu_N=\Law(X_N)$ and $\eta_N(B)=\E[I_N\ind_{\{X_N\in B\}}]$ for Borel sets $B\subset\R$.
Under Assumptions~\ref{ass:base} and~\ref{ass:weighted-occupancy}, their joint subsequential limits satisfy
\[
 D\mu=\eta-\beta\mu\quad\text{on }\R,
 \qquad A\eta=(1-A)\mu\quad\text{on }(1,\infty).
\]
These identities have the form~\eqref{eq:scalar-two-identities} with $\ell=1$, $b=\beta$, $\phi=A$, and $\psi=1-A$.
Lemma~\ref{nds:coefficient-bounds} gives $g(1+y)=1/y+O(1)$ as $y\downarrow0$, verifying~\eqref{eq:scalar-boundary-pole}.
Theorem~\ref{thm:scalar-identification} therefore gives $\mu(\dd x)=p_{g,\beta}(x)\,\dd x$ and $\eta(\dd x)=g(x)p_{g,\beta}(x)\,\dd x$ on $(1,\infty)$, with no atom at one for either measure.

\paragraph{Application to the super-Halfin--Whitt regime.}
For $Y_N=\Delta_N(S_N-N)/N$, let $\nu_N=\Law(Y_N)$ and $\theta_N(B)=\E[(I_N/\Delta_N)\ind_{\{Y_N\in B\}}]$.
Under Assumption~\ref{ass:base} and~\eqref{shw:scaling}, the weighted occupancy collapse and population balance give
\[
 D\nu=\theta-\nu\quad\text{on }\R,
 \qquad y\,\theta(\dd y)=\nu(\dd y)\quad\text{on }(0,\infty)
\]
for every joint subsequential limit.
Here $\ell=0$, $b=1$, $\phi(y)=y$, and $\psi=1$, so $g(y)=1/y$.
Theorem~\ref{thm:scalar-identification} gives the population density $ye^{-y}$ and the idleness-weighted density $e^{-y}$ on $(0,\infty)$, with no atom at zero for either measure.

\subsection{Drift averaging and interchange of limits}
\label{sec:nds-corrector}

To average idleness along the population path, we use the quadratic function $G_N(I)=I(I+1)/(2N)$.
Its generator identity~\eqref{quad:generator} supplies a corrector for stochastic averaging~\cite{Kurtz1992}: it controls the accumulated difference between instantaneous idleness and its limiting conditional mean.
In NDS, we retain the coefficient $A$ from Assumption~\ref{ass:weighted-occupancy} and put $g=(1-A)/A$.
Write $\widehat X_N(t)=S_N(Nt)/N$.
For a compactly supported function $h$ in $(1,\infty)$, products involving $h$ are extended by zero outside this interval.
For a fixed positive coefficient $a$, the leading term of the quadratic generator satisfies $1-a-Ia=-a\{I-(1-a)/a\}$.
Multiplication by $-h(X_N)/A(X_N)$ therefore converts that generator into the centered drift error $h(X_N)\{I_N-g(X_N)\}$ after replacing $a_N$ by $A(X_N)$.
The factor $I_N+1$ in \eqref{path:weighted-coefficient} is precisely the weight produced by this replacement.
Consequently one occupancy estimate controls both stationary drift identification and averaging along the population path.

Lemma~\ref{path:common-corrector} in Appendix~\ref{app:common-technical} makes this argument precise: for every $h\in C_c^1((1,\infty))$ and finite $T$,
\[
 \E\sup_{t\le T}\left|
 \int_0^t h(\widehat X_N(s))
 \{I_N(Ns)-g(\widehat X_N(s))\}\,\dd s
 \right|\longrightarrow0
\]
under stationarity and under initial laws uniformly dominated by stationarity.
Lemma~\ref{quad:maximal} bounds the expected maximum of the normalized corrector $G_N(I_N(c_Nt))/c_N$ on bounded time intervals, where $c_N$ is the time acceleration.
Under stationarity or uniformly dominated initial laws, Proposition~\ref{prop:common-moments} verifies its conditions \eqref{quad:max-smallness} and \eqref{quad:averaging-smallness} for $c_N=N$ in NDS and $c_N=N/\Delta_N^2$ under \eqref{shw:scaling}.
In the Halfin--Whitt regime, idleness remains a diffusion coordinate, and the proof instead uses the reflected queueing equations.

\paragraph{Transfer to nonstationary initial laws.}

A sequence of processes is \emph{$C$-tight} if it is tight and every subsequential limit has continuous paths.

To pass from stationary estimates to process limits, we first transfer path tightness and averaging bounds to the prescribed initial laws, then identify the limit through its martingale problem or reflected queueing equations.
The first step uses domination of the full initial queue configuration.
Under $\zeta_N\le C\Pi_N$, every nonnegative functional of the queueing path has expectation bounded by $C$ times its stationary expectation, because the transition kernels after time zero are unchanged.
This transfers stationary $C$-tightness and averaging estimates with remainders that vanish in $L^1$ to the specified initial laws.

A martingale approximation relative to the full finite-system filtration identifies the process through a well-posed limiting martingale problem.
For each test function, its change along the scaled process, less the time integral of its limiting generator, must equal a martingale plus a remainder whose expected supremum on every bounded time interval tends to zero under stationarity.
Reweighting the initial law preserves these martingales and transfers the $L^1$ remainder bound, so well-posedness of the limiting martingale problem identifies the process from the limiting initial law.
The quadratic corrector supplies this approximation in NDS and the super-Halfin--Whitt regime.
Proposition~\ref{unif:stationary-transfer} in Appendix~\ref{app:initialization} states and proves the general transfer result.
In the Halfin--Whitt regime, domination transfers the vanishing errors in the reflected queueing equations, and continuity of their solution map identifies the process limit.

Once the stationary law and finite-time process limits have been identified, finite-system ergodicity and ergodicity of the limiting diffusion give the two orders of distributional limits.
Uniform integrability of the scaled $p$th moments upgrades the fixed-time and iterated limits to Wasserstein convergence of order $p$.
Full-state domination gives the same moment-tail bound at every time, and the stationary second-moment estimates in Proposition~\ref{prop:common-moments} ensure the required uniform integrability for $1\le p<2$.
In particular, the coordinate means have the same iterated limits.

\section{Proofs of the main results}
\label{sec:main-proofs}

The stationary proofs combine the estimates of Section~\ref{sec:common} with the weighted occupancy approximation in NDS and collapse above two jobs when spare capacity grows.
We then verify the remaining process conditions for interchange and use the full NDS profiles to derive laws for the number of jobs ahead and for customer delay.

\subsection{Stationary limits in NDS}
\label{sec:nds-proof}

We prove Theorem~\ref{thm:regular-nds} by relating the ordinary and idleness-weighted population measures $\mu_N$ and $\eta_N$ from Section~\ref{sec:scalar-stationary}.
For fixed $\Delta=\beta>0$, the population generator under Assumption~\ref{ass:base} yields a basic adjoint relationship between these measures~\cite{BravermanDaiMiyazawa2017}; see also the invariant-measure characterization of Echeverr\'ia~\cite{Echeverria1982}.

\begin{lemma}[Population and idleness balance equations]
\label[lemma]{lem:population-bar}
Suppose Assumption~\ref{ass:base} holds.
Both families are tight, and any joint subsequential limits $\mu,\eta$ have masses $1,\beta$, are supported on $[1,\infty)$, and satisfy
\begin{equation}
 D\mu=\eta-\beta\mu
 \quad\text{in the sense of distributions on }\R.
 \label{eq:whole-line-bar}
\end{equation}
Under Assumption~\ref{ass:weighted-occupancy}, they also satisfy
\begin{equation}
 A\eta=(1-A)\mu\quad\text{on }(1,\infty).
 \label{nds:coefficient-measures}
\end{equation}
\end{lemma}

\begin{proof}
Tightness and preservation of the measure masses follow from \eqref{eq:NDS-moments}.
The bound $X_N\ge1-I_N/N$ and $\E I_N/N\to0$ imply that $\mu$ is supported on $[1,\infty)$.
For $X_N<1-\varepsilon$, necessarily $I_N>\varepsilon N$, so
\[
 \eta_N((\! -\infty,1-\varepsilon))
 \le\frac{\E I_N^2}{\varepsilon N}\longrightarrow0.
\]
The weighted limit has the same support.

Let $f$ be bounded and smooth, with compactly supported derivative.
Taylor expansion of the two population jumps gives
\[
 0=\E\!\left[(I_N-\beta)f'(X_N)
       +\left(1-\frac{\beta+I_N}{2N}\right)f''(X_N)\right]+O_f(N^{-1}).
\]
We pass to the limit and take $f'=h\in C_c^\infty(\R)$ to obtain
\[
 -\mu(h')=\eta(h)-\beta\mu(h).
\]
This is \eqref{eq:whole-line-bar}.

For $h\in C_c^1((1,\infty))$, we apply~\eqref{quad:stationary-coefficient} with $h$ depending only on population.
Under Assumption~\ref{ass:weighted-occupancy}, replacing $a_N$ by $A(X_N)$ produces an error at most
\[
 \|h\|_\infty\E\left[(I_N+1)|a_N-A(X_N)|\ind_{\{X_N\in\operatorname{supp}h\}}\right]\longrightarrow0.
\]
Passing to the limit proves~\eqref{nds:coefficient-measures}.
\end{proof}

\begin{proof}[Proof of Theorem~\ref{thm:regular-nds}]
\emph{Stationary population.}
Lemma~\ref{lem:population-bar} gives tightness and the two balance equations for every joint subsequential limit of $\mu_N,\eta_N$.
Lemma~\ref{nds:coefficient-bounds} gives the occupancy bounds~\eqref{nds:occupancy-bounds} and verifies the boundary condition of Theorem~\ref{thm:scalar-identification}.
That theorem identifies the limiting measures as $\mu(\dd x)=p_{g,\beta}(x)\,\dd x$ and $\eta(\dd x)=g(x)p_{g,\beta}(x)\,\dd x$ on $(1,\infty)$, with no atom at one for either measure.
Uniqueness identifies the full sequence, and the uniform second population moment in Proposition~\ref{prop:common-moments} gives convergence in $W_p$ for every $1\le p<2$.

\emph{Conditional idleness.}
The bound $\E I_N=\beta$ makes the integer-valued idleness sequence tight.
Assumption~\ref{ass:weighted-occupancy}, stationary population convergence, and $\E I_N/N\to0$ give
\[
 \frac{M_{1,N}}N-\{1-A(X_N)\}
 =A(X_N)-\frac{Z_{2,N}}N-\frac{I_N}N\longrightarrow0
 \quad\text{in probability}.
\]
Here $A$ is extended continuously by zero on $(-\infty,1]$, and localization to compact subsets of $(1,\infty)$ is removed using $\Pp(X>1)=1$.
Let $\psi$ be a finitely supported function on $\N_0$ and let $h\in C_c^1((1,\infty))$.
The stationary generator of $\psi(I_N)h(X_N)$, divided by $N$, gives
\[
 \E\!\left[h(X_N)\left\{
 \frac{M_{1,N}}{N}\bigl(\psi(I_N+1)-\psi(I_N)\bigr)
 +\frac{\lambda_N}{N}\ind_{\{I_N>0\}}\bigl(\psi(I_N-1)-\psi(I_N)\bigr)
 \right\}\right]=O_{h,\psi}(N^{-1}).
\]
All population-shift errors have this order because $\psi$ is bounded, the total event rate is at most $2N$, and each population increment is $1/N$.
Passing to a subsequential joint limit and disintegrating over $X=x$ gives the invariant equation for a birth--death chain with birth rate $1-A(x)$ and death rate one away from zero.
Since $A(x)>0$ for $x>1$, its invariant probabilities are $A(x)\{1-A(x)\}^j$, including the degenerate case $A(x)=1$.
This identifies every subsequential joint limit.
The total-mass identity in Theorem~\ref{thm:scalar-identification} gives $\E I=\int g(x)p_{g,\beta}(x)\,\dd x=\beta=\E I_N$.
Weak convergence of nonnegative random variables together with convergence of their first moments gives uniform integrability of $I_N$.

\emph{Mean delay and immediate service.}
Little's law~\cite{Little1961} gives $\E T_N^{\mathrm w}=\E W_N/(N-\beta)$ and $\E T_N^{\mathrm s}=\E S_N/(N-\beta)$.
We use $W_N/N=X_N-1+I_N/N$, $\E I_N=\beta$, and the first-moment convergence established above.
PASTA~\cite{Wolff1982} gives $\Pp(T_N^{\mathrm w}=0)=\Pp(I_N>0)$.
The joint convergence yields the last assertion because the idleness coordinate is discrete.
\end{proof}

\begin{remark}[Idleness second moments]
\label[remark]{rem:idleness-moment}
Under the hypotheses of Theorem~\ref{thm:regular-nds}, the limiting idleness $I$ satisfies $\E I=\beta$ and $\E I^2=\infty$, and $\E I_N^2\to\infty$.
The coefficient bounds in Lemma~\ref{nds:coefficient-bounds} and the boundary conclusion of Theorem~\ref{thm:scalar-identification} give $A(1+y)\sim y$ and $p_{g,\beta}(1+y)\sim cy$ as $y\downarrow0$, for some $c>0$.
The conditional geometric law gives
\[
 \E[I^2\mid X=x]
 =\frac{\{1-A(x)\}\{2-A(x)\}}{A(x)^2}.
\]
Hence $\E[I^2\mid X=1+y]p_{g,\beta}(1+y)\sim2c/y$, whose integral diverges at zero, proving $\E I^2=\infty$.
Applying $I_N\Rightarrow I$ to the bounded continuous tests $i\mapsto i^2\wedge R$ and then letting $R\to\infty$ gives $\E I_N^2\to\infty$.
The product estimate \eqref{eq:main-weighted-estimate} controls the mixed idleness--tail moments directly, without requiring a uniform second moment of idleness.
\end{remark}

\subsection{Verifying the NDS occupancy approximation}
\label{sec:nds-verification}

We verify Assumption~\ref{ass:weighted-occupancy} directly for one-job priority and through ordinary profile convergence and Proposition~\ref{prop:weighted-occupancy-criterion} for fixed fallback routing.
For $\Delta_N=\beta>0$, we compare missing and excess jobs at a queue-length threshold.
For $k\ge1$, define
\[
 V_k=\sum_{i=1}^N(k-Q_i)^+,
 \qquad b_k(q)=\Pp(q_\iota\ge k\mid q).
\]
Then $S=kN-V_k+H_k$, and the exact generators are
\begin{equation}
 \cL H_k=\lambda b_k-Z_{k+1},
 \qquad
 \cL V_k=N-I-Z_{k+1}-\lambda(1-b_k).
 \label{eq:threshold-generators}
\end{equation}
A service completion contributes to the first identity only above level $k$, and to the second only at a nonempty queue of length at most $k$.
An arrival contributes according to whether its chosen queue has at least $k$ jobs.

Exact priority below level $k$ is the condition $V_k>0\Rightarrow b_k=0$.
Idle priority is its $k=1$ instance, and one-job priority adds its $k=2$ instance.
JSQ satisfies the condition at every threshold.
Under exact priority, arrivals cannot increase $H_k$ while $V_k>0$.
Under fallback $F$, the routing probability is $b_k=F(Z_k/N)$ whenever $I=0$.

\paragraph{One-job priority.}

The occupancy error measures the coexistence of missing jobs below level two and excess jobs above it.
The threshold identities give
\begin{equation}
 M_1-N(2-X)^+=\min(V_2,H_2)-2I.
 \label{pri:overlap-identity}
\end{equation}
Exact priority prevents arrivals from increasing $H_2$ while $V_2>0$.
When many queues have at least three jobs, the positive deficit has negative drift; when few such queues exist, the quadratic excess potential controls their contribution to $H_2$.
Lemma~\ref{pri:marked} combines these bounds to give $\E[(I+1)\min(V_2,H_2)]=O_\beta(N^{5/6})$.
Proposition~\ref{pri:weighted-occupancy} uses \eqref{pri:overlap-identity} to convert this overlap bound into the weighted one-job approximation.
The identity $M_1=N-I-Z_2$ then gives \eqref{path:weighted-coefficient} with $A(x)=\min\{x-1,1\}$, proving Example~\ref{ex:priority-nds}.

\begin{remark}[The role of exact routing information]
\label[remark]{rem:threshold-scope}
If arrivals can bypass a queue below level $k$, the positive-deficit drift in \eqref{eq:threshold-generators} contains the additional term $\lambda b_k\ind_{\{V_k>0\}}$.
Such arrivals can also create a long queue during a positive deficit excursion.
Extending this threshold argument requires control of both this drift error and those entries, with the weights used in the averaging estimates.
\end{remark}

\paragraph{Fallback routing and profile attraction.}

For a fixed fallback function $F$ as in Example~\ref{ex:fallback-nds}, we obtain the occupancy coefficient from the mean-preserving profile flow.
The routing probabilities in \eqref{eq:F-routing} determine the occupancy transition rates.
Define the space of finite-mean tail profiles by
\[
\mathcal Q=\left\{q=(q_k)_{k\ge1}:1\ge q_1\ge q_2\ge\cdots\ge0,\quad \sum_{k\ge1}q_k<\infty\right\}.
\]
We set $q_0=1$ and write $m(q)=\sum_{k\ge1}q_k$.
The associated probability distribution on $\N_0$ has masses $p_0=1-q_1$ and $p_k=q_k-q_{k+1}$ for $k\ge1$.
The distance $\norm{q-r}_1$ is the Wasserstein distance of order one between the associated probability distributions.
Under stationarity, the population per server changes by an amount that vanishes in probability on each fixed interval of physical time, while the occupancy profile evolves on that time scale.
Writing $\rho=\lambda_N/N$, the linear idleness balance replaces the factor $\rho\ind_{\{I_N=0\}}$ indicating fallback routing by $u_1^N$ in the integrated profile equations.
Thus the limiting profile evolves at a frozen population level, with arrival coefficient $q_1$.
Arrivals under fallback and service completions give the profile equation
\begin{equation}\label{prof:ode}
\dot q_k=q_1\{F(q_{k-1})-F(q_k)\}-q_k+q_{k+1},\qquad k\ge1.
\end{equation}
Its conserved mean is the excess population per server.
Appendix~\ref{reg:section} justifies this averaging for stationary finite systems.

Recall the profile $q^F(y)$ introduced in Example~\ref{ex:fallback-nds}.
For $y>0$, its defining equations are
\begin{equation}\label{prof:equilibrium}
q^F_1(y)=a_F(y)\in(0,1),\qquad
q^F_{k+1}(y)=a_F(y)F(q^F_k(y)),\qquad
\sum_{k\ge1}q^F_k(y)=y.
\end{equation}
At $y=0$, $q^F(0)=0$.

\begin{theorem}[Attraction of the mean-preserving profile flow]\label[theorem]{thm:profile}
Equation~\eqref{prof:ode} defines a unique global continuous semiflow $S_t$ on $\mathcal Q$ that preserves $m(q)$ and is nonexpansive in $\ell^1$.
For each $y\ge0$, there is a unique equilibrium with mean $y$, namely $q^F(y)$.
Every $q\in\mathcal Q$ satisfies
\begin{equation}\label{prof:attraction}
\lim_{t\to\infty}\norm{S_tq-q^F(m(q))}_1=0.
\end{equation}
The convergence is uniform over compact subsets of $\mathcal Q$ in $\ell^1$.
Moreover, $\norm{q^F(y)-q^F(v)}_1=|y-v|$ for all $y,v\ge0$.
\end{theorem}

The proof uses conservation of the mean and contraction to prove convergence to equilibrium, together with a moment bound that prevents the conserved mean from escaping to high queue lengths.
The difference of two profile trajectories evolves under a conservative Markov operator, which gives strict contraction when their means agree.
A uniform lower bound on $p_0=1-q_1$, the profile probability of zero waiting jobs, gives the moment bound needed for compactness of the profile orbit.
Appendix~\ref{app:profile-flow} develops these arguments and proves Theorem~\ref{thm:profile}.

For stationary finite systems, the ordinary tail bound in Lemma~\ref{reg:tail-tightness} preserves the mean $x-1$ in each joint population--profile limit, while Lemma~\ref{reg:ordinary-profile} makes that limit invariant under the profile flow.
Attraction then places every ordinary limit on the equilibrium graph $u=q^F(x-1)$.
Since $q_1^F(x-1)=a_F(x-1)$, Proposition~\ref{prop:weighted-occupancy-criterion} converts this ordinary occupancy convergence and the tail bound into~\eqref{path:weighted-coefficient} with $A(x)=a_F(x-1)$.
Theorem~\ref{thm:regular-nds} then identifies the population law and conditional idleness and gives uniform integrability of $I_N$.
Proposition~\ref{reg:weighted-law} uses these conclusions to obtain the ordinary and idleness-weighted joint profile limits.

\paragraph{Comparisons of fallback policies.}

Let $\pi_{F,\beta}$ denote the probability law with density $p_{g_F,\beta}$ in \eqref{eq:generic-density}.
Theorem~\ref{thm:regular-nds} identifies this law as the limiting stationary population under the fallback rule $F$.
The equilibrium recursion therefore orders stationary congestion directly through the routing function.
For $a\in[0,1)$, write $v_0^F(a)=1$, $v_{k+1}^F(a)=aF(v_k^F(a))$, and $M_F(a)=\sum_{k\ge1}v_k^F(a)$ for the resulting profile mean.

\begin{proposition}[Fallback information and stationary population order]\label[proposition]{prof:ordering}
For fallback functions $F$ and $G$ as in Example~\ref{ex:fallback-nds}, if $F\le G$, then $g_F\le g_G$, and $\pi_{F,\beta}$ is smaller than $\pi_{G,\beta}$ in the monotone likelihood ratio order, hence in the usual stochastic order $\le_{\mathrm{st}}$.
\end{proposition}

\begin{proof}
If $F\le G$, induction in the equilibrium recursion gives $v_k^F(a)\le v_k^G(a)$, hence $M_F(a)\le M_G(a)$ and $a_F(y)\ge a_G(y)$.
Thus $g_F\le g_G$.
Finally,
\[
\frac{\dd}{\dd x}\log\frac{p_{g_G,\beta}(x)}{p_{g_F,\beta}(x)}=g_G(x)-g_F(x)\ge0,
\]
which proves the stated likelihood ratio order.
Its implication for usual stochastic order is standard; see Shaked and Shanthikumar~\cite[Section~1.C]{ShakedShanthikumar2007}.
\end{proof}

\subsection{Stationary limits in the super-Halfin--Whitt regime}
\label{sec:shw-proof}
We now let spare capacity grow according to \eqref{shw:scaling} and recall $Y_N=\Delta_N(S_N-N)/N$.
We prove Theorem~\ref{thm:superhw} using the occupancy collapse supplied by Assumption~\ref{ass:base} on this scale.

\paragraph{Collapse above two jobs.}

Proposition~\ref{prop:common-moments} and \eqref{eq:tail-cubic} give $\E[P_2\ind_{\{I>0\}}]\le2\E W=O(N/\Delta_N)$.
This estimate first yields ordinary collapse above two jobs, and Theorem~\ref{lem:tail-transfer} then gives the idleness-weighted collapse needed in the population balance.

\begin{lemma}[Ordinary and weighted collapse above two jobs]
\label[lemma]{shw:collapse}
Under Assumption~\ref{ass:base} and \eqref{shw:scaling},
\begin{equation}
\E H_2=o(N/\Delta_N),\qquad
\E[IH_2]=o(N).
\label{shw:collapse-bounds}
\end{equation}
\end{lemma}

\begin{proof}
Suppress the subscript on $\Delta_N$ and put $B=\E[P_2\ind_{\{I>0\}}]$ and $T=\E[IW]$.
Stationarity of $IP_2$, followed by $I\le N$, gives
\begin{equation}
\E[(N-Z_2)P_2]
\le(\lambda+N)B+2T=O(N^2/\Delta).
\label{shw:weighted-p2}
\end{equation}
There are no simultaneous jumps of $I$ and $Z_2$.
When $I>0$, the only changes of $Z_2$ are service completions at queues of length exactly two.
Writing $M_2$ for the number of such queues, stationarity of $IZ_2$ gives the exact identity
\begin{equation}
\lambda\E[Z_2\ind_{\{I=0\}}]
=\E Z_2^2+\E[IZ_2]+\E[IM_2]-\Delta\,\E Z_2
\le\E W^2+2T.
\label{shw:iz2-identity}
\end{equation}
The stationary $H_2$ balance and \cref{ass:base} imply
\begin{equation}
\E Z_3
\le\frac\lambda N\E[Z_2\ind_{\{I=0\}}]
\le\frac{\E W^2+2T}{N}
=O(N/\Delta^2+1).
\label{shw:z3-bound}
\end{equation}

On the event $W\le N/2$, one has $N-Z_2\ge N-W\ge N/2$ deterministically.
Thus \eqref{shw:weighted-p2} gives
\begin{equation}
\E[P_2\ind_{\{W\le N/2\}}]=O(N/\Delta).
\label{shw:local-p2}
\end{equation}
Since $H_2^2\le Z_3P_2$, Cauchy--Schwarz and \eqref{shw:z3-bound}--\eqref{shw:local-p2} yield
\[
\E[H_2\ind_{\{W\le N/2\}}]
=O\left(\frac{N}{\Delta^{3/2}}+\sqrt{\frac N\Delta}\right).
\]
On the complementary event, $H_2\le W\le2W^2/N$, so
\[
\E[H_2\ind_{\{W>N/2\}}]=O(N/\Delta^2).
\]
Both bounds are $o(N/\Delta)$ under \eqref{shw:scaling}.

Applying \eqref{eq:main-weighted-estimate} with $k=2$ gives
\[
(\E[IH_2])^2
\le\Delta(2N+\Delta^2)\E H_2.
\]
The right side is $o(N^2)$, so taking square roots gives $\E[IH_2]=o(N)$.
\end{proof}

\paragraph{Population and idleness-weighted limits.}

The idleness-weighted limit describes the congestion levels at which unused capacity contributes to the stationary balance.
It also prevents the time-integrated idleness drift from concentrating at zero in the process limit.
The population law and its weighted counterpart are identified together.

\begin{proposition}[The stationary law weighted by idleness]
\label[proposition]{shw:weighted-law}
Under Assumption~\ref{ass:base} and \eqref{shw:scaling}, for every bounded continuous $h:\R\to\R$,
\begin{equation}
 \E\left[\frac{I_N}{\Delta_N}h(Y_N)\right]
 \longrightarrow\int_0^\infty h(y)e^{-y}\,\dd y.
 \label{shw:weighted-limit}
\end{equation}
\end{proposition}

\begin{proof}
Recall the population law $\nu_N=\Law(Y_N)$ and the normalized idleness-weighted measure $\theta_N$ from Section~\ref{sec:scalar-stationary}.
The identity $\E I_N=\Delta_N$ gives $\theta_N(\R)=1$.
Proposition~\ref{prop:common-moments} implies
\begin{equation}
\E Y_N^2=O(1),\qquad
\E[(I_N/\Delta_N)Y_N^2]
\le\frac{\Delta_N}{N^2}\{\E[I_NW_N^2]+\E I_N^3\}=O(1).
\label{shw:measure-tightness}
\end{equation}
Hence both sequences of measures are tight.
Their limits are supported on $[0,\infty)$, because $S-N=W-I\ge-I$ gives
\[
\E(Y_N)^-\le\Delta_N^2/N\longrightarrow0,\qquad
\E[(I_N/\Delta_N)(Y_N)^-]\le\E I_N^2/N\longrightarrow0.
\]
We show that every joint subsequential limit $(\nu,\theta)$ satisfies $y\,\theta(\dd y)=\nu(\dd y)$ and $D\nu=\theta-\nu$ on $\R$.
The first identity comes from the quadratic idleness balance after substituting the weighted occupancy collapse; the second comes from the population generator.
Together they give the two measure identities required by Theorem~\ref{thm:scalar-identification}.

The identity $W=Z_2+H_2$ yields
\begin{equation}
\frac{IZ_2}N
=\frac I\Delta Y_N+\frac{I^2}N-\frac{IH_2}N.
\label{shw:collapse-substitution}
\end{equation}
The last two terms vanish in $L^1$ by Proposition~\ref{prop:common-moments} and Lemma~\ref{shw:collapse}.
Fix a bounded twice continuously differentiable $h$ with bounded derivatives.
For the first measure identity, it suffices to show that
\[
\E[h(Y_N)IZ_2/N]-\E h(Y_N)\longrightarrow0.
\]
We use the quadratic idleness corrector $G_N(I)=I(I+1)/(2N)$ and its generator identity \eqref{quad:generator}.
Taylor's formula gives
\[
|\cL h(Y_N)|
\le C_h\left\{\frac{\Delta I}{N}+\frac{\Delta^2}{N}\right\}.
\]
Consequently,
\[
\left|\E[G_N\cL h(Y_N)]\right|
\le C_h\left\{
\frac{\Delta(\E I^3+\E I^2)}{N^2}
+\frac{\Delta^2(\E I^2+\E I)}{N^2}
\right\}
=o(1).
\]
By Proposition~\ref{prop:common-moments}, the leading terms have orders $\Delta^2/N$ and $(\Delta/\sqrt N)^3$, both of which vanish under \eqref{shw:scaling}.
The absolute cross term in the generator product rule is at most
\[
C_h\frac{\Delta(\E I+1)}N
=O_h(\Delta^2/N+\Delta/N)=o(1),
\]
because an idleness-changing jump changes $G_N$ by at most $(I+1)/N$ and $Y_N$ by $\Delta/N$.
The remaining terms in \eqref{quad:generator} satisfy
\[
\frac{\E Z_2}N=O(\Delta^{-1}),\qquad
\frac{(\Delta+1)\E I+\E I^2}{N}
=O(\Delta^2/N+\Delta/\sqrt N)=o(1).
\]
Stationarity of $G_Nh(Y_N)$ therefore gives this limit.

It follows that any subsequential pair $(\nu,\theta)$ satisfies
\begin{equation}
y\,\theta(\dd y)=\nu(\dd y).
\label{shw:limiting-measure-identity}
\end{equation}
The uniform weighted second moment in \eqref{shw:measure-tightness} justifies passage of the factor $y$ to the limit.

For $f\in C_c^3(\R)$, we multiply the stationary generator of $f(Y_N)$ by $N/\Delta_N^2$.
Taylor's formula gives
\[
0=\E\left[
\left(\frac I\Delta-1\right)f'(Y_N)
+\frac{2N-\Delta-I}{2N}f''(Y_N)
\right]+O_f(\Delta/N).
\]
The expected difference of the second-order coefficient from one is $\Delta/N$.
Thus
\begin{equation}
\int f''\,\dd\nu+\int f'\,\dd(\theta-\nu)=0.
\label{shw:limiting-bar}
\end{equation}
Equation~\eqref{shw:limiting-bar} gives $D(D\nu-\theta+\nu)=0$, so $D\nu-\theta+\nu$ is a constant distribution.
It vanishes on $(-\infty,0)$ and therefore equals zero.
Thus $D\nu=\theta-\nu$ on $\R$.

The preceding tightness and the two limiting measure identities verify the super-Halfin--Whitt application of Theorem~\ref{thm:scalar-identification} in Section~\ref{sec:scalar-stationary}.
The identified density is $p(y)=ye^{-y}$, and the weighted measure is $\theta(\dd y)=e^{-y}\,\dd y$ on $(0,\infty)$, with no atom at $y=0$.
These limits are unique, proving \eqref{shw:weighted-limit}.
\end{proof}

\begin{proof}[Proof of Theorem~\ref{thm:superhw}]
The measure identity~\eqref{shw:limiting-measure-identity} and Proposition~\ref{shw:weighted-law} identify every subsequential population limit as the law of $G$, so $Y_N\Rightarrow G$.
The uniform second moment in \eqref{shw:measure-tightness} gives convergence in $W_p$ for every $1\le p<2$, and hence $\E Y_N\to2$.
Finally,
\[
\E\left|\frac{\Delta W}N-Y_N\right|=\frac{\Delta^2}N\to0,
\qquad
\E\left|\frac{\Delta Z_2}N-Y_N\right|
\le\frac{\Delta^2}N+\frac{\Delta\E H_2}N\to0.
\]
These give the remaining weak limits, which strengthen to $W_p$ convergence for $1\le p<2$ by Proposition~\ref{prop:common-moments} and $Z_2\le W$.
\end{proof}

\subsection{Stationary limits in the Halfin--Whitt regime}
\label{sec:hw-proof}

We prove Theorem~\ref{hw:stationary} by identifying the stationary path limits of $(U_N,V_N)=(-I_N,Z_{2,N})/\sqrt N$ when $\Delta_N=\beta\sqrt N$.
Idleness remains visible on this scale, so the limit is governed by the reflected equations \eqref{hw:limit-u}--\eqref{hw:limit-v}.
The following lemma removes jobs beyond the second position.

\begin{lemma}[Higher queues in the Halfin--Whitt regime]
\label[lemma]{hw:collapse}
Under Assumption~\ref{ass:base} and $\Delta_N=\beta\sqrt N$ with fixed $\beta>0$,
\[
 \E Z_{3,N}=O(1),\qquad \E H_{2,N}=O(N^{1/4}).
\]
\end{lemma}

\begin{proof}
Proposition~\ref{prop:common-moments} and Theorem~\ref{lem:tail-transfer}, with $\Delta=\beta\sqrt N$, give
\begin{equation}
 \E I^2+\E W^2=O_\beta(N),\qquad
 \E[IW]=O_\beta(N),\qquad
 \E[P_2\ind_{\{I>0\}}]=O_\beta(\sqrt N).
 \label{hw:moments}
\end{equation}
The argument for ordinary collapse in Lemma~\ref{shw:collapse} remains valid with $\Delta=\beta\sqrt N$.
In detail, \eqref{shw:iz2-identity} gives
\[
 \lambda\E[Z_2\ind_{\{I=0\}}]\le\E W^2+2\E[IW]=O(N),
 \qquad \E Z_3=O(1).
\]
Also, $\E[(N-Z_2)P_2]=O(N^{3/2})$, so $\E[P_2\ind_{\{W\le N/2\}}]=O(\sqrt N)$.
It follows that
\[
 \E[H_2\ind_{\{W\le N/2\}}]\le
 \sqrt{\E Z_3\,\E[P_2\ind_{\{W\le N/2\}}]}=O(N^{1/4}),
\]
whereas $\E[H_2\ind_{\{W>N/2\}}]\le2\E W^2/N=O(1)$.
\end{proof}

\begin{proof}[Proof of Theorem~\ref{hw:stationary}]
We use Lemma~\ref{hw:collapse} and \eqref{hw:moments} to identify each stationary path limit.
Its one-time law is invariant for the reflected diffusion, whose unique invariant law $\Pi_\beta$ then determines the stationary limit.
Let $A$ count all arrivals, let $A_0$ count arrivals that see $I=0$, and let $A_{\ge2}$ count arrivals routed to queues with at least two jobs.
Let $D_j$ count service completions from queues containing exactly $j$ jobs immediately before completion.
The routing assumptions imply
\[
 \E A_{\ge2}(T)
 \le\frac\lambda N\,T\E[Z_2\ind_{\{I=0\}}]=O_T(1).
\]
Put $R_N=A_0/\sqrt N$.
The exact queueing equations are
\begin{align*}
 U_N(t)&=U_N(0)+M_N(t)-\beta t+
       \int_0^t\{-U_N(s)+V_N(s)\}\,\dd s-R_N(t),\\
 V_N(t)&=V_N(0)-\int_0^tV_N(s)\,\dd s+R_N(t)+\varepsilon_N(t),
\end{align*}
where
\[
 M_N(t)=\frac{A(t)-\lambda t}{\sqrt N}
 -\frac{D_1(t)-\int_0^tM_1(s)\,\dd s}{\sqrt N}
\]
and
\[
 \varepsilon_N(t)
 =-\frac{A_{\ge2}(t)}{\sqrt N}
 -\frac{D_2(t)-\int_0^tM_2(s)\,\dd s}{\sqrt N}
 +\int_0^t\frac{Z_3(s)}{\sqrt N}\,\dd s.
\]
Stationarity, \eqref{hw:moments}, and Doob's inequality give
\[
 \E\sup_{t\le T}|\varepsilon_N(t)|\longrightarrow0.
\]
The martingale $M_N$ has jumps of size $1/\sqrt N$ and predictable quadratic variation
\[
 \langle M_N\rangle(t)=\int_0^t\frac{\lambda+M_1(s)}N\,\dd s
 \longrightarrow2t
\]
uniformly in probability on bounded intervals, since $M_1=N-I-Z_2$ and $\E(I+Z_2)=O(\sqrt N)$.
The martingale functional central limit theorem therefore gives $M_N\Rightarrow\sqrt2 B$, jointly with any convergent subsequence of initial states; the Brownian motion is independent of the limiting initial state~\cite[Theorem~2.1]{Whitt2007}.

The process $R_N$ is exactly the upper Skorokhod regulator for the first equation.
Indeed, it increases by $1/\sqrt N$ only when an arrival sees $U_N=0$, cancelling that arrival's positive increment in the unconstrained first coordinate.
The equations fit Lemma~1 of Eschenfeldt and Gamarnik~\cite{EschenfeldtGamarnik2018} with upper second-coordinate boundary $B=\infty$, initial vector $(U_N(0),V_N(0))$, and input paths $(M_N(t)-\beta t,\varepsilon_N(t))$.
Its continuity conclusion applies because the inputs converge to continuous paths; the second regulator is identically zero.
Consequently every subsequential stationary path limit solves \eqref{hw:limit-u}--\eqref{hw:limit-v}.
Its one-time law is invariant, and uniqueness of $\Pi_\beta$ identifies the limit.
The uniform second moments in \eqref{hw:moments} give convergence in $W_p$ for every $1\le p<2$.
Finally, $S-N=-I+Z_2+H_2$ and $\E H_2/\sqrt N\to0$ give the population conclusion.
\end{proof}

\subsection{Process limits and interchange of limits}
\label{sec:process-proof}

For the NDS and super-Halfin--Whitt limits in Theorem~\ref{thm:common-interchange}, we apply Proposition~\ref{unif:stationary-transfer} using stationary $C$-tightness, full-filtration martingale approximations, and well-posedness of the limiting diffusions.
The quadratic corrector controls their accumulated idleness, while the Halfin--Whitt argument transfers the reflected queueing equations through their continuous solution map.

\begin{proof}[Proof of Theorem~\ref{thm:common-interchange}]
For NDS, Appendix~\ref{app:stationary-process} uses Lemma~\ref{path:common-corrector} and localization away from $x=1$ to establish stationary path convergence and the full-filtration martingale approximation with an $L^1$ remainder.
The coefficient bounds in Lemma~\ref{nds:coefficient-bounds} and the local Lipschitz continuity of $A$ verify the hypotheses of Lemma~\ref{int:ergodicity}, which gives a well-posed martingale problem on $(1,\infty)$ and convergence in total variation to $\pi_{g,\beta}$.

For the super-Halfin--Whitt window, Appendix~\ref{app:shw-process} proves stationary path convergence and the martingale approximation with clock $c_N=N/\Delta_N^2$.
The moment substitutions in Appendix~\ref{app:quadratic-corrector} verify the conditions of Lemma~\ref{quad:maximal} on this clock.
The quadratic corrector identifies $\int_0^t\widehat Y_N(s)\,\dd A_N(s)\to t$, where $A_N(t)=\int_0^tI_N(Ns/\Delta_N^2)/\Delta_N\,\dd s$.
On $[0,T]$, Proposition~\ref{shw:weighted-law} makes the expected accumulated idleness drift below population level $\varepsilon$ asymptotically at most $T(1-e^{-\varepsilon})$, which tends to zero as $\varepsilon\downarrow0$.
This excludes concentration of the limiting drift at zero, while $Y\,\dd A=\dd t$ gives the reciprocal term $\dd A=\dd t/Y$ in the positive interior.
Lemma~\ref{int:ergodicity}, applied with $g(x)=1/(x-1)$, $\beta=1$, and the translation $X=1+Y$, gives a well-posed martingale problem on $(0,\infty)$ and convergence in total variation to $\pi_{\mathrm{shw}}$.

In the Halfin--Whitt regime, we use the reflected queueing equations from the proof of Theorem~\ref{hw:stationary}, driven by $M_N$ with error $\varepsilon_N$.
Under \eqref{eq:common-dominated-initialization}, full-path domination transfers $\E\sup_{t\le T}|\varepsilon_N(t)|\to0$ and the uniform convergence in probability of $\langle M_N\rangle(t)$ to $2t$ on bounded intervals.
The time-zero change of measure preserves the martingale property and predictable quadratic variation of $M_N$, whose jumps have size $N^{-1/2}$.
The joint martingale functional central limit theorem therefore gives $M_N\Rightarrow\sqrt2 B$, with $B$ independent of the limiting initial state.
The regulator constraints hold pathwise, so continuity of the reflected solution map~\cite[Lemma~1]{EschenfeldtGamarnik2018} identifies the limit with initial law $\nu$.
Exponential ergodicity is Corollary~1 of Braverman~\cite{Braverman2020}.

Proposition~\ref{unif:stationary-transfer} therefore gives process convergence in NDS and the super-Halfin--Whitt regime, and the reflected solution map gives it in the Halfin--Whitt regime.
Finite-system ergodicity follows from Lemma~\ref{lem:foundations} and the aperiodic continuous-time chain; its integer-time skeleton has a positive holding probability~\cite[Theorem~13.0.1]{MeynTweedie1993Book}.
For the moment conclusions, the stationary second-moment bounds are Proposition~\ref{prop:common-moments}, \eqref{shw:measure-tightness}, and \eqref{hw:moments}.
For the stationary initial observation $\mathcal X_N(0)$ in each regime and $1\le p<2$,
\[
 \E\!\left[\|\mathcal X_N(0)\|^p
 \ind_{\{\|\mathcal X_N(0)\|>R\}}\right]
 \le R^{p-2}\E\|\mathcal X_N(0)\|^2.
\]
This proves uniform integrability of the $p$th moments and hence stationary $W_p$ convergence~\cite[Definition~6.8 and Theorem~6.9]{Villani2009}.
Full-state domination gives the same bound uniformly over time under the specified initial laws.
The uniform-integrability and ergodicity argument in the proof of Proposition~\ref{unif:stationary-transfer} therefore gives fixed-time and iterated $W_p$ limits and convergence of the corresponding coordinate means in all three regimes.
Convergence for continuous functions bounded in absolute value by $C'(1+\|x\|^a)$, for some $C'<\infty$ and $0<a<2$, follows from the same uniform-integrability argument.
\end{proof}

\subsection{From occupancy profiles to customer delay in NDS}
\label{sec:delay-proof}

We return to NDS, with $\lambda_N=N-\beta$ and fixed $\beta>0$.
Proposition~\ref{reg:weighted-law} already identifies the full stationary profile under fallback $F$.
For one-job priority, Theorem~\ref{thm:regular-nds} identifies population congestion through the occupancy condition~\eqref{path:weighted-coefficient}.
Customer-delay distributions require the ordinary occupancy profile at all population levels.
We obtain these profiles for JSQ and fixed-threshold policies with a specified fallback function, then use conditional level balance to derive laws for the number of jobs ahead and for delay in all three policy families.

\paragraph{JSQ and threshold-priority profiles.}

Fix the threshold-priority rule from Section~\ref{sec:delay-distribution}, with integer $r\ge1$ and fallback function $F$ as in Example~\ref{ex:fallback-nds}, both fixed across $N$.
Retain the JSQ and threshold-priority profiles $q^{\mathrm{JSQ}}(x)$ and $q^{r,F}(x)$ defined in \eqref{ext:jsq-profile} and \eqref{ext:priority-profile}.
Both profiles have mean $x-1$.
The JSQ profiles are coordinatewise ordered by $x$, so their $\ell^1$ distance equals the difference between their population parameters.
The threshold profiles are continuous in $\ell^1$ by their finite lower coordinates and the Lipschitz identity in Theorem~\ref{thm:profile}.
Their first $r$ coordinates record the queues directly controlled by shortest-queue priority; the remaining coordinates follow the profile associated with $F$ for the excess population above level $r+1$.

Recall $u_m^N=Z_{m+1}/N$, and let $X_*$ have density \eqref{eq:priority-density}.

\begin{proposition}[Full profiles under threshold priority]
\label[proposition]{ext:profiles}
Fix $\beta>0$ and let $\lambda_N=N-\beta$.
For JSQ,
\[
 (X_N,u^N)\Rightarrow(X_*,q^{\mathrm{JSQ}}(X_*))
 \quad\text{in }\R\times\ell^1.
\]
For threshold priority with fixed $r\ge1$ and fixed fallback function $F$ as in Example~\ref{ex:fallback-nds},
\[
 (X_N,u^N)\Rightarrow(X_*,q^{r,F}(X_*))
 \quad\text{in }\R\times\ell^1.
\]
\end{proposition}

The proof in Appendix~\ref{app:threshold-profiles} applies Lemma~\ref{pri:marked} at every exact-priority threshold $j\ge2$, while the bound at $j=1$ follows from $\min(V_1,H_1)\le I$ and $\E I=\beta$.
It gives $\E\min(V_j,H_j)/N=O(N^{-1/3})$ whenever routing gives exact priority below $j$.
This determines the lower profile coordinates, while Lemmas~\ref{reg:regeneration} and~\ref{reg:ordinary-profile}, applied with cutoff $\ell=r+1$, identify the upper profile through the attraction in Theorem~\ref{thm:profile}.
The distinction between the two parts reflects which queue lengths the dispatcher compares exactly.

\paragraph{Jobs ahead and Erlang mixtures.}

Routing determines the queue an arrival joins, so we derive the law of the number of jobs ahead through conditional level balance.
FCFS and memorylessness then give the Erlang mixtures.
Let $D_N$ be the number of jobs already at the server selected by a stationary arrival, including any job in service.
Let $T_N^{\mathrm w}$ and $T_N^{\mathrm s}$ be its waiting and sojourn times.
Use the service-time sums $E_n$, $n\in\N_0$, introduced in Section~\ref{sec:delay-distribution}.

Let $q(x)=(q_m(x))_{m\ge1}$, $x\ge1$, be a tail profile that is continuous as an $\ell^1$-valued function and satisfies $\sum_{m\ge1}q_m(x)=x-1$.
Given a population variable $X\ge1$, define $D$, $T^{\mathrm w}$, and $T^{\mathrm s}$ through
\[
\Pp(D\ge m\mid X=x)=q_m(x),\quad m\ge1,
\]
and
\[
\Law(T^{\mathrm w}\mid X,D=d)=\Law(E_d),
\qquad
\Law(T^{\mathrm s}\mid X,D=d)=\Law(E_{d+1}),
\qquad d\in\N_0.
\]

\begin{lemma}[Arrival-level balance]
\label[lemma]{ext:arrival-balance}
Suppose that Assumption~\ref{ass:base} holds and that
\[
(X_N,u^N)\Rightarrow(X,q(X))
\quad\text{in }\R\times\ell^1,
\qquad \E X_N\to\E X.
\]
Then $(X_N,D_N)\Rightarrow(X,D)$.
The laws of $D_N$, $T_N^{\mathrm w}$, and $T_N^{\mathrm s}$ also converge in Wasserstein distance of order one to those of $D$, $T^{\mathrm w}$, and $T^{\mathrm s}$, respectively.
\end{lemma}

\begin{proof}
PASTA identifies the arrival's pre-routing state with the stationary state; its nonanticipation condition holds because routing uses the current queue lengths and independent arrival marks \cite{Wolff1982}.
Put $b_m(Q)=\Pp(Q_\iota\ge m\mid Q)$.
The exact level balance is $\cL H_m=\lambda b_m-Z_{m+1}$.
Multiplying this balance by a smooth function of population retains the information needed to identify the conditional tails of $D_N$.
For $h\in C_c^2(\R)$, the prelimit calculation in \eqref{nds:coefficient-level-balance} uses only Assumption~\ref{ass:base} and Proposition~\ref{prop:common-moments}.
Since $\rho=1-\beta/N$, that calculation gives
\begin{equation}
\left|\E[h(X_N)b_m(Q)]-\E[h(X_N)u_m^N]\right|=O_h(N^{-1}).
\label{ext:conditional-level-balance}
\end{equation}
The routing bound~\eqref{eq:domination} implies $\E D_N\le\E X_N$, so $(X_N,D_N)$ is tight.
Passing to a subsequential limit in \eqref{ext:conditional-level-balance} identifies every conditional tail and hence the joint law.

Under FCFS, the jobs ahead of the arrival require $D_N$ successive rate-one exponential completions; memorylessness applies to the job already in service.
Later arrivals do not affect this waiting time.
Conditional on $X_N$ and $D_N=d$, the waiting and sojourn times therefore have laws $\Law(E_d)$ and $\Law(E_{d+1})$ in every finite system.
Little's law~\cite{Little1961} and $\E I_N=\beta$ give
\[
\E D_N=\E T_N^{\mathrm w}
=\frac{\E S_N-N+\E I_N}{N-\beta}
=\frac{\E S_N}{N-\beta}-1
\longrightarrow\E(X-1)=\E D.
\]
Convergence of the mixing laws gives weak convergence of the Erlang mixtures, whose means are $\E D_N$ and $\E D_N+1$.
The first-moment criterion for Wasserstein convergence gives all three conclusions.
\end{proof}

Proposition~\ref{reg:weighted-law} verifies the lemma's hypotheses for fallback $F$ with $q_m(x)=q_m^F(x-1)$; Proposition~\ref{ext:profiles} supplies the two priority profiles.
For each of these classes, $\E[D\mid X=x]=\E[T^{\mathrm w}\mid X=x]=x-1$.

\paragraph{Delay comparisons.}

Let $X_F$ denote the limiting population under the policy with fallback $F$ in Example~\ref{ex:fallback-nds}, identified in Theorem~\ref{thm:regular-nds}.

\begin{proof}[Proof of Proposition~\ref{ext:delay-order}]
Fix the stop-loss function $\phi_t(w)=(w-t)^+$ and put $f_t(n)=\E\phi_t(E_n)$ for $n\in\N_0$.
The sequence $f_t$ is nondecreasing and convex, the random-sum convexity used in Shaked and Shanthikumar~\cite[Theorem~3.A.13]{ShakedShanthikumar2007}.
For an independent rate-one exponential $E$, convexity makes $w\mapsto\E[\phi_t(w+E)-\phi_t(w)]$ nondecreasing, and the natural coupling gives $E_{n+1}\ge E_n$.
Thus $f_t(n+1)-f_t(n)$ is nondecreasing in $n$.
Let $\bar f_t$ be its piecewise-linear interpolation.
For every nonnegative integer-valued $D$ with mean $y$,
\[
 \E f_t(D)=\E\bar f_t(D)\ge\bar f_t(y).
\]
Under JSQ, the conditional law of $D$ given $X=x$ is supported on at most two adjacent integers and attains equality for $y=x-1$.
All the priority policies have $\E[D\mid X=x]=x-1$ and share the population law of JSQ.
Conditioning on that population proves JSQ's convex-order minimality for waiting time.

Let $G_y$ have tails $\Pp(G_y\ge k)=q_k^F(y)$, so $\E G_y=y$, with $G_0=0$ and $a_F(0)=0$.
For $x>r+2$, put $y=x-r-1$.
Conditional on $X=x$, the numbers of jobs ahead under thresholds $r$ and $r+1$ are distributed as $r+G_y$ and $r+1+G_{y-1}$, respectively.
It therefore suffices in this case to prove $1+G_{y-1}\le_{\mathrm{cx}}G_y$ for $y\ge1$, where both variables have mean $y$.
For $y\ge1$, put $a=a_F(y)$, $b=a_F(y-1)$, and
\[
 d_i=q_{i+1}^F(y)-q_i^F(y-1),\qquad i\ge1.
\]
Since $a_F$ is increasing, $a\ge b$.
The profile recursion implies that $d_i\ge0$ entails $d_{i+1}\ge0$.
Moreover, the differences are absolutely summable and
\[
 \sum_{i\ge1}d_i=(y-a)-(y-1)=1-a>0.
\]
Hence the sequence has at most one change from negative to nonnegative terms, and every tail sum is nonnegative.
For each integer $k\ge1$, this gives
\[
 \E(G_y-k)^+-\E(1+G_{y-1}-k)^+
 =\sum_{i\ge k}d_i\ge0.
\]
Their means agree, and stop-loss expectations interpolate linearly between integer thresholds, so $1+G_{y-1}\le_{\mathrm{cx}}G_y$.

This establishes the comparison for $x>r+2$.
If $x\le r+2$, the conditional law of $D$ under threshold $r+1$ is supported on at most two adjacent integers and is again minimal among integer-valued distributions with that mean.
Applying the convex function $f_t$ to these conditional comparisons and then integrating over $X_*$ proves the waiting-time order between consecutive thresholds.
The equal means convert the stop-loss inequalities into convex order.
Replacing $f_t(n)$ by $\E\phi_t(E_{n+1})$ proves all the sojourn-time comparisons in the same way.

For the limit as $r\to\infty$, we couple the numbers of jobs ahead under threshold priority and JSQ to agree whenever $X_*\le r+1$, where the profiles already coincide.
On the complementary event, their conditional expected sum is $2(X_*-1)$.
Using the same sequence of service exponentials and writing $T$ for either waiting or sojourn time gives
\[
 W_1\bigl(\Law(T_{r,F}),\Law(T_{\mathrm{JSQ}})\bigr)
 \le2\E[(X_*-1)\ind_{\{X_*>r+1\}}]\longrightarrow0.
\]

Finally, \eqref{nds:idleness-bounds} gives
\[
 \frac{\dd}{\dd x}\log\frac{p_{g_F,\beta}(x)}{p_{g_*,\beta}(x)}
 =g_F(x)-g_*(x)\ge0.
\]
Thus $X_*\le_{\mathrm{lr}}X_F$, and hence $X_*\le_{\mathrm{st}}X_F$.
For $1<x\le2$, the conditional law of $D$ under one-job priority is supported on at most two adjacent integers with mean $x-1$ and is therefore smaller than $G_{x-1}$ in convex order.
For $x>2$, its law is that of $1+G_{x-2}$, which is smaller than $G_{x-1}$ in convex order by the comparison proved above.
Set $\psi_t(x)=\E f_t(G_{x-1})$ for $x>1$.
The equilibrium profiles are coordinatewise nondecreasing in their mean by Lemma~\ref{prof:massmap}, so $G_y$ is stochastically nondecreasing in $y$ and $\psi_t$ is nondecreasing.
The bound $f_t(n)\le n+(-t)^+$ and the finite population means make the following expectations finite:
\[
 \E\phi_t(T_{1,F}^{\mathrm w})
 \le\E\psi_t(X_*)
 \le\E\psi_t(X_F)
 =\E\phi_t(T_F^{\mathrm w}).
\]
This proves increasing convex order for waiting time.
Replacing $f_t(n)$ by $\E\phi_t(E_{n+1})$ gives the sojourn-time comparison.
\end{proof}

\paragraph{Growing priority thresholds.}

The preceding limit sends $r\to\infty$ after the stationary many-server limit.
Corollary~\ref{ext:growing-thresholds} allows the threshold to grow with $N$.
The estimate in Lemma~\ref{pri:marked}, whose constant is independent of the threshold, controls the accumulated error in the first profile coordinates, while the population moment bound controls the remaining tail.

\begin{proof}[Proof of Corollary~\ref{ext:growing-thresholds}]
Every policy in the sequence gives one-job priority, so Proposition~\ref{pri:weighted-occupancy} and Theorem~\ref{thm:regular-nds} give convergence of $X_N$ to $X_*$ in Wasserstein distance of order one, while \eqref{eq:NDS-moments} bounds the second moments uniformly.
For every $2\le j\le r_N+1$, exact priority below level $j$ holds, so Lemma~\ref{pri:marked} gives $\E\min(V_j,H_j)\le C_\beta N^{2/3}$ uniformly in the threshold and fallback.
For $j=1$, the same normalized bound follows from $\min(V_1,H_1)\le I$ and $\E I=\beta$.
Write $\varepsilon_j=\min(V_j,H_j)/N$.
Then
\[
 \E\varepsilon_j\le C_\beta N^{-1/3},
 \qquad
 \frac{H_j}N=(X_N-j)^++\varepsilon_j,
 \qquad j\le r_N+1.
\]
Since $u_m^N=(H_m-H_{m+1})/N$, for every integer $1\le J\le r_N$,
\begin{align*}
 \E\|u^N-q^{\mathrm{JSQ}}(X_N)\|_1
 &\le 2\sum_{j=1}^{J+1}\E\varepsilon_j
       +\E\varepsilon_{J+1}
       +2\E(X_N-J-1)^+\\
 &\le C_\beta\left\{(J+1)N^{-1/3}+\frac1{J+1}\right\}.
\end{align*}
The last step uses the uniform second population moment.
For all sufficiently large $N$, choosing $J=\min\{r_N,\lfloor N^{1/6}\rfloor\}$ gives the bound
\begin{equation}
 \E\|u^N-q^{\mathrm{JSQ}}(X_N)\|_1
 \le C_\beta\left(N^{-1/6}+\frac1{r_N}\right).
 \label{ext:growing-threshold-profile-bound}
\end{equation}
The constant $C_\beta$ in \eqref{ext:growing-threshold-profile-bound} is independent of the thresholds, fallbacks, and system size.
The comparison evaluates the JSQ profile at the same finite-system population $X_N$; convergence of the population law follows separately from Theorem~\ref{thm:regular-nds}.
The map $x\mapsto q^{\mathrm{JSQ}}(x)$ is one-Lipschitz in $\ell^1$ on the extended domain.
Coupling the populations therefore proves joint $W_1$ convergence; the bound $|X_N|+\|u^N\|_1\le2X_N$ and \eqref{eq:NDS-moments} strengthen it to $W_p$ convergence for every $1\le p<2$.
Lemma~\ref{ext:arrival-balance} then gives convergence of $D_N$ and both customer-delay laws.
\end{proof}

\section{Conclusion}
\label{sec:conclusion}

Which queue-length information matters for stationary performance depends on spare capacity and on the performance measure.
In the super-Halfin--Whitt and Halfin--Whitt regimes, idle-server information suffices: JIQ and every other policy in the basic routing class have the JSQ scaled stationary limits.
With fixed spare capacity, JIQ has a more congested population limit, while one-job priority recovers the JSQ population law and mean delays.
For a fixed fallback rule in the family studied here, each additional shortest-queue priority level beyond one job reduces expected delay in excess of any target in the limit, while preserving the limiting population law and mean delays.

In NDS, stationary idleness balance and the uniform tail estimates turn ordinary occupancy convergence and ordinary queue-tail control into the weighted approximation required by the population drift.
The same approximation identifies the stationary population law and justifies dynamic averaging.

Extending the process limits from the uniformly dominated initial laws treated in Theorem~\ref{thm:common-interchange} to deterministic initial configurations requires direct control of the initial occupancy evolution.
For the NDS delay laws identified here, stronger tail estimates would extend convergence of means and expected excess delays to variances and higher moments.

Two further extensions concern NDS sampling: ordinary shortest-of-sample routing, which can miss idle servers, and idle-first shortest-of-sample routing with a growing number of independent uniform samples.
The former requires estimates that account for the configurations in which shorter queues are bypassed; the latter calls for ordinary profile convergence and tail control uniform over the changing sampling rule.

General service distributions and heterogeneous servers require additional state variables describing available capacity; the vacation models of Atar, Karmakar, and Lipshutz~\cite{AtarKarmakarLipshutz2022} illustrate how availability and congestion can remain coupled when their time scales match.
When availability equilibrates faster, the averaging problem is to identify its conditional mean given congestion and control the weighted error in its contribution to the population drift.

\clearpage
\appendix

\section{Stationary bounds and idleness averaging}
\label{app:common-technical}

Finite-system stability and integrability justify the stationary generator identities used in Section~\ref{sec:common}.
The moment and occupancy bounds support stationary identification, while the idleness corrector estimates control the accumulated drift errors in Section~\ref{sec:nds-corrector}.

\subsection{Stability and stationary generator identities}

\begin{proof}[Proof of Lemma~\ref{lem:foundations}]
The total transition rate is at most $2N$, so the chain is nonexplosive.
A service-only path reaches the empty state from every finite state with positive probability.
The routing bound~\eqref{eq:domination} gives positive probability of choosing some shortest queue at every arrival.
An arrival-only path choosing shortest queues reaches a prescribed common height at all servers, after which a service-only path reaches any prescribed state.
Thus the chain is irreducible.

For $1<z<N/\lambda$, put $\Phi_z(q)=\sum_i z^{q_i}$.
The routing bound~\eqref{eq:domination} and the service transitions give
\[
 \cL\Phi_z
 \le (z-1)\left(\frac\lambda N-\frac1z\right)\Phi_z
       +N(1-z^{-1}).
\]
The coefficient of $\Phi_z$ is negative and its sublevel sets are finite.
The Foster--Lyapunov criterion gives positive recurrence and a finite stationary $\Phi_z$ moment~\cite[Theorem~4.2]{MeynTweedie1993III}.
Consequently all polynomial moments of $S$ are finite for fixed $N$.
For a polynomially growing nonnegative test $f$, the increment of $f\wedge R$ is bounded in absolute value by the increment of $f$.
The stationary absolute jump rate of that increment is integrable, so dominated convergence gives $\E\cL f=0$.
Differences of nonnegative polynomial tests are treated separately.
\end{proof}

\subsection{Population and idleness moment bounds}
\label{app:population-moments}

\begin{proof}[Proof of Proposition~\ref{prop:common-moments}]
The bounded test $I^3-I$ has upward increment $3I(I+1)$ and downward increment $-3I(I-1)$.
Using $M_1\le N-I$ gives
\[
 \E I^3\le\Delta\E I^2+2N\Delta.
\]
Since $(\E I^2)^2\le\E I\,\E I^3$, solving the resulting quadratic proves the second-moment bound.
The third-moment bound follows from $\Delta^2\le B^2N$.

Put $M=\E W$ and $T=\E[IW]$.
Equation~\eqref{eq:waiting-balance} gives $\Delta M\le T+N+\Delta^2$.
Applying \cref{lem:tail-transfer} with $k=1$ gives
\[
 T^2\le\Delta(2N+\Delta^2)M.
\]
Because $\Delta^2\le B^2N$, substitution of \eqref{eq:waiting-balance} yields
\[
 T^2\le(2N+\Delta^2)(T+N+\Delta^2)
       \le(2+B^2)N\{T+(1+B^2)N\}.
\]
The positive root is bounded by a constant depending on $B$ times $N$.
This proves $T=O_B(N)$ and $M=O_B(N/\Delta)$.
Equation \eqref{eq:tail-square} then proves $\E[IW^2]=O_B(N^2/\Delta)$.

Let $U=S-N=W-I$.
The stationary generator identity for $U^3$ and $\E(I-\Delta)=0$ give
\[
 \Delta\E U^2
 =\E[IU^2]+(2N-\Delta)\E U-\E[IU].
\]
Since $(W-I)^2\le W^2+I^2$, $\E U\le\E W$, and $-\E[IU]\le\E I^2$,
\[
 \Delta\E U^2
 \le\E[IW^2]+\E I^3+2N\E W+\E I^2
 \le C_B N^2/\Delta.
\]
Finally, $W^2\le2U^2+2I^2$ and $\Delta^2\le B^2N$ prove \eqref{eq:waiting-second}.
The bounds in \eqref{eq:NDS-moments} follow from $X\le1+W/N$ and $IX^2\le2I+2IW^2/N^2$.

For the NDS profile bound~\eqref{eq:NDS-profile-moment}, we use the product $IP_2$.
There are no simultaneous jumps of $I$ and $P_2$, and $\cL P_2=-2W+Z_2$ whenever $I>0$.
Its stationary generator identity, justified by Lemma~\ref{lem:foundations}, therefore gives
\[
 \E[M_1P_2]=\lambda\E[P_2\ind_{\{I>0\}}]+\E[I(2W-Z_2)].
\]
Since $N-Z_2=M_1+I$ and $I\le N$,
\[
 \E[(N-Z_2)P_2]
 \le(\lambda+N)\E[P_2\ind_{\{I>0\}}]+2\E[IW].
\]
For fixed $\beta$, \eqref{eq:tail-cubic} gives $\E[P_2\ind_{\{I>0\}}]\le2\E W=O_\beta(N)$, and \eqref{eq:waiting-moments} gives $\E[IW]=O_\beta(N)$.
These estimates prove \eqref{eq:NDS-profile-moment}.
\end{proof}

\subsection{Bounds on the NDS occupancy fraction}

\begin{lemma}[Bounds on the occupancy fraction]
\label[lemma]{nds:coefficient-bounds}
Suppose that Assumptions~\ref{ass:base} and~\ref{ass:weighted-occupancy} hold and $\lambda_N=N-\beta$ with fixed $\beta>0$.
Then, for every $x>1$,
\[
 \frac{x-1}{x}\le A(x)\le\min\{x-1,1\}.
\]
The corresponding bounds on $g=(1-A)/A$ are given by~\eqref{nds:idleness-bounds}.
\end{lemma}

\begin{proof}
\emph{Upper bound.}
Write $a_N=Z_{2,N}/N$ and $g=(1-A)/A$.
Take any joint subsequential limit $(\mu,\eta)$ of the ordinary and idleness-weighted population measures in Lemma~\ref{lem:population-bar}.
The two balance equations in that lemma and the interior conclusion of Theorem~\ref{thm:scalar-identification} give a strictly positive continuous density $p$ on $(1,\infty)$.
The pointwise bound $a_N\le W_N/N=X_N-1+I_N/N$ and the ordinary part of~\eqref{path:weighted-coefficient} imply $A(x)\le x-1$ for $\mu$-almost every $x>1$.
Continuity of $A$ and positivity of $p$ extend this inequality to every $x>1$.
Together with $A\le1$, this proves the upper bound.

\emph{Local tail control.}
To prove the lower bound, we compare the limiting occupancy tails with a geometric sequence of ratio $A(x)$.
We first preserve the occupancy mean locally wherever $A<1$, so that the tail comparison can be summed.
Proposition~\ref{prop:common-moments} gives the bound \eqref{eq:NDS-profile-moment}.
For a compact set $K\subset\{x>1:A(x)<1\}$, choose $\delta>0$ and $R<\infty$ such that $A\le1-\delta$ on $K$ and $K\subset(1,R]$.
On $\{X_N\in K,\ |a_N-A(X_N)|\le\delta/2\}$, one has $N-Z_{2,N}\ge\delta N/2$.
On its complement within $\{X_N\in K\}$, we use $H_m/N\le X_N\le R$.
Since $H_m\le P_2/(m-1)$ for $m\ge2$,~\eqref{eq:NDS-profile-moment} gives
\[
 \E[(H_{m,N}/N)\ind_{\{X_N\in K\}}]
 \le\frac{2C_\beta}{\delta(m-1)}
 +R\Pp\{X_N\in K,\ |a_N-A(X_N)|>\delta/2\}.
\]
The final probability tends to zero by~\eqref{path:weighted-coefficient}.
Consequently,
\begin{equation}
 \lim_{m\to\infty}\limsup_{N\to\infty}
 \E[(H_{m,N}/N)\ind_{\{X_N\in K\}}]=0.
 \label{nds:coefficient-local-tails}
\end{equation}

\emph{Lower bound.}
Set $u_m^N=Z_{m+1,N}/N$ for $m\ge1$ and $u_0^N=1$.
After a further diagonal subsequence, the finite measures
\[
 \nu_{m,N}(B)=\E[u_m^N\ind_{\{X_N\in B\}}],\qquad m\ge0,
\]
converge weakly to measures $\nu_m\le\mu$, with $\nu_0=\mu$.
We next show that
\begin{equation}
 \nu_m\le A\nu_{m-1}\quad\text{on }(1,\infty),\qquad m\ge1.
 \label{nds:coefficient-tail-recursion}
\end{equation}
Recall $b_m(Q)=\Pp(Q_\iota\ge m\mid Q)$, the conditional probability that the selected queue contains at least $m$ jobs.
Under Assumption~\ref{ass:base}, for $h\in C_c^2(\R)$, the population generator satisfies $|\cL h(X)|\le C_h(I+1)/N$.
The product rule for $h(X)H_m$, together with $\E[(I+1)H_m]=O_\beta(N)$ from Proposition~\ref{prop:common-moments}, gives
\begin{equation}
 \E[h(X_N)u_m^N]
 =\rho\E[h(X_N)b_m(Q)]+O_h(N^{-1}),
 \qquad \rho=\lambda/N.
 \label{nds:coefficient-level-balance}
\end{equation}
The product of the jump increments of $h(X)$ and $H_m$ has magnitude at most $C_h/N$, so its contribution to this identity is also $O_h(N^{-1})$ after division by $N$.
For the remaining coefficient argument, take $h$ nonnegative and supported in $(1,\infty)$.
The linear idleness identity is
\[
 \cL I=N\{\rho\ind_{\{I=0\}}-a_N\}+\beta-I.
\]
For $f_N=h(X_N)u_{m-1}^N$, its generator and the sum of products of its jump increments with those of $I$ are uniformly bounded.
Stationarity of $I f_N$ and $\E I_N=\beta$ therefore yield
\begin{equation}
 \E\!\left[h(X_N)u_{m-1}^N
       \{\rho\ind_{\{I_N=0\}}-a_N\}\right]=O_h(N^{-1}).
 \label{nds:coefficient-linear-average}
\end{equation}
Idle priority and~\eqref{eq:domination} give $b_m\le\ind_{\{I=0\}}u_{m-1}^N$.
Combining~\eqref{nds:coefficient-level-balance} and~\eqref{nds:coefficient-linear-average} for nonnegative $h$, replacing $a_N$ by $A(X_N)$, and passing to the limit proves~\eqref{nds:coefficient-tail-recursion}.
Thus $\nu_m\le A^m\mu$ in the interior.

The identity $\sum_{m\ge1}u_m^N=W_N/N=X_N-1+I_N/N$ and the local tail bound~\eqref{nds:coefficient-local-tails} imply
\[
 (x-1)\mu(\dd x)=\sum_{m\ge1}\nu_m(\dd x)
 \le\frac{A(x)}{1-A(x)}\mu(\dd x)
 \quad\text{on }\{x>1:A(x)<1\}.
\]
For the equality, one first sums finitely many coordinates and then lets the truncation level tend to infinity; the omitted sum is $H_{m+1,N}/N$ and vanishes locally by~\eqref{nds:coefficient-local-tails}.
Positivity of $p$ and continuity extend the resulting inequality to every point where $A<1$.
The lower bound $A(x)\ge(x-1)/x$ follows there and holds directly where $A=1$.
Transforming both bounds by $g=(1-A)/A$ proves \eqref{nds:idleness-bounds}.
\end{proof}

\subsection{Maximal bounds for the idleness corrector}
\label{app:quadratic-corrector}

Recall the occupancy fraction $a_N(t)=Z_{2,N}(t)/N$, the quadratic idleness function $G_N(I)=I(I+1)/(2N)$, and its generator in \eqref{quad:generator}.
The next estimates make the quadratic corrector negligible on the rescaled clock and control the integrated idleness balance.
For a time acceleration $c_N>0$ and any initial law of the queueing system, write
\[
 m_{1,N}(T)=\sup_{s\le T}\E I_N(c_Ns),
 \qquad
 m_{2,N}(T)=\sup_{s\le T}\E(I_N(c_Ns)+1)^2.
\]

For a locally square-integrable martingale $M$, $\langle M\rangle$ denotes its predictable quadratic variation.

\begin{lemma}[Quadratic corrector maximum and averaging]
\label[lemma]{quad:maximal}
Under \cref{ass:base}, suppose that, for every finite $T$,
\begin{equation}
 \frac{m_{2,N}(T)}{Nc_N}\longrightarrow0,
 \qquad
 \frac{\Delta_Nm_{1,N}(T)}N\longrightarrow0.
 \label{quad:max-smallness}
\end{equation}
Then
\begin{equation}
 \E\sup_{t\le T}
 \frac{I_N(c_Nt)(I_N(c_Nt)+1)}{2Nc_N}
 \longrightarrow0.
 \label{quad:max-conclusion}
\end{equation}
If, in addition,
\begin{equation}
 \frac{m_{2,N}(T)+(\Delta_N+1)m_{1,N}(T)}N\longrightarrow0,
 \label{quad:averaging-smallness}
\end{equation}
then
\begin{equation}
 \E\sup_{t\le T}\left|
 \int_0^t a_N(c_Ns)I_N(c_Ns)\,\dd s
 -\int_0^t\{1-a_N(c_Ns)\}\,\dd s
 \right|\longrightarrow0.
 \label{quad:averaging-conclusion}
\end{equation}
\end{lemma}

\begin{proof}
Set $J_N(t)=G_N(I_N(c_Nt))/c_N$, and let $K_N$ be its compensated martingale.
The idleness jump rates sum to at most $2N$, and each jump of $G_N$ has magnitude at most $(I+1)/N$.
Consequently,
\[
 \E\langle K_N\rangle(T)
 \le\frac{2T}{Nc_N}m_{2,N}(T)\longrightarrow0.
\]
Equation~\eqref{quad:generator} also gives the one-sided drift bound $\cL G_N\le1+\Delta_NI/N$.
For a deterministic grid of mesh $h$, we integrate this bound from the grid point immediately preceding each $t\le T$ to obtain
\[
 J_N(t)\le\max_{j\le\lceil T/h\rceil}J_N(jh\wedge T)
 +h+\frac{\Delta_N}{N}\int_0^T I_N(c_Ns)\,\dd s
 +2\sup_{s\le T}|K_N(s)|.
\]
For fixed $h$, the expected maximum over the grid is at most the sum of its expected values, each bounded by $m_{2,N}(T)/(2Nc_N)$.
Doob's inequality and \eqref{quad:max-smallness} make the other terms vanish except for $h$.
We first let $N\to\infty$ and then $h\downarrow0$ to prove \eqref{quad:max-conclusion}.
Finally, integrating \eqref{quad:generator} expresses the difference of the two integrals in \eqref{quad:averaging-conclusion} as $J_N(0)-J_N(t)+K_N(t)$ plus the integrated last term of that generator.
The bounds on the maxima of $J_N$ and $|K_N|$ over $[0,T]$, together with \eqref{quad:averaging-smallness}, prove the second conclusion.
\end{proof}

Under stationarity, Proposition~\ref{prop:common-moments} gives $m_{1,N}(T)=\Delta_N$ and $m_{2,N}(T)\le\Delta_N^2+\sqrt{2N}\Delta_N+2\Delta_N+1$.
In NDS, where $\Delta_N=\beta$, these bounds verify \eqref{quad:max-smallness} and \eqref{quad:averaging-smallness} both on the original time scale, $c_N=1$, and on the diffusion time scale, $c_N=N$.

In the super-Halfin--Whitt regime, put $\varepsilon_N=\Delta_N/\sqrt N$ and $c_N=N/\Delta_N^2$.
The two quantities in \eqref{quad:max-smallness} are then $O(\varepsilon_N^3)$ and $\varepsilon_N^2$, while the quantity in \eqref{quad:averaging-smallness} is $O(\varepsilon_N)$, so all three vanish.

The corresponding upper bounds, multiplied by the domination constant, hold under initial laws uniformly dominated by stationarity, because the moment bounds transfer at every time.

\subsection{Idleness averaging in NDS}

Fix $\beta>0$ and let $\lambda_N=N-\beta$.
Retain the coefficient $A$ from Assumption~\ref{ass:weighted-occupancy} and put $g=(1-A)/A$.
For compactly supported $h$, all products are extended by zero outside $(1,\infty)$ as in Section~\ref{sec:nds-corrector}.

\begin{lemma}[Interior NDS averaging]
\label[lemma]{path:common-corrector}
Under Assumptions~\ref{ass:base} and~\ref{ass:weighted-occupancy}, every $h\in C_c^1((1,\infty))$ satisfies
\begin{equation}
 \E[h(X_N)\{I_N-g(X_N)\}]\longrightarrow0.
 \label{path:stationary-average}
\end{equation}
For stationary initialization and every finite $T$,
\begin{equation}
 \E\sup_{t\le T}\left|
 \int_0^t h(\widehat X_N(s))
 \{I_N(Ns)-g(\widehat X_N(s))\}\,\dd s
 \right|\longrightarrow0.
 \label{path:common-average}
\end{equation}
The latter conclusion also holds under any sequence of full-state initial laws uniformly dominated by the stationary laws.
\end{lemma}

\begin{proof}
The function $h/A$ has a bounded globally Lipschitz extension that vanishes near and below one.
With this extension, define
\begin{equation}
 \Psi_N(Q)=-\frac{h(X_N)}{A(X_N)}\frac{I_N(I_N+1)}{2N}.
 \label{path:common-corrector-function}
\end{equation}
For each fixed $N$, the corrector is bounded because $I_N\le N$, and the bounded total jump rate justifies Dynkin's formula and the stationary identity $\E\cL\Psi_N=0$.
Applying the generator product rule to \eqref{path:common-corrector-function} and using \eqref{quad:generator} gives
\begin{equation}
 \cL\Psi_N
 =h(X_N)\{I_N-g(X_N)\}
 +\frac{h(X_N)}{A(X_N)}(I_N+1)\{a_N-A(X_N)\}+R_N,
 \qquad |R_N|\le C_{h,A,\beta}\frac{(I_N+1)^2}{N}.
 \label{path:corrector}
\end{equation}
To obtain the remainder bound in \eqref{path:corrector}, each population jump changes the Lipschitz coefficient by at most $C_{h,A}/N$, the total jump rate is at most $2N$, and the quadratic factor is at most $(I_N+1)^2/(2N)$.
The cross terms in the generator product rule satisfy the same bound.
The expectation of the coefficient-replacement term tends to zero by \eqref{path:weighted-coefficient}, and Proposition~\ref{prop:common-moments} gives $\E|R_N|\to0$.
Stationarity of $\Psi_N$ therefore proves \eqref{path:stationary-average}.

Each jump of $\Psi_N$ has magnitude at most $C_{h,A}(I_N+1)/N$.
Its Dynkin martingale consequently satisfies
\[
 \E\left\langle N^{-1}M_N^\Psi(N\cdot)\right\rangle(T)
 \le C_{h,A}T\frac{\E(I_N+1)^2}{N^2}\longrightarrow0.
\]
Doob's inequality controls its maximum in $L^2$.
Lemma~\ref{quad:maximal}, with $c_N=N$, makes the expected maximum of $|\Psi_N|/N$ along the queueing path up to time $NT$ tend to zero.
This controls the corrector values at times $0$ and $Nt$ in Dynkin's formula.
Dynkin's formula at time $Nt$, divided by $N$, and stationarity of the integrated errors prove \eqref{path:common-average}.
Full-state domination bounds every nonnegative path functional by a fixed multiple of its stationary expectation, giving the final assertion.
\end{proof}

\section{Attraction of the mean-preserving profile flow}
\label{app:profile-flow}

To prove Theorem~\ref{thm:profile}, we combine strict contraction at a fixed mean with compactness of profile trajectories.
Fix $F$ as in Example~\ref{ex:fallback-nds}, and let $L$ be a Lipschitz constant for $F$.
For a profile $q\in\mathcal Q$, we use $q_0=1$ and recall the associated probability masses $p_0=1-q_1$ and $p_k=q_k-q_{k+1}$ for $k\ge1$.

\begin{lemma}[Well-posedness, mean preservation, and contraction]\label[lemma]{prof:contraction}
Equation~\eqref{prof:ode} defines a unique global continuous semiflow $S_t$ on $\mathcal Q$ that preserves $m(q)$.
The semiflow preserves coordinatewise order and is nonexpansive in $\ell^1$.
If $q\ne r$ and $m(q)=m(r)$, then
\[
\norm{S_tq-S_tr}_1<\norm{q-r}_1,\qquad t>0.
\]
\end{lemma}

\begin{proof}
\emph{Existence and mean preservation.}
Let $b$ denote the right side of~\eqref{prof:ode}.
For $q,r\in\mathcal Q$, splitting the difference of the arrival terms gives
\[
\norm{b(q)-b(r)}_1\le(2L+3)\norm{q-r}_1.
\]
The term involving $q_1-r_1$ has norm at most $|q_1-r_1|$, because the nonnegative increments $F(r_{k-1})-F(r_k)$ sum to one.
Euler steps preserve $\mathcal Q$ for step size $h\le1/(L+1)$.
Indeed, their probability-mass updates satisfy
\[
p_0^{\mathrm{new}}\ge(1-hL)p_0+hp_1,
\qquad
p_k^{\mathrm{new}}\ge[1-h(L+1)]p_k,\quad k\ge1,
\]
and preserve total probability.
The tail drift is absolutely summable and satisfies
\[
\sum_{k\ge1}b_k(q)
=q_1\sum_{k\ge1}\{F(q_{k-1})-F(q_k)\}-q_1=0.
\]
Thus the Euler construction and the Lipschitz bound give a unique global semiflow on $\mathcal Q$, with exact preservation of the mean.

\emph{Order and contraction.}
For two trajectories $q(t)$ and $r(t)$, put $z=q-r$, $a=q_1$, and
\[
c_k=
\begin{cases}
\dfrac{F(q_k)-F(r_k)}{q_k-r_k},&q_k\ne r_k,\\
0,&q_k=r_k,
\end{cases}
\qquad
h_k=F(r_{k-1})-F(r_k).
\]
Then $0\le c_k\le L$, $h_k\ge0$, and $\sum_{k\ge1}h_k=1$.
With $c_0=z_0=0$, the exact difference equation is
\begin{equation}\label{prof:difference}
\dot z_k=a(c_{k-1}z_{k-1}-c_kz_k)+h_kz_1-z_k+z_{k+1}.
\end{equation}
Its coefficient matrix, in column-vector convention, has nonnegative off-diagonal entries and zero column sums.
For a column $j\ge2$, its downward and upward rates are respectively $1$ and $ac_j$.
Column $1$ has rates $ac_1+h_2$ to index $2$ and $h_k$ to index $k\ge3$, with total rate $ac_1+F(r_1)$.
Every exit rate is therefore at most $L+1$.
The bounded time-dependent Markov evolution $U(t,s)$ associated with~\eqref{prof:difference} is positive and preserves the sum of nonnegative $\ell^1$ vectors, so it is nonexpansive in $\ell^1$.
This proves order preservation and nonexpansion for the nonlinear flow.

For every $j\ge1$ and $t>s$, the path consisting of $j-1$ consecutive downward jumps gives
\[
U_{1j}(t,s)\ge e^{-(L+1)(t-s)}\frac{(t-s)^{j-1}}{(j-1)!}>0.
\]
If $z(s)$ is nonzero and has zero sum, its positive and negative parts both have positive mass.
Their images under $U(t,s)$ consequently overlap at index $1$, which gives strict contraction of their difference in $\ell^1$.
\end{proof}

For $a\in[0,1)$, recall the sequence $v_0^F(a)=1$, $v_{k+1}^F(a)=aF(v_k^F(a))$ and its mean $M_F(a)=\sum_{k\ge1}v_k^F(a)$.

\begin{lemma}[Equilibrium parametrization]\label[lemma]{prof:massmap}
The equilibria are given by~\eqref{prof:equilibrium}, and their parametrization by the mean satisfies
\[
\norm{q^F(y)-q^F(v)}_1=|y-v|.
\]
\end{lemma}

\begin{proof}
The bound $F(z)\le z$ gives $0\le v_k^F(a)\le a^k$.
The profile mean $M_F(a)$ is continuous and strictly increasing on $[0,1)$, with $M_F(0)=0$.
Continuity follows from geometric domination on every compact subinterval of $[0,1)$, and strict increase follows from monotonicity of the recursion and $v_1^F(a)=a$.
For each fixed $k$, continuity of $F$ and $F(1)=1$ imply $v_k^F(a)\to1$ as $a\uparrow1$, so $M_F(a)\to\infty$.
Thus $a_F(y)=M_F^{-1}(y)$ is well defined.

At any finite-mean equilibrium, $q_{k+1}-q_1F(q_k)$ is constant in $k$.
This constant is zero because $q_k\to0$ and $F(q_k)\to0$.
The value $q_1=1$ would force every coordinate to equal one, whereas $q_1=0$ gives the zero profile.
This proves uniqueness and~\eqref{prof:equilibrium}.
Finally, the equilibrium profiles are coordinatewise ordered as their means increase, so their $\ell^1$ distance is the difference of their means.
\end{proof}

\begin{proof}[Proof of Theorem~\ref{thm:profile}]
We first prove attraction for finite-support initial laws using an exponential-moment bound to obtain relative compactness, then extend the result by truncation and nonexpansion.

\emph{Compactness for finite-support profiles.}
We represent the probability masses of a single profile by a birth-and-death chain.
Fix a trajectory whose initial profile has finite support and positive mean $y$, and let $p(t)$ be its associated probability law.
For $p_j(t)>0$, define
\[
\lambda_j(t)=q_1(t)\frac{F(q_j(t))-F(q_{j+1}(t))}{p_j(t)},\qquad j\ge0,
\]
and set $\lambda_j(t)=0$ when $p_j(t)=0$.
These rates belong to $[0,L]$.
The forward equations of the time-inhomogeneous birth-and-death chain $J(t)$ with upward rates $\lambda_j(t)$ and unit downward rate away from zero are precisely the equations for $p(t)$.
In particular, $\E J(t)=y$.
We obtain a lower bound on $p_0(t)$ to keep the arrival coefficient $q_1(t)=1-p_0(t)$ uniformly below one.
Put $K=\lceil2y\rceil$.
At every time $s$, at least half the mass lies on $\{0,\ldots,K\}$.
Using independent potential-birth and potential-service clocks of rates $L$ and one, respectively, yields
\begin{equation}\label{prof:zero-bound}
p_0(t)\ge\delta:=\tfrac12e^{-L}\Pp\{\operatorname{Poisson}(1)\ge K\}>0,
\qquad t\ge1.
\end{equation}
Indeed, no potential birth and at least $K$ potential services during the next unit interval bring the auxiliary chain to zero from every state in this set.

All finite-time exponential moments of $J(t)$ are finite, since $J(t)$ is bounded above by $J(0)$ plus a rate-$L$ Poisson process.
Define $R_t$ through its tails $\Pp(R_t\ge k)=F(q_k(t))$.
The condition $F\le\mathrm{id}$ gives stochastic domination of $R_t$ by $J(t)$.
Writing $a(t)=q_1(t)$ and $M_\theta(t)=\E e^{\theta J(t)}$, the birth-and-death generator therefore gives
\begin{align*}
M_\theta'(t)
&=a(t)(e^\theta-1)\E e^{\theta R_t}
 +(e^{-\theta}-1)\{M_\theta(t)-p_0(t)\}\\
&\le\{a(t)(e^\theta-1)+e^{-\theta}-1\}M_\theta(t)+(1-e^{-\theta}).
\end{align*}
By~\eqref{prof:zero-bound}, $a(t)\le c:=1-\delta<1$ for $t\ge1$.
Choose $0<\theta<-\log c$ and put $\kappa=(e^\theta-1)(e^{-\theta}-c)>0$.
Then $M_\theta'\le-\kappa M_\theta+(1-e^{-\theta})$ on $[1,\infty)$, so $\sup_{t\ge1}M_\theta(t)<\infty$.
Consequently the orbit is relatively compact in $\ell^1$, because
\[
\sup_{t\ge1}\sum_{k>m}q_k(t)
\le\frac{e^{-\theta(m+1)}}{1-e^{-\theta}}\sup_{t\ge1}M_\theta(t)\longrightarrow0
\qquad(m\to\infty).
\]

\emph{Attraction at a fixed mean.}
The distance $\norm{S_tq-q^F(y)}_1$ is nonincreasing, hence has a limit $\ell$.
Choose $t_n\to\infty$ along which $S_{t_n}q$ converges in $\ell^1$ to $v$.
Preservation of the mean gives $m(v)=y$, and continuity of the flow gives
\[
\norm{v-q^F(y)}_1=\ell=\norm{S_hv-q^F(y)}_1,\qquad h>0.
\]
Strict contraction from Lemma~\ref{prof:contraction} forces $v=q^F(y)$ and $\ell=0$.

\emph{Extension and uniformity.}
For a general initial law of finite mean $y$, let $q^{[M]}$ be the tail profile of $J(0)\wedge M$, and write $y_M$ for its mean.
Then $\norm{q-q^{[M]}}_1=y-y_M\to0$.
Nonexpansion and Lemma~\ref{prof:massmap} give
\[
 \norm{S_tq-q^F(y)}_1
 \le2(y-y_M)+\norm{S_tq^{[M]}-q^F(y_M)}_1.
\]
We first let $t\to\infty$, using attraction for the finite-support initial law, and then let $M\to\infty$ to prove~\eqref{prof:attraction}.
The zero-mean case is immediate.
Finally, Lemma~\ref{prof:massmap} gives, for any $q,r\in\mathcal Q$,
\[
\norm{S_tq-q^F(m(q))}_1
\le2\norm{q-r}_1+\norm{S_tr-q^F(m(r))}_1.
\]
A finite covering of a compact subset of $\mathcal Q$ proves compact-uniform convergence.
\end{proof}

\section{Stationary profiles under fixed fallback routing}\label{reg:section}

Throughout this appendix, $\lambda_N=N-\beta$ with fixed $\beta>0$, and $F$ is a fixed fallback function satisfying the conditions in Example~\ref{ex:fallback-nds}.
For idle priority with fallback $F$, we verify Assumption~\ref{ass:weighted-occupancy} and identify the ordinary and idleness-weighted stationary profile limits.
Let $L$ be a Lipschitz constant for $F$.
Write $\rho=\lambda_N/N$, $u_k^N=Z_{k+1,N}/N$ for $k\ge1$, $u_0^N=1$, and $a_N=u_1^N$.
The population and the excess-tail profile satisfy
\begin{equation}
 \|u^N\|_1=\frac{W_N}N=X_N-1+\frac{I_N}N.
 \label{reg:mean-identity}
\end{equation}
Let $X_F$ have density $p_{g_F,\beta}$ in \eqref{eq:generic-density}.

\begin{proposition}[Stationary occupancy under fallback $F$]
\label[proposition]{reg:weighted-law}
Assumption~\ref{ass:weighted-occupancy} holds with $A(x)=a_F(x-1)$, and
\[
 (X_N,u^N)\Rightarrow\bigl(X_F,q^F(X_F-1)\bigr)
 \qquad\text{in }\R\times\ell^1.
\]
For every continuous $h:\R\times\ell^1\to\R$ satisfying $|h(x,u)|\le C(1+|x|)$,
\begin{equation}
 \E[I_Nh(X_N,u^N)]
 \longrightarrow
 \E\left[g_F(X_F)h\bigl(X_F,q^F(X_F-1)\bigr)\right].
 \label{reg:weighted-convergence}
\end{equation}
\end{proposition}

\subsection{Short queues under fixed-threshold priority}

Fix an integer cutoff $\ell\ge1$.
An arrival joins a shortest queue whenever the minimum queue length is below $\ell$ and otherwise uses the fallback tail probabilities $F(Z_m/N)$.
These policies satisfy Assumption~\ref{ass:base}, and $\ell=1$ gives idle priority with fallback $F$.
The cutoff and the fallback function are fixed as $N$ varies; constants in these arguments may depend on $\ell$, $F$, and $\beta$.
The following estimate combines with quadratic moment bounds to control the tails of the upper profiles, both for idle priority and for the threshold-priority profiles in Appendix~\ref{app:threshold-profiles}.

\begin{lemma}[Short queues on bounded population sets]
\label[lemma]{reg:regeneration}
For every fixed $R<\infty$, there are positive constants $d_R,c_R,C_R$ such that, for all sufficiently large $N$,
\begin{equation}
\Pp\{X\le R,\ Z_{\ell+1}>(1-d_R)N\}\le C_Re^{-c_RN}.
\label{reg:regeneration-bound}
\end{equation}
\end{lemma}

\begin{proof}
Fix $R\ge0$.
We use uniform tie breaking within each queue-length class, which preserves the occupancy generator.
While fallback is active, a queue of length $m$ receives arrivals at rate at most
\[
\frac\lambda{M_m}
\left\{F(Z_m/N)-F(Z_{m+1}/N)\right\}
\le\frac{\lambda L}N\le L
\]
whenever $M_m>0$.
We represent fallback arrivals to each queue by thinning an independent rate-$L$ proposal clock and use independent rate-one potential-service clocks.
Priority arrivals can raise a queue's length only from below $\ell$ to at most $\ell$.

Set $K_R=\lceil2(R+1)\rceil$.
On $X(0)\le R+1$, at least $N/2$ queues initially have length at most $K_R$.
For each such queue, consider the event that its fallback proposal clock has no event during $[0,1]$ and its potential-service clock has at least $K_R$ events.
Conditional on the initial state, these events are independent and each has probability
\[
p_R=e^{-L}\Pp\{\operatorname{Poisson}(1)\ge K_R\}>0.
\]
Every successful queue has length at most $\ell$ at time one.
Above $\ell$ it receives no priority arrival, and after its length reaches $\ell$, priority arrivals cannot raise it above that level.
A Chernoff bound therefore gives
\begin{equation}
\sup_{q:\,\sum_jq_j\le(R+1)N}
\Pp_q\{Z_{\ell+1}(1)>(1-p_R/4)N\}
\le e^{-p_RN/16}.
\label{reg:conditional-regeneration}
\end{equation}

To impose the population restriction at time one, we use the standard construction with independent external arrivals $A$ and aggregate potential services $D$.
Their rates are $N-\beta$ and $N$, and
\[
S(1)\ge S(0)+A(1)-D(1).
\]
Consequently,
\begin{equation}
\Pp\{X(1)\le R,\ X(0)>R+1\}
\le\Pp\{A(1)-D(1)<-N\}
\le C_\beta e^{-c_\beta N}.
\label{reg:population-localization}
\end{equation}
The localization estimate transfers the bound from an initially bounded population to the terminal population event used in stationarity.
Combining \eqref{reg:conditional-regeneration} and \eqref{reg:population-localization}, and then using stationarity, proves \eqref{reg:regeneration-bound} with $d_R=p_R/4$.
\end{proof}

\subsection{Concentration on equilibrium profiles}

For every stationary system satisfying Assumption~\ref{ass:base} at fixed $\beta$, the population per server varies negligibly on each fixed physical-time interval $[0,T]$.
The total population has drift $I-\beta$, and its martingale has predictable quadratic variation at most $2Nt$.
Stationarity and $\E I=\beta$ therefore give
\begin{equation}
\E\sup_{0\le t\le T}|X_N(t)-X_N(0)|
\le\sqrt{8T/N}+2\beta T/N\longrightarrow0.
\label{reg:frozen-population}
\end{equation}

For the same fixed cutoff $\ell$, write
\[
v_k^N(t)=\frac{Z_{\ell+k,N}(t)}N,\qquad k\ge1,\qquad v_0^N(t)=1,
\]
so $\|v^N(t)\|_1=H_{\ell,N}(t)/N$.
Suppose that the stationary profiles satisfy
\begin{equation}
\lim_{m\to\infty}\limsup_{N\to\infty}\E[H_{m,N}/N]=0,
\qquad
\E\left|\frac{H_{\ell,N}}N-(X_N-\ell)^+\right|\longrightarrow0,
\label{reg:profile-inputs}
\end{equation}
and that, for every finite $T$,
\begin{equation}
\sup_{0\le t\le T}\frac{|V_{\ell,N}(t)-V_{\ell,N}(0)|}N
\xrightarrow{\Pp}0.
\label{reg:deficit-increments}
\end{equation}
The tail condition yields tightness in $\ell^1$, and the mean approximation fixes the mass of the limiting profile.
The deficit estimate allows us to replace the fallback arrival coefficient $\rho\ind_{\{V_\ell=0\}}$ by $v_1^N$ in the limiting profile equations.

\begin{lemma}[Concentration on the equilibrium profile]
\label[lemma]{reg:ordinary-profile}
The laws of $(X_N,v^N)$ are tight in $\R\times\ell^1$.
Every subsequential limit $\mu$ satisfies
\begin{equation}
x\ge1,
\qquad
v=q^F((x-\ell)^+)
\quad\text{for }\mu\text{-almost every }(x,v).
\label{reg:ordinary-profile-support}
\end{equation}
\end{lemma}

\begin{proof}
The tail bound in \eqref{reg:profile-inputs}, the bound $0\le v_k^N\le1$, and the population second-moment bound in \eqref{eq:NDS-moments} give joint tightness.
The mean approximation in \eqref{reg:profile-inputs} implies $\|v\|_1=(x-\ell)^+$ under every joint limit.
Moreover, $\E(1-X_N)^+\le\E I_N/N=\beta/N$, so every limit has $x\ge1$.

\emph{Averaging the fallback arrival rate.}
We identify the evolution of a stationary profile on a fixed physical-time interval $[0,T]$.
The exact deficit identities are
\begin{equation}
\cL V_\ell=\beta-I-Z_{\ell+1}+\lambda\ind_{\{V_\ell=0\}},
\qquad
\lambda\ind_{\{V_\ell=0\}}-Z_{\ell+1}=I-\beta+\cL V_\ell.
\label{reg:linear-idleness}
\end{equation}
Let $\mathcal M_{V,\ell}$ be the compensated martingale for $V_\ell$; its predictable quadratic variation is at most $2Nt$.
Integrating \eqref{reg:linear-idleness} gives
\[
B_N(t):=\int_0^t\{\rho\ind_{\{V_\ell(s)=0\}}-v_1^N(s)\}\,\dd s
=\frac{V_\ell(t)-V_\ell(0)-\mathcal M_{V,\ell}(t)}N
+\frac1N\int_0^t(I(s)-\beta)\,\dd s.
\]
The first difference vanishes uniformly in probability by \eqref{reg:deficit-increments}.
The normalized martingale has expected squared maximum at most $8T/N$, and the expected absolute value of the integrated drift is at most $2\beta T/N$.
Consequently $\sup_{t\le T}|B_N(t)|\to0$ in probability.

For each fixed $k\ge1$, the occupancy equation is
\[
v_k^N(t)=v_k^N(0)
+\int_0^t\Bigl[
\rho\ind_{\{V_\ell(s)=0\}}
\{F(v_{k-1}^N(s))-F(v_k^N(s))\}
-v_k^N(s)+v_{k+1}^N(s)
\Bigr]\,\dd s+\mathcal M_k^N(t).
\]
Here $\mathcal M_k^N$ is the compensated martingale for $v_k^N$, and $\E\langle\mathcal M_k^N\rangle(T)\le2T/N$.
The convention $v_0^N=1$ makes the equation exact at $k=1$, because the arrival term is active only when every queue has at least $\ell$ jobs.
Let $b_k^N(t)=F(v_{k-1}^N(t))-F(v_k^N(t))$.
Its total variation on $[0,T]$ has expectation at most $4LT$, since a coordinate jumps by $1/N$ and the total event rate is at most $2N$.
Integration by parts yields
\[
\sup_{t\le T}\left|
\int_0^t b_k^N(s)
\{\rho\ind_{\{V_\ell(s)=0\}}-v_1^N(s)\}\,\dd s
\right|
\le\sup_{t\le T}|B_N(t)|
\{1+\operatorname{TV}_{[0,T]}(b_k^N)\}
\xrightarrow{\Pp}0.
\]
We may therefore replace $\rho\ind_{\{V_\ell=0\}}$ by $v_1^N$ in each integrated profile equation, with an error that vanishes uniformly in probability on $[0,T]$.

\emph{Convergence to the profile flow.}
The tail control also holds uniformly on $[0,T]$.
Only arrivals can increase $H_m$, and their expected rate is at most $\lambda\E[Z_m/N]$ by~\eqref{eq:domination}.
Consequently,
\begin{equation}
\E\sup_{0\le t\le T}\frac{H_m(t)}N
\le\E\frac{H_m}N+\rho T\E\frac{Z_m}N
\le\E\frac{H_m}N+\frac{CT}m.
\label{reg:path-tail-bound}
\end{equation}
Here $mZ_m\le S$ and the uniform first population moment justify the last inequality.
The bounded drifts and vanishing martingale maxima give tightness of every finite coordinate projection, with continuous limits.
Equation~\eqref{reg:path-tail-bound} makes $\sup_{t\le T}\sum_{k>m}v_k^N(t)$ uniformly small in probability as $m\to\infty$.
Approximation by these finite projections therefore gives tightness in $D([0,T],\ell^1)$, with continuous limit paths.
Every limiting profile follows \eqref{prof:ode}, while \eqref{reg:frozen-population} makes the limiting population constant.
The initial mean is $(x-\ell)^+$ by \eqref{reg:profile-inputs}, and Theorem~\ref{thm:profile} preserves this mean.

\emph{Stationary concentration.}
Stationarity of the finite systems makes the limiting joint law invariant under $(x,v)\mapsto(x,S_tv)$, as in Bena\"im and Le Boudec~\cite[Theorem~1]{BenaimLeBoudec2011}.
We apply invariance to the bounded continuous function
\[
(x,v)\longmapsto\min\{1,\|v-q^F((x-\ell)^+)\|_1\}.
\]
Theorem~\ref{thm:profile} makes its value along the limiting flow tend to zero, because $\|v\|_1=(x-\ell)^+$.
Dominated convergence therefore proves \eqref{reg:ordinary-profile-support}.
\end{proof}

\subsection{Stationary laws under idle priority}

We now take $\ell=1$, so $v^N=u^N$ and $V_1=I$.

\begin{lemma}[Ordinary occupancy tails]
\label[lemma]{reg:tail-tightness}
For each fixed $R<\infty$,
\begin{equation}
\sup_N\E[(P_2/N)\ind_{\{X\le R\}}]<\infty.
\label{reg:local-p2}
\end{equation}
Moreover,
\begin{equation}
\lim_{m\to\infty}\limsup_{N\to\infty}\E[H_m/N]=0.
\label{reg:ordinary-tail-limit}
\end{equation}
\end{lemma}

\begin{proof}
Proposition~\ref{prop:common-moments} gives \eqref{eq:NDS-profile-moment} under Assumption~\ref{ass:base}.
On $Z_2\le(1-d_R)N$, the factor $N-Z_2$ is at least $d_RN$.
On the complementary event with $X\le R$, one has $P_2\le W^2\le S^2\le R^2N^2$.
Thus \eqref{reg:local-p2} follows from \eqref{eq:NDS-profile-moment} and Lemma~\ref{reg:regeneration} with $\ell=1$.

For $m\ge2$, the pointwise bounds $H_m\le P_2/(m-1)$ and $H_m\le S$ imply
\[
\E[H_m/N]
\le\frac{C_R}{m-1}+\E[X\ind_{\{X>R\}}]
\le\frac{C_R}{m-1}+\frac CR,
\]
where the final inequality uses the uniform second population moment in Proposition~\ref{prop:common-moments}.
We take the limit superior as $N\to\infty$, then let $m\to\infty$ and finally $R\to\infty$ to obtain~\eqref{reg:ordinary-tail-limit}.
\end{proof}

\begin{proof}[Proof of Proposition~\ref{reg:weighted-law}]
Lemma~\ref{reg:tail-tightness} gives the tail bound in \eqref{reg:profile-inputs}.
Since $H_1=W$, \eqref{reg:mean-identity} gives
\[
\E\left|\frac{H_1}N-(X_N-1)^+\right|
\le\frac{\E I}N=\frac\beta N\longrightarrow0,
\]
which verifies its mean approximation.
Lemma~\ref{quad:maximal}, with physical-time clock $c_N=1$, and the common NDS moments give
\[
\E\sup_{0\le t\le T}\frac{I_N(t)(I_N(t)+1)}{2N}
\longrightarrow0.
\]
In particular, $\sup_{t\le T}I_N(t)/N\to0$ in probability, so \eqref{reg:deficit-increments} holds.
Lemma~\ref{reg:ordinary-profile} therefore places every ordinary joint limit on $u=q^F(x-1)$.

Since $a_F(y)=q_1^F(y)$, Lemma~\ref{prof:massmap} gives
\[
 |a_F(y)-a_F(v)|
 \le\|q^F(y)-q^F(v)\|_1=|y-v|,
 \qquad y,v\ge0.
\]
Thus $A(x)=a_F(x-1)$ is positive and one-Lipschitz on $(1,\infty)$, as required by Assumption~\ref{ass:weighted-occupancy}.
Ordinary profile concentration, tested against continuous cutoffs equal to one on a compact $K\subset(1,\infty)$, gives
\[
 \E\bigl[|a_N-A(X_N)|\ind_{\{X_N\in K\}}\bigr]\longrightarrow0.
\]
Together with~\eqref{reg:ordinary-tail-limit}, Proposition~\ref{prop:weighted-occupancy-criterion} proves Assumption~\ref{ass:weighted-occupancy}.
Theorem~\ref{thm:regular-nds} now gives $(X_N,I_N)\Rightarrow(X_F,I)$, where $\E[I\mid X_F]=g_F(X_F)$, and uniform integrability of $(I_N)$.

The ordinary profile laws are tight in $\R\times\ell^1$, and every subsequential limit of $(X_N,u^N,I_N)$ has $u=q^F(x-1)$ and the joint population--idleness law just identified.
Thus
\[
 (X_N,u^N,I_N)\Rightarrow\bigl(X_F,q^F(X_F-1),I\bigr),
\]
which also proves the ordinary joint convergence.
Uniform integrability of $I_N$ yields~\eqref{reg:weighted-convergence} first for bounded continuous $h$.
The bound $\sup_N\E[I_NX_N^2]<\infty$ in Proposition~\ref{prop:common-moments} then allows a population cutoff to be removed and gives the stated convergence for continuous $h$ of at most linear growth in $x$.
\end{proof}

\section{Occupancy estimates uniform over priority thresholds}\label{pri:section}

Throughout this appendix, $\lambda_N=N-\beta$ with fixed $\beta>0$ and the routing rule satisfies Assumption~\ref{ass:base}.
We first establish estimates at any threshold where arrivals have exact priority for shorter queues.
Their constants are independent of the threshold, so the same argument applies to one-job priority, fixed higher thresholds, and growing thresholds.

For an integer $j\ge2$ at which exact priority holds, write
\[
 V=V_j,\qquad H=H_j,\qquad P=P_{2,j},\qquad
 Z=Z_{j+1},\qquad \Omega=\min(V,H).
\]
Thus $V>0$ implies that every arrival joins a queue shorter than $j$.
Here $\Omega$ measures the coexistence of missing jobs below $j$ and excess jobs above it.
Since $S=jN-V+H$, we have $H=(S-jN)^++\Omega$, so bounding $\Omega$ controls the error in recovering $H$ from $S$.

\begin{lemma}[Uniform estimates at exact-priority thresholds]
\label[lemma]{pri:marked}
For all sufficiently large $N$, the following bounds hold simultaneously at every exact-priority threshold $j\ge2$, with a constant $C_\beta$ independent of $j$, $N$, and the routing policy:
\[
 \E[P\ind_{\{V>0\}}]\le2\E H\le C_\beta N,
 \qquad
 \E\Omega\le C_\beta N^{2/3},
 \qquad
 \E[(I+1)\Omega]\le C_\beta N^{5/6}.
\]
\end{lemma}

\begin{proof}
While $V>0$, arrivals do not change $H$ or $P$.
Applying the generator to $P_{3,j}-H_j$, as in the proof of Theorem~\ref{lem:tail-transfer}, and using exact priority below $j$ gives
\[
 B+(1-\lambda/N)\E[P\ind_{\{V=0\}}]
 \le\E H+(\lambda/N)\E[H\ind_{\{V=0\}}]\le2\E H,
 \qquad B=\E[P\ind_{\{V>0\}}].
\]
Since $H\le W$, Proposition~\ref{prop:common-moments} gives $B\le C_\beta N$ uniformly in $j$.

\emph{Ordinary overlap.}
We control the deficit when many queues extend above the priority threshold.
Choose an integer $K<N$ with $d_K=K-\beta\ge1$, and put
\[
 \mathcal A_K=\{Z\ge K\},\qquad
 \rho_K=\frac{N-K}{N-\beta}.
\]
The upward rate of $V$ is $N-I-Z$, and its downward rate while positive is $\lambda$.
Moreover, $Z$ can increase only at $V=0$, and such an increase leaves $V$ equal to zero.
Applying stationarity to $\ind_{\{V\ge v,Z\ge K\}}$ for $v\ge1$ therefore gives
\[
 \lambda\Pp(V=v,\mathcal A_K)
 \le(N-K)\Pp(V=v-1,\mathcal A_K).
\]
Entries into $\mathcal A_K$ occur at $V=0$ and contribute zero, whereas exits contribute nonpositively.
Summing over higher levels and iterating yields
\[
 \Pp(V\ge v,\mathcal A_K)
 \le\rho_K^v\Pp(\mathcal A_K),\qquad v\ge0.
\]
Summation of the tails then gives
\begin{equation}
 \E[V\ind_{\mathcal A_K}]
 \le\frac{\rho_K}{1-\rho_K}\Pp(\mathcal A_K),
 \qquad
 \E[V^2\ind_{\mathcal A_K}]
 \le\frac{\rho_K(1+\rho_K)}{(1-\rho_K)^2}\Pp(\mathcal A_K).
 \label{pri:marked-v-moments}
\end{equation}

On $Z<K$, we use $H^2\le KP$ and $\Omega=0$ on $V=0$.
On $\mathcal A_K$, we use $\Omega\le V$.
Then
\[
 \frac{\E\Omega}{N}
 \le\frac{\sqrt{KB}}N+\frac1{K-\beta}
 \le C_\beta\sqrt{K/N}+\frac1{K-\beta}.
\]
Taking $K=\lceil N^{1/3}\rceil$ proves the ordinary overlap estimate.

\emph{Idleness-weighted overlap.}
For the weighted bound, we again take an arbitrary integer $K<N$ with $K-\beta\ge1$.
Let $D=\{V>\lambda/2\}$, $m=\lfloor\lambda/2\rfloor+1$, and $e_K=e^{-d_K/2}$.
Because $1-\rho_K=d_K/\lambda$, the same geometric tail gives
\begin{equation}
 \E[V\ind_{D\cap\mathcal A_K}]
 \le\left(m-1+\frac1{1-\rho_K}\right)
       \rho_K^m\Pp(\mathcal A_K)
 \le2Ne_K\Pp(\mathcal A_K).
 \label{pri:marked-large-deficit}
\end{equation}
The geometric tail makes this bound uniform in $j$.
For $j\ge2$, $M_1\le V$.
Applying the stationary generator identity to $I(I+1)\ind_{\mathcal A_K}$ and $I(I+1)(I+2)\ind_{\mathcal A_K}$ gives
\[
 \lambda\E[I\ind_{\mathcal A_K}]\le\E[V(I+1)\ind_{\mathcal A_K}],
 \qquad
 \lambda\E[I(I+1)\ind_{\mathcal A_K}]
 \le\E[VI(I+1)\ind_{\mathcal A_K}]+2\E[V(I+1)\ind_{\mathcal A_K}].
\]
There are no positive entry terms because $I=0$ whenever $Z$ increases.
Splitting the first inequality over $D$ and its complement, and using $I\le N$, gives
\[
 \frac\lambda2\E[I\ind_{\mathcal A_K}]
 \le\E[V\ind_{\mathcal A_K}]+N\E[V\ind_{D\cap\mathcal A_K}].
\]
Equations~\eqref{pri:marked-v-moments} and~\eqref{pri:marked-large-deficit}, with $\lambda\ge N/2$, then give
\begin{equation}
 \E[I\ind_{\mathcal A_K}]
 \le(2/d_K+8Ne_K)\Pp(\mathcal A_K).
 \label{pri:marked-i}
\end{equation}
For $J_K=\E[I(I+1)\ind_{\mathcal A_K}]$, the second inequality gives
\[
 J_K\le2\E[I\ind_{\mathcal A_K}]
       +\frac4\lambda\E[V\ind_{\mathcal A_K}]
       +\frac{8N^3+8N^2}{\lambda}e_K\Pp(\mathcal A_K).
\]
Combining \eqref{pri:marked-i} with the deficit bound \eqref{pri:marked-v-moments} gives
\begin{equation}
 \E[I^2\ind_{\mathcal A_K}]
 \le(8/d_K+48N^2e_K)\Pp(\mathcal A_K).
 \label{pri:marked-i2}
\end{equation}

On $Z<K$, the weighted contribution, divided by $N$, is bounded by
\[
 \frac1N\sqrt{\E(I+1)^2}\sqrt{KB}
 \le C_\beta N^{-1/4}\sqrt K.
\]
When $K$ grows as any fixed positive power of $N$, \eqref{pri:marked-v-moments} and \eqref{pri:marked-i2} give
\[
 \frac1N\E[(I+1)V\ind_{\mathcal A_K}]
 \le\frac1N\left\{\E[V\ind_{\mathcal A_K}]
 +\sqrt{\E[I^2\ind_{\mathcal A_K}]\E[V^2\ind_{\mathcal A_K}]}\right\}
 \le C_\beta/K.
\]
Taking $K=\lceil N^{1/6}\rceil$ proves the weighted overlap estimate.
\end{proof}

At threshold $j=2$, identity~\eqref{pri:overlap-identity} converts the overlap bound into the weighted one-job approximation used to verify Assumption~\ref{ass:weighted-occupancy} in Section~\ref{sec:nds-verification}.

\begin{proposition}[Weighted one-job approximation]
\label[proposition]{pri:weighted-occupancy}
Uniformly over the one-job priority class,
\[
 \E\left[(I+1)\left|\frac{M_1}{N}-(2-X)^+\right|\right]
 =O_\beta(N^{-1/6}).
\]
\end{proposition}

\begin{proof}
We multiply the absolute-value bound from \eqref{pri:overlap-identity} by $(I+1)/N$.
Lemma~\ref{pri:marked} bounds the overlap term by $C_\beta N^{-1/6}$, and Proposition~\ref{prop:common-moments} bounds the remaining term by $2\E[I(I+1)]/N=O_\beta(N^{-1/2})$.
\end{proof}

\section{Process limits and interchange of limits}
\label{app:initialization}

For NDS and the super-Halfin--Whitt regime, we prove stationary process convergence and martingale approximations whose errors remain small under the initial laws in \eqref{eq:common-dominated-initialization}.
Proposition~\ref{unif:stationary-transfer} gives the general transfer argument and the resulting interchange of limits.
For an observation map $\phi_N$ and time scale $c_N>0$, put $Y_N(t)=\phi_N(Q^N(c_Nt))$ and write
\[
 \mathcal L_N f(q)=c_N\cL(f\circ\phi_N)(q).
\]
This operator applies the full queueing generator to a test function of the observation and is evaluated at the full state $q$.
With limiting generator $\mathcal A$ and test class $\mathcal D$, the common martingale approximation is
\begin{equation}
 f(Y_N(t))-f(Y_N(0))
 -\int_0^t\mathcal Af(Y_N(s))\,\dd s
 =M_N^f(t)+R_N^f(t),
 \qquad
 \E_{\Pi_N}\sup_{t\le T}|R_N^f(t)|\longrightarrow0,
 \label{unif:martingale-approximation}
\end{equation}
where $M_N^f$ is an integrable martingale relative to the full finite-system filtration under the stationary initial law $\Pi_N$, and the bound on the expected supremum is required for every $f\in\mathcal D$ and every finite $T$.
For the queueing verifications, the tests are smooth and compactly supported in the limiting state space, with smooth extensions to the ambient Euclidean space.

\subsection{Well-posedness and ergodicity of the NDS diffusion}

Let $A:(1,\infty)\to(0,1]$ be locally Lipschitz and satisfy \eqref{nds:occupancy-bounds}, put $g=(1-A)/A$, and let $P_t^g$ denote the transition semigroup of
\[
 \dd X(t)=\sqrt2\,\dd B(t)+\{g(X(t))-\beta\}\,\dd t,
 \qquad X(t)>1.
\]
Write $\pi_{g,\beta}(\dd x)=p_{g,\beta}(x)\,\dd x$.
The diffusion has generator $\mathcal Af=f''+(g-\beta)f'$ on $\mathcal D=C_c^\infty((1,\infty))$.

\begin{lemma}[Well-posedness and ergodicity]
\label[lemma]{int:ergodicity}
For every initial probability law $\nu$ on $(1,\infty)$, this SDE has a unique nonexplosive solution, its $(\mathcal A,\mathcal D)$ martingale problem is well posed, and
\[
 \|\nu P_t^g-\pi_{g,\beta}\|_{\mathrm{TV}}\longrightarrow0
 \qquad\text{as }t\to\infty.
\]
If $\nu\le C\pi_{g,\beta}$ for a finite constant $C$, the convergence also holds in Wasserstein distance of every finite order $p\ge1$.
\end{lemma}

\begin{proof}
The positivity and local Lipschitz continuity of $A$ make $g=(1-A)/A$ locally Lipschitz on $(1,\infty)$.
Local strong existence and pathwise uniqueness follow.

The scale density is proportional to
\[
 s'(x)=\exp\left\{\beta x-\int_2^xg(u)\,\dd u\right\}.
\]
The bounds in \eqref{nds:idleness-bounds} show that $s'(1+y)$ is of order $1/y$ as $y\downarrow0$.
At infinity, $s'(x)$ is bounded below by a positive constant times $e^{\beta x}/(x-1)$.
Thus the scale integrals diverge at $x=1$ and at infinity.
The speed density, proportional to $1/s'$, is integrable, is of order $x-1$ near one, and has an exponential upper tail.
The scale function maps $(1,\infty)$ onto $\R$, and the speed measure has finite mass.
The boundary classification and recurrence theorems therefore show that the diffusion does not reach $x=1$ in finite time and is nonexplosive and positive recurrent, with invariant probability $\pi_{g,\beta}$~\cite[Theorems~23.12 and~23.15]{Kallenberg2002}.
Local Lipschitz continuity also gives uniqueness for the stopped martingale problem on every compact subinterval, by its equivalence with the stopped SDE~\cite[Chapters~4 and~5]{EthierKurtz1986}.
The boundary and nonexplosion conclusions allow these stopped problems to exhaust the full time interval, proving well-posedness.
Regular positive recurrent diffusions converge in total variation to their invariant probability from every initial law~\cite[Lemma~23.17]{Kallenberg2002}.
If $\nu\le C\pi_{g,\beta}$, invariance gives $\nu P_t^g\le C\pi_{g,\beta}$ for every $t$.
The invariant density has moments of every order, so the $p$th moments of this dominated family are uniformly integrable for every finite $p\ge1$, and total-variation convergence implies Wasserstein-$p$ convergence.
\end{proof}

\subsection{Stationary process convergence in NDS}
\label{app:stationary-process}

We work under stationary initialization and Assumptions~\ref{ass:base} and~\ref{ass:weighted-occupancy}, with $\Delta_N=\beta>0$.
We establish convergence of $\widehat X_N$ to the diffusion~\eqref{path:sde} and verify the martingale approximation~\eqref{unif:martingale-approximation}.

\paragraph{Localization and path convergence.}

The exact population decomposition is
\begin{equation}
\widehat X_N(t)=\widehat X_N(0)+\widehat M_N(t)
+\int_0^t(I_N(Ns)-\beta)\,\dd s,
\qquad
\langle\widehat M_N\rangle(t)
=\int_0^t\left(2-\frac{\beta+I_N(Ns)}N\right)\,\dd s.
\label{path:population}
\end{equation}
The jumps of $\widehat M_N$ have size $1/N$.
Stationarity and $\E I=\beta$ give convergence of its predictable quadratic variation to $2t$, uniformly in probability on bounded intervals.

Assumption~\ref{ass:weighted-occupancy} permits the use of Lemma~\ref{path:common-corrector} with the same coefficient $A$ and conditional mean idleness $g=(1-A)/A$ that identify the stationary population law.
We first identify the paths on compact population intervals and then use boundary inaccessibility to remove the localization.

For $m\ge3$, put $D_m=(1+1/m,m)$ and $E_m=(1+1/(2m),2m)$, so $\overline D_m\subset E_m$.
We localize the drift on the larger interval $E_m$ and test departures from the inner interval $D_m$.
Let $\tau_N^{(m)}$ be the first exit of $\widehat X_N$ from $E_m$, and choose $h\in C_c^1((1,\infty))$ equal to one on $\overline E_m$.
Stopping the common averaging estimate gives
\begin{equation}
 \label{path:stopped-average}
 R_N^{(m)}(t)
 =\int_0^{t\wedge\tau_N^{(m)}}
       \{I_N(Ns)-g(\widehat X_N(s))\}\,\dd s,
 \qquad
 \E\sup_{t\le T}|R_N^{(m)}(t)|\longrightarrow0.
\end{equation}

Write $b=g-\beta$, and extend its restriction to $\overline E_m$ to a bounded globally Lipschitz function $b^{(m)}$ by composing $b$ with projection onto $\overline E_m$.
Let $\widetilde X_N^{(m)}$ be the unique solution of
\[
 \widetilde X_N^{(m)}(t)
 =\widehat X_N(0)+\widehat M_N(t)
 +\int_0^t b^{(m)}(\widetilde X_N^{(m)}(s))\,\dd s
 +R_N^{(m)}(t).
\]
Equation~\eqref{path:population} and Gronwall's inequality give $\widetilde X_N^{(m)}=\widehat X_N$ through $\tau_N^{(m)}$, including the exit jump.
Applying the martingale functional central limit theorem jointly with the time-zero variable in the full queueing filtration~\cite[Theorem~2.1]{Whitt2007}, and using Theorem~\ref{thm:regular-nds}, gives $(\widehat X_N(0),\widehat M_N)\Rightarrow(X(0),\sqrt2 B)$, where $X(0)$ has law $\pi_{g,\beta}$ and $B$ is independent of $X(0)$.
Continuous dependence of the integral equation on its forcing path, together with~\eqref{path:stopped-average}, therefore gives $\widetilde X_N^{(m)}\Rightarrow\widetilde X^{(m)}$, where $\dd\widetilde X^{(m)}=\sqrt2\,\dd B+b^{(m)}(\widetilde X^{(m)})\,\dd t$ and $\widetilde X^{(m)}(0)=X(0)$.

We couple $\widetilde X^{(m)}$ with the diffusion $X$ in~\eqref{path:sde} using this same initial state and Brownian motion.
Pathwise uniqueness makes them agree until leaving $E_m$.
The path event
\[
 \mathcal H_m
 =\left\{z:\inf_{t\le T}z(t)\le1+1/m
       \text{ or }\sup_{t\le T}z(t)\ge m\right\}
\]
is closed in the Skorokhod topology and occurs for an original path if and only if it occurs for its auxiliary path, because $\overline D_m\subset E_m$.
Portmanteau's theorem consequently gives
\[
 \limsup_{N\to\infty}\Pp\{\widehat X_N\in\mathcal H_m\}
 \le\Pp\{X\in\mathcal H_m\}.
\]
The right-hand side tends to zero as $m\to\infty$, because Lemma~\ref{int:ergodicity} gives continuous paths that remain in $(1,\infty)$ throughout $[0,T]$.
The probability that either original path differs from its auxiliary path on $[0,T]$ is bounded by its probability of $\mathcal H_m$, so convergence together proves the asserted global path convergence.
Finally, every prelimit process is stationary, so its limiting diffusion is stationary as well.

\paragraph{Martingale approximation.}

Here $\phi_N(q)=S(q)/N$ and $c_N=N$.
For $f\in C_c^3((1,\infty))$, the limiting generator is
\[
 \mathcal Af(x)=f''(x)+\{g(x)-\beta\}f'(x).
\]
Both $f$ and $\mathcal Af$ extend by zero to bounded continuous functions on $\R$.
Taylor expansion of the exact population generator leaves the averaging term
\[
 \int_0^t f'(\widehat X_N(s))\{I_N(Ns)-g(\widehat X_N(s))\}\,\dd s.
\]
Lemma~\ref{path:common-corrector}, with $h=f'$, makes its maximum converge to zero in $L^1$.
The error in the diffusion coefficient and the third-order Taylor remainder have vanishing $L^1$ integrals by Proposition~\ref{prop:common-moments}.
The Dynkin martingale of $f(S_N/N)$ is a martingale for the full queueing filtration, so \eqref{unif:martingale-approximation} follows.

\subsection{Stationary process convergence in the super-Halfin--Whitt regime}
\label{app:shw-process}

Under Assumption~\ref{ass:base} and \eqref{shw:scaling}, we establish convergence of the stationary processes $\widehat Y_N$ to the diffusion~\eqref{shw:limit-sde} and verify the martingale approximation~\eqref{unif:martingale-approximation}.
Write $\Delta=\Delta_N$ and $c_N=N/\Delta^2$.
In the calculations below, $I(t)$, $Z_k(t)$, and $H_2(t)$ denote queue statistics at physical time $c_Nt$.

\paragraph{Accumulated idleness and path tightness.}

We define the accumulated idleness drift by
\[
 A_N(t)=\int_0^t\frac{I(s)}{\Delta_N}\,\dd s.
\]
The exact population equation is
\begin{equation}
 \widehat Y_N(t)=\widehat Y_N(0)+M_N(t)+A_N(t)-t,
 \label{shw:path-population}
\end{equation}
where $M_N$ has jumps of size $\Delta/N$ and predictable quadratic variation
\[
 \langle M_N\rangle(t)
 =\int_0^t\left\{2-\frac\Delta N-\frac{I(s)}N\right\}\,\dd s.
\]
Its predictable quadratic variation converges uniformly in probability on bounded intervals to $2t$.
The martingale functional central limit theorem gives $M_N\Rightarrow\sqrt2 B$, jointly with the initial state; the limiting Brownian motion is independent of that state~\cite[Theorem~2.1]{Whitt2007}.

We use the quadratic idleness estimate to identify the limiting accumulated idleness.
By the moment bounds in Appendix~\ref{app:quadratic-corrector}, Lemma~\ref{quad:maximal} applies to
\[
 J_N(t)=\frac{\Delta_N^2}{2N^2}I(t)(I(t)+1)
\]
and yields
\begin{equation}
 \E\sup_{t\le T}J_N(t)\longrightarrow0,
 \qquad
 \E\sup_{t\le T}\left|
 \int_0^t\frac{I(s)Z_2(s)}N\,\dd s-t\right|\longrightarrow0.
 \label{shw:path-corrector-maximum}
\end{equation}
For the second assertion we also used $\E Z_2/N=O(\Delta_N^{-1})$.
Equation~\eqref{shw:collapse-substitution}, Proposition~\ref{prop:common-moments}, and \eqref{shw:collapse-bounds} show that the integrated difference between $I Z_2/N$ and $I\widehat Y_N/\Delta_N$ vanishes in $L^1$.
Consequently
\begin{equation}
 \E\sup_{t\le T}\left|
 \int_0^t\widehat Y_N(s)\,\dd A_N(s)-t\right|\longrightarrow0.
 \label{shw:path-averaging}
\end{equation}

We also control idleness accumulated near zero, where the population factor in \eqref{shw:path-averaging} is small.
Recall the stationary idleness-weighted probability measure $\theta_N$ from Section~\ref{sec:scalar-stationary}.
For every fixed $\eta>0$, \eqref{shw:weighted-limit} and stationarity give
\begin{equation}
 \limsup_{N\to\infty}\E
 \int_0^T\ind_{\{\widehat Y_N(s)\le\eta\}}\,\dd A_N(s)
 \le T(1-e^{-\eta}).
 \label{shw:path-boundary-drift}
\end{equation}
Also, $\widehat Y_N\ge-\Delta I/N$ implies
\[
 \E\int_0^T(\widehat Y_N(s))^-\,\dd A_N(s)
 \le\frac{T}{N}\E I_N^2\longrightarrow0.
\]
Put $C_N(t)=\int_0^t\widehat Y_N(s)\,\dd A_N(s)$ and $e_N=\sup_{t\le T}|C_N(t)-t|$.
For $0\le s\le t\le T$, positivity of $\dd A_N$ gives
\[
 A_N(t)-A_N(s)
 \le\int_0^T\ind_{\{\widehat Y_N(u)\le\eta\}}\,\dd A_N(u)
 +\frac1\eta\left\{t-s+2e_N
 +\int_0^T(\widehat Y_N(u))^-\,\dd A_N(u)\right\}.
\]
Equations~\eqref{shw:path-averaging} and~\eqref{shw:path-boundary-drift} prove the uniform modulus-of-continuity criterion by first letting $N\to\infty$, then $t-s\downarrow0$, and finally $\eta\downarrow0$.
Since $\E A_N(T)=T$, the processes $A_N$ are $C$-tight.
Equation~\eqref{shw:path-population} then makes $\widehat Y_N$ $C$-tight as well.

\paragraph{Identification of the limiting diffusion.}

Consider a joint subsequential limit $(Y,A,M)$.
The martingale identities pass to the joint limiting filtration, so $M=\sqrt2 B$ remains Brownian in that filtration.
Since $A_N$ is nondecreasing and all limiting paths are continuous, convergence of the Stieltjes integrals in \eqref{shw:path-averaging} gives
\[
 Y(t)=Y(0)+\sqrt2 B(t)+A(t)-t,
 \qquad
 \int_0^tY(s)\,\dd A(s)=t.
\]
Moreover, $\widehat Y_N\ge-\sqrt{2J_N}$ and \eqref{shw:path-corrector-maximum} imply $Y\ge0$.
Stationary convergence gives $\Law(Y(0))=\pi_{\mathrm{shw}}$, so $Y(0)>0$ almost surely.
On each compact subinterval of $(0,\infty)$, the measure identity $Y\,\dd A=\dd t$ implies $\dd A=\dd t/Y$, and the stopped process solves \eqref{shw:limit-sde}.
This is the diffusion of Lemma~\ref{int:ergodicity} with $g(x)=1/(x-1)$, $\beta=1$, and the translation $X=1+Y$.
Local uniqueness and the fact that this diffusion stays in $(0,\infty)$ for all finite times identify the full process as the compact intervals expand.
This proves stationary process convergence.

\paragraph{Martingale approximation.}

Put $\phi_N(q)=\Delta_N(S(q)-N)/N$ and $\delta_N=\Delta_N/N$.
The joint $C$-tightness just proved and \eqref{shw:path-averaging} give
\[
 \sup_{t\le T}|C_N(t)-t|\longrightarrow0
 \quad\text{in probability}.
\]
Also,
\[
 \operatorname{TV}_{[0,T]}(C_N)
 \le \sup_{t\le T}|\widehat Y_N(t)|A_N(T)
\]
is tight.
For $f\in C_c^3((0,\infty))$, extend $f$ by zero and put $k(y)=f'(y)/y$, again extended by zero near and below zero.
This is a bounded continuous function, and $yk(y)=f'(y)$ for all $y\in\R$.
Consequently
\begin{equation}
 \sup_{t\le T}\left|
 \int_0^t f'(\widehat Y_N(s))\,\dd A_N(s)
 -\int_0^t\frac{f'(\widehat Y_N(s))}{\widehat Y_N(s)}\,\dd s
 \right|
 =\sup_{t\le T}\left|\int_0^t k(\widehat Y_N(s))\,\dd\{C_N(s)-s\}\right|
 \longrightarrow0
 \label{unif:shw-test-average}
\end{equation}
in probability.
Here the quotient is zero off the support of $f'$.
We approximate $k(\widehat Y_N)$ by its values on a deterministic mesh.
For a fixed mesh, the step-function integrals vanish uniformly because $C_N$ converges uniformly to the identity.
The approximation error is bounded by the modulus of continuity of $k(\widehat Y_N)$ on that mesh times $\operatorname{TV}(C_N)+T$.
$C$-tightness and the preceding variation bound make this error vanish as the mesh decreases.

Applying the full queueing generator to the population test gives
\[
 \mathcal L_N f(q)
 =f''(y)+\left(\frac{I(q)}\Delta-1\right)f'(y)+e_N^f(q),
 \qquad y=\phi_N(q),
\]
where
\[
 \E_{\Pi_N}\int_0^T|e_N^f(Q^N(c_Ns))|\,\dd s
 \le C_fT\frac{\Delta+\E_{\Pi_N}I}{N}
 \longrightarrow0.
\]
The third-order Taylor error is absorbed in this bound.
Thus \eqref{unif:shw-test-average} makes the remainder against
\[
 \mathcal Af(y)=f''(y)+(1/y-1)f'(y)
\]
converge uniformly to zero in probability.
To obtain the bound on the expected supremum in \eqref{unif:martingale-approximation}, we establish uniform integrability of the remainder maxima.
Because $f$ is compactly supported in $(0,\infty)$, both $f$ and $\mathcal Af$ are bounded.
The Dynkin decomposition therefore reduces this to an $L^2$ bound for the test-function martingale maxima.
Let $M_N^f$ be the corresponding full-filtration Dynkin martingale.
Each population jump has size $\delta_N$, and the original total jump rate affecting population is at most $2N$.
On the accelerated clock $N/\Delta^2$,
\[
 \langle M_N^f\rangle(T)
 \le2T\|f'\|_\infty^2.
\]
Doob's inequality bounds its maximum uniformly in $L^2$.
Uniform integrability then strengthens the remainder convergence to \eqref{unif:martingale-approximation}.

\subsection{Transfer to initial laws dominated by stationarity}

The following proposition formulates the transfer argument for observations of a general Markov process.
Let $Z_N$ be a Markov process with invariant law $\Pi_N$, let $c_N>0$, and define the scaled observation $Y_N(t)=\phi_N(Z_N(c_Nt))\in\R^d$.
Write $\mu_N=\Pi_N\circ\phi_N^{-1}$.
All initial laws below use the same transition kernels after time zero.
Let $(\mathcal A,\mathcal D)$ be a well-posed martingale problem on a state space $E\subseteq\R^d$ whose solutions have continuous paths.
For every $f\in\mathcal D$, assume that $f$ and $\mathcal Af$ have bounded continuous extensions to $\R^d$, used in the prelimit expressions.

We use \eqref{unif:martingale-approximation} with integrable martingales for the full filtration of $Z_N$ on the scaled clock under $\Pi_N$ and require its remainder bound for every finite $T$.
Consider any sequence of initial laws satisfying
\[
 \zeta_N\le C\Pi_N,\qquad
 \zeta_N\circ\phi_N^{-1}\Rightarrow\nu,
\]
with a constant $C$ independent of $N$.
For Wasserstein convergence of order $p\ge1$, we also use the moment condition
\begin{equation}
 \lim_{R\to\infty}\sup_N
 \int_{\{\|y\|>R\}}\|y\|^p\,\mu_N(\dd y)=0.
 \label{unif:moment-ui}
\end{equation}

\begin{proposition}[Transfer of process limits and interchange]
\label[proposition]{unif:stationary-transfer}
Suppose $\mu_N\Rightarrow\mu$, \eqref{unif:martingale-approximation} holds for every $f\in\mathcal D$, and, under $\Pi_N$, the processes $Y_N$ are tight on every finite time interval, with all subsequential limits continuous and $E$-valued.
Then $\mu$ is invariant, and under $\zeta_N$ the processes $Y_N$ converge to the solution of $(\mathcal A,\mathcal D)$ started from $\nu$.

If additionally each finite system converges in total variation to $\Pi_N$, the limiting process converges in total variation to $\mu$ from every initial law dominated by a finite multiple of $\mu$, and \eqref{unif:moment-ui} holds for some $p\ge1$, then fixed-time convergence and both iterated limits hold in Wasserstein distance of order $p$:
\begin{equation}
 \lim_{N\to\infty}\lim_{t\to\infty}\Law_{\zeta_N}(Y_N(t))
 =\mu
 =\lim_{t\to\infty}\lim_{N\to\infty}\Law_{\zeta_N}(Y_N(t)).
 \label{unif:two-orders}
\end{equation}
In particular, for $p=1$, the corresponding coordinate means have the same iterated limits.
\end{proposition}

\begin{proof}
Let $P_s^N$ be the transition semigroup of the full prelimit process $Z_N$ in its original time, and let $P_t$ be the transition semigroup of the limiting process.
On every finite time interval, the full path law under $\zeta_N$ is bounded by $C$ times the full path law under $\Pi_N$.
Indeed, its Radon--Nikodym derivative is the initial-state density $\dd\zeta_N/\dd\Pi_N$ evaluated at $Z_N(0)$.
Consequently tightness transfers, and every subsequential observation-path limit under $\zeta_N$ is dominated by $C$ times a corresponding stationary path limit.
It therefore has continuous paths in $E$.
The same inequality transfers the remainder estimate in \eqref{unif:martingale-approximation}.
Because the change of measure is measurable at time zero, $M_N^f$ remains an integrable martingale under $\zeta_N$.

Fix $0\le t_1<\cdots<t_k\le s<t$ and a bounded continuous $H:\R^{dk}\to\R$.
We multiply the increment of \eqref{unif:martingale-approximation} from $s$ to $t$ by $H(Y_N(t_1),\ldots,Y_N(t_k))$ and take expectations under $\zeta_N$.
The martingale term has expectation zero, and the absolute remainder is at most
\[
 2C\|H\|_\infty\,\E_{\Pi_N}\sup_{u\le t}|R_N^f(u)|\longrightarrow0.
\]
Along a convergent subsequence, continuity of the limiting paths and bounded continuity of $f$ and $\mathcal Af$ permit passage to the limit in the remaining terms.
The limiting process solves the $(\mathcal A,\mathcal D)$-martingale problem with initial law $\nu$.
Well-posedness identifies its law and proves process convergence; this is the usual martingale-problem identification argument~\cite[Chapter~4]{EthierKurtz1986}.
We apply the same reasoning with $\zeta_N=\Pi_N$.
Stationarity passes to the process limit, so $\mu$ is invariant.

Invariance of $\Pi_N$ gives $\zeta_NP_s^N\le C\Pi_N$ for every original time $s$.
Thus \eqref{unif:moment-ui} makes the $p$th powers of the scaled norms uniformly integrable over both $N$ and time under $\zeta_N$.
Weak convergence at each fixed time consequently strengthens to Wasserstein-$p$ convergence.
For fixed $N$, total-variation convergence to $\Pi_N$ and domination give Wasserstein-$p$ convergence of the observed laws to $\mu_N$.
The same uniform-integrability condition also gives $\mu_N\to\mu$ in Wasserstein distance of order $p$.

Passing the initial domination inequality to the limit gives $\nu\le C\mu$.
Invariance of $\mu$ then gives $\nu P_t\le C\mu$ for every limiting-process time $t$.
The probability $\mu$ has a finite $p$th moment by \eqref{unif:moment-ui}.
Therefore the limiting process's total-variation convergence, together with this domination, gives Wasserstein-$p$ convergence to $\mu$.
These four convergences prove \eqref{unif:two-orders}.
\end{proof}

\section{Stationary profiles under JSQ and threshold priority}
\label{app:threshold-profiles}

The threshold estimates in Appendix~\ref{pri:section} determine the coordinates controlled by exact priority.
For a fixed threshold and fallback function $F$, attraction of the profile flow in Theorem~\ref{thm:profile} identifies the remaining coordinates.
These profiles complete the customer-delay arguments in Section~\ref{sec:delay-proof}.

\begin{proof}[Proof of Proposition~\ref{ext:profiles}]
Write $\Omega_j=\min(V_j,H_j)$, so $S=jN-V_j+H_j$ gives $H_j=N(X-j)^++\Omega_j$.
At every threshold $j\ge2$ where exact priority holds, Lemma~\ref{pri:marked} gives
\begin{equation}
 \frac{\E\Omega_j}{N}\le C_\beta N^{-1/3}.
 \label{ext:threshold-overlap}
\end{equation}
The same bound holds at $j=1$ because $\Omega_1\le I$ and $\E I=\beta$.
The constant $C_\beta$ in \eqref{ext:threshold-overlap} is independent of the threshold and the fallback rule.

\emph{JSQ.}
JSQ has the required priority at every fixed $j$.
Equation~\eqref{ext:threshold-overlap} gives $H_j/N-(X_N-j)^+\to0$ in $L^1$.
Since $Z_{m+1}=H_m-H_{m+1}$, this identifies each coordinate in \eqref{ext:jsq-profile}.
It also gives
\[
 \limsup_{N\to\infty}\E[H_j/N]
 =\E(X_*-j)^+\longrightarrow0\quad(j\to\infty),
\]
using the first-moment convergence of $X_N$.
Thus the contribution of high levels to the profile mean vanishes, proving the JSQ conclusion in $\ell^1$.

\emph{Fixed-threshold priority.}
Put $\ell=r+1$, $V=V_\ell$, $H=H_\ell$, and $P=P_{2,\ell}$.
Equation~\eqref{ext:threshold-overlap} applies at every $j\le\ell$.
It identifies the first $r$ profile coordinates and gives
\begin{equation}
 \frac{\min(V,H)}N\longrightarrow0,\qquad
 \frac VN-(\ell-X_N)^+\longrightarrow0,\qquad
 \frac HN-(X_N-\ell)^+\longrightarrow0
 \quad\text{in }L^1.
 \label{ext:priority-mass}
\end{equation}
We next control the tails of the upper profile so that its total mass is preserved in the limit.
Lemma~\ref{pri:marked} gives $B=\E[P\ind_{\{V>0\}}]=O(N)$.
The variables $V$ and $P$ never change at the same jump, and while $V>0$ the drift of $P$ is $-2H+Z_{\ell+1}$.
Stationarity of $VP$ gives
\[
 \E[(N-I-Z_{\ell+1})P]
 =\lambda B+\E[V(2H-Z_{\ell+1})].
\]
Using $I\le N\ind_{\{V>0\}}$, $V\le\ell N$, and $\E H=O(N)$ yields
\begin{equation}
 \E[(N-Z_{\ell+1})P]
 \le(\lambda+N)B+2\ell N\E H=O(N^2).
 \label{ext:priority-local-quadratic}
\end{equation}

The threshold rule satisfies the cutoff conditions of Lemma~\ref{reg:regeneration} with $\ell=r+1$.
Using \eqref{ext:priority-local-quadratic} and $P\le S^2$, the localization argument in the proof of Lemma~\ref{reg:tail-tightness} applies at this cutoff and gives
\[
 \E[(P/N)\ind_{\{X_N\le R\}}]\le C_R.
\]
Here $C_R$ may depend on the fixed threshold $r$ and fallback function $F$.
For $m>\ell$, the bounds $H_m\le P/(m-\ell)$ and $H_m\le S$, together with the common second population moment, give
\[
 \E[H_m/N]\le\frac{C_R}{m-\ell}+\frac CR.
\]
We take the limit superior as $N\to\infty$, then let $m\to\infty$ and finally $R\to\infty$ to obtain ordinary profile-tail tightness.

Let $v_k^N=Z_{\ell+k,N}/N$ for $k\ge1$ and set $v_0^N=1$.
To apply Lemma~\ref{reg:ordinary-profile}, we also verify that the normalized deficit $V/N$ changes negligibly on fixed physical-time intervals.
The process $V/N$ is $C$-tight on fixed physical-time intervals: $0\le V/N\le\ell$, its jumps are $1/N$, its drift is bounded, and its martingale has predictable quadratic variation at most $2t/N$.
By stationarity, \eqref{ext:priority-mass}, and the frozen-population estimate \eqref{reg:frozen-population}, every subsequential limit is the constant path $(\ell-X_*)^+$, jointly with the limiting population $X_*$.
Consequently,
\[
 \sup_{t\le T}\frac{|V(t)-V(0)|}N\longrightarrow0
 \quad\text{in probability},
\]
which verifies \eqref{reg:deficit-increments}.
The ordinary tail bound above and the mean relation \eqref{ext:priority-mass} verify \eqref{reg:profile-inputs}.
Lemma~\ref{reg:ordinary-profile} now identifies every limiting upper profile as $q^F((X_*-\ell)^+)$.
Together with the first $r$ coordinates and the population limit under one-job priority from Example~\ref{ex:priority-nds}, this proves \eqref{ext:priority-profile} and the joint $\ell^1$ convergence.
\end{proof}

\bibliographystyle{abbrvnat}
\bibliography{routing_info}
\end{document}